
\ifx\shlhetal\undefinedcontrolsequence\let\shlhetal\relax\fi
\def\fmtname{AmS-TeX}

\def\fmtversion{2.2}
\catcode`\@=11
\ifx\amstexloaded@\relax\catcode`\@=\active
  \endinput\else\let\amstexloaded@\relax\fi
\newlinechar=`\^^J
\def\W@{\immediate\write\sixt@@n}
\def\CR@{\W@{^^J\fmtname - Version \fmtversion^^J}}
\CR@ \everyjob{\CR@}
\message{Loading definitions for}
\message{misc utility macros,}
\toksdef\toks@@=2
\long\def\rightappend@#1\to#2{\toks@{\\{#1}}\toks@@
 =\expandafter{#2}\xdef#2{\the\toks@@\the\toks@}\toks@{}\toks@@{}}
\def\alloclist@{}
\newif\ifalloc@
\def\showallocations{{\def\\{\immediate\write\m@ne}\alloclist@}\alloc@true}
\def\alloc@#1#2#3#4#5{\global\advance\count1#1by\@ne
 \ch@ck#1#4#2\allocationnumber=\count1#1
 \global#3#5=\allocationnumber
 \edef\next@{\string#5=\string#2\the\allocationnumber}%
 \expandafter\rightappend@\next@\to\alloclist@}
\newcount\count@@
\newcount\count@@@
\def\FN@{\futurelet\next}
\def\DN@{\def\next@}
\def\DNii@{\def\nextii@}
\def\RIfM@{\relax\ifmmode}
\def\RIfMIfI@{\relax\ifmmode\ifinner}
\def\setboxz@h{\setbox\z@\hbox}
\def\wdz@{\wd\z@}
\def\boxz@{\box\z@}
\def\setbox@ne{\setbox\@ne}
\def\wd@ne{\wd\@ne}
\def\iterate{\body\expandafter\iterate\else\fi}
\def\err@#1{\errmessage{AmS-TeX error: #1}}
\newhelp\defaulthelp@{Sorry, I already gave what help I could...^^J
Maybe you should try asking a human?^^J
An error might have occurred before I noticed any problems.^^J
``If all else fails, read the instructions.''}
\def\Err@{\errhelp\defaulthelp@\err@}
\def\eat@#1{}
\def\in@#1#2{\def\in@@##1#1##2##3\in@@{\ifx\in@##2\in@false\else\in@true\fi}%
 \in@@#2#1\in@\in@@}
\newif\ifin@
\def\space@.{\futurelet\space@\relax}
\space@. %
\newhelp\athelp@
{Only certain combinations beginning with @ make sense to me.^^J
Perhaps you wanted \string\@\space for a printed @?^^J
I've ignored the character or group after @.}
{\catcode`\~=\active 
 \lccode`\~=`\@ \lowercase{\gdef~{\FN@\at@}}}
\def\at@{\let\next@\at@@
 \ifcat\noexpand\next a\else\ifcat\noexpand\next0\else
 \ifcat\noexpand\next\relax\else
   \let\next\at@@@\fi\fi\fi
 \next@}
\def\at@@#1{\expandafter
 \ifx\csname\space @\string#1\endcsname\relax
  \expandafter\at@@@ \else
  \csname\space @\string#1\expandafter\endcsname\fi}
\def\at@@@#1{\errhelp\athelp@ \err@{\Invalid@@ @}}
\def\atdef@#1{\expandafter\def\csname\space @\string#1\endcsname}
\newhelp\defahelp@{If you typed \string\define\space cs instead of
\string\define\string\cs\space^^J
I've substituted an inaccessible control sequence so that your^^J
definition will be completed without mixing me up too badly.^^J
If you typed \string\define{\string\cs} the inaccessible control sequence^^J
was defined to be \string\cs, and the rest of your^^J
definition appears as input.}
\newhelp\defbhelp@{I've ignored your definition, because it might^^J
conflict with other uses that are important to me.}
\def\define{\FN@\define@}
\def\define@{\ifcat\noexpand\next\relax
 \expandafter\define@@\else\errhelp\defahelp@                               
 \err@{\string\define\space must be followed by a control
 sequence}\expandafter\def\expandafter\nextii@\fi}                          
\def\undefined@@@@@@@@@@{}
\def\preloaded@@@@@@@@@@{}
\def\next@@@@@@@@@@{}
\def\define@@#1{\ifx#1\relax\errhelp\defbhelp@                              
 \err@{\string#1\space is already defined}\DN@{\DNii@}\else
 \expandafter\ifx\csname\expandafter\eat@\string                            
 #1@@@@@@@@@@\endcsname\undefined@@@@@@@@@@\errhelp\defbhelp@
 \err@{\string#1\space can't be defined}\DN@{\DNii@}\else
 \expandafter\ifx\csname\expandafter\eat@\string#1\endcsname\relax          
 \global\let#1\undefined\DN@{\def#1}\else\errhelp\defbhelp@
 \err@{\string#1\space is already defined}\DN@{\DNii@}\fi
 \fi\fi\next@}

\def\predefine#1#2{\let#1#2}
\def\undefine#1{\let#1\undefined}
\message{page layout,}
\newdimen\captionwidth@
\captionwidth@\hsize
\advance\captionwidth@-1.5in
\def\pagewidth#1{\hsize#1\relax
 \captionwidth@\hsize\advance\captionwidth@-1.5in}
\def\pageheight#1{\vsize#1\relax}
\def\hcorrection#1{\advance\hoffset#1\relax}
\def\vcorrection#1{\advance\voffset#1\relax}
\message{accents/punctuation,}

\let\graveaccent\`
\let\acuteaccent\'
\let\tildeaccent\~
\let\hataccent\^
\let\underscore\_
\let\B\=
\let\D\.
\let\ic@\/
\def\/{\unskip\ic@}
\def\textfonti{\the\textfont\@ne}
\def\t#1#2{{\edef\next@{\the\font}\textfonti\accent"7F \next@#1#2}}
\def~{\unskip\nobreak\ \ignorespaces}
\def\.{.\spacefactor\@m}
\atdef@;{\leavevmode\null;}
\atdef@:{\leavevmode\null:}
\atdef@?{\leavevmode\null?}
\edef\@{\string @}
\def\relaxnext@{\let\next\relax}
\atdef@-{\relaxnext@\leavevmode
 \DN@{\ifx\next-\DN@-{\FN@\nextii@}\else
  \DN@{\leavevmode\hbox{-}}\fi\next@}%
 \DNii@{\ifx\next-\DN@-{\leavevmode\hbox{---}}\else
  \DN@{\leavevmode\hbox{--}}\fi\next@}%
 \FN@\next@}
\def\srdr@{\kern.16667em}
\def\drsr@{\kern.02778em}
\def\sldl@{\drsr@}
\def\dlsl@{\srdr@}
\atdef@"{\unskip\relaxnext@
 \DN@{\ifx\next\space@\DN@. {\FN@\nextii@}\else
  \DN@.{\FN@\nextii@}\fi\next@.}%
 \DNii@{\ifx\next`\DN@`{\FN@\nextiii@}\else
  \ifx\next\lq\DN@\lq{\FN@\nextiii@}\else
  \DN@####1{\FN@\nextiv@}\fi\fi\next@}%
 \def\nextiii@{\ifx\next`\DN@`{\sldl@``}\else\ifx\next\lq
  \DN@\lq{\sldl@``}\else\DN@{\dlsl@`}\fi\fi\next@}%
 \def\nextiv@{\ifx\next'\DN@'{\srdr@''}\else
  \ifx\next\rq\DN@\rq{\srdr@''}\else\DN@{\drsr@'}\fi\fi\next@}%
 \FN@\next@}

\def\textfontii{\the\textfont\tw@}
\def\lbrace@{\delimiter"4266308 }
\def\rbrace@{\delimiter"5267309 }
\def\{{\RIfM@\lbrace@\else{\textfontii f}\spacefactor\@m\fi}
\def\}{\RIfM@\rbrace@\else
 \let\@sf\empty\ifhmode\edef\@sf{\spacefactor\the\spacefactor}\fi
 {\textfontii g}\@sf\relax\fi}
\let\lbrace\{
\let\rbrace\}
\def\AmSTeX{{\textfontii A\kern-.1667em%
  \lower.5ex\hbox{M}\kern-.125emS}-\TeX\spacefactor1000 }
\message{line and page breaks,}
\def\vmodeerr@#1{\Err@{\string#1\space not allowed between paragraphs}}
\def\mathmodeerr@#1{\Err@{\string#1\space not allowed in math mode}}
\def\linebreak{\RIfM@\mathmodeerr@\linebreak\else
 \ifhmode\unskip\unkern\break\else\vmodeerr@\linebreak\fi\fi}

\newskip\saveskip@
\def\allowlinebreak{\RIfM@\mathmodeerr@\allowlinebreak\else
 \ifhmode\saveskip@\lastskip\unskip
 \allowbreak\ifdim\saveskip@>\z@\hskip\saveskip@\fi
 \else\vmodeerr@\allowlinebreak\fi\fi}
\def\nolinebreak{\RIfM@\mathmodeerr@\nolinebreak\else
 \ifhmode\saveskip@\lastskip\unskip
 \nobreak\ifdim\saveskip@>\z@\hskip\saveskip@\fi
 \else\vmodeerr@\nolinebreak\fi\fi}
\def\newline{\relaxnext@
 \DN@{\RIfM@\expandafter\mathmodeerr@\expandafter\newline\else
  \ifhmode\ifx\next\par\else
  \expandafter\unskip\expandafter\null\expandafter\hfill\expandafter\break\fi
  \else
  \expandafter\vmodeerr@\expandafter\newline\fi\fi}%
 \FN@\next@}
\def\dmatherr@#1{\Err@{\string#1\space not allowed in display math mode}}
\def\nondmatherr@#1{\Err@{\string#1\space not allowed in non-display math
 mode}}
\def\onlydmatherr@#1{\Err@{\string#1\space allowed only in display math mode}}
\def\nonmatherr@#1{\Err@{\string#1\space allowed only in math mode}}
\def\mathbreak{\RIfMIfI@\break\else
 \dmatherr@\mathbreak\fi\else\nonmatherr@\mathbreak\fi}
\def\nomathbreak{\RIfMIfI@\nobreak\else
 \dmatherr@\nomathbreak\fi\else\nonmatherr@\nomathbreak\fi}
\def\allowmathbreak{\RIfMIfI@\allowbreak\else
 \dmatherr@\allowmathbreak\fi\else\nonmatherr@\allowmathbreak\fi}
\def\pagebreak{\RIfM@
 \ifinner\nondmatherr@\pagebreak\else\postdisplaypenalty-\@M\fi
 \else\ifvmode\removelastskip\break\else\vadjust{\break}\fi\fi}
\def\nopagebreak{\RIfM@
 \ifinner\nondmatherr@\nopagebreak\else\postdisplaypenalty\@M\fi
 \else\ifvmode\nobreak\else\vadjust{\nobreak}\fi\fi}
\def\nonvmodeerr@#1{\Err@{\string#1\space not allowed within a paragraph
 or in math}}
\def\vnonvmode@#1#2{\relaxnext@\DNii@{\ifx\next\par\DN@{#1}\else
 \DN@{#2}\fi\next@}%
 \ifvmode\DN@{#1}\else
 \DN@{\FN@\nextii@}\fi\next@}
\def\newpage{\vnonvmode@{\vfill\break}{\nonvmodeerr@\newpage}}
\def\smallpagebreak{\vnonvmode@\smallbreak{\nonvmodeerr@\smallpagebreak}}
\def\medpagebreak{\vnonvmode@\medbreak{\nonvmodeerr@\medpagebreak}}
\def\bigpagebreak{\vnonvmode@\bigbreak{\nonvmodeerr@\bigpagebreak}}
\def\NoBlackBoxes{\global\overfullrule\z@}
\def\BlackBoxes{\global\overfullrule5\p@}
\def\Invalid@#1{\def#1{\Err@{\Invalid@@\string#1}}}
\def\Invalid@@{Invalid use of }
\message{figures,}
\Invalid@\caption
\Invalid@\captionwidth
\newdimen\smallcaptionwidth@
\def\topspace{\mid@false\ins@}
\def\midspace{\mid@true\ins@}
\newif\ifmid@
\def\captionfont@{}
\def\ins@#1{\relaxnext@\allowbreak
 \smallcaptionwidth@\captionwidth@\gdef\thespace@{#1}%
 \DN@{\ifx\next\space@\DN@. {\FN@\nextii@}\else
  \DN@.{\FN@\nextii@}\fi\next@.}%
 \DNii@{\ifx\next\caption\DN@\caption{\FN@\nextiii@}%
  \else\let\next@\nextiv@\fi\next@}%
 \def\nextiv@{\vnonvmode@
  {\ifmid@\expandafter\midinsert\else\expandafter\topinsert\fi
   \vbox to\thespace@{}\endinsert}
  {\ifmid@\nonvmodeerr@\midspace\else\nonvmodeerr@\topspace\fi}}%
 \def\nextiii@{\ifx\next\captionwidth\expandafter\nextv@
  \else\expandafter\nextvi@\fi}%
 \def\nextv@\captionwidth##1##2{\smallcaptionwidth@##1\relax\nextvi@{##2}}%
 \def\nextvi@##1{\def\thecaption@{\captionfont@##1}%
  \DN@{\ifx\next\space@\DN@. {\FN@\nextvii@}\else
   \DN@.{\FN@\nextvii@}\fi\next@.}%
  \FN@\next@}%
 \def\nextvii@{\vnonvmode@
  {\ifmid@\expandafter\midinsert\else
  \expandafter\topinsert\fi\vbox to\thespace@{}\nobreak\smallskip
  \setboxz@h{\noindent\ignorespaces\thecaption@\unskip}%
  \ifdim\wdz@>\smallcaptionwidth@\centerline{\vbox{\hsize\smallcaptionwidth@
   \noindent\ignorespaces\thecaption@\unskip}}%
  \else\centerline{\boxz@}\fi\endinsert}
  {\ifmid@\nonvmodeerr@\midspace
  \else\nonvmodeerr@\topspace\fi}}%
 \FN@\next@}
\message{comments,}
\def\newcodes@{\catcode`\\12\catcode`\{12\catcode`\}12\catcode`\#12%
 \catcode`\%12\relax}
\def\oldcodes@{\catcode`\\0\catcode`\{1\catcode`\}2\catcode`\#6%
 \catcode`\%14\relax}
\def\comment{\newcodes@\endlinechar=10 \comment@}
{\lccode`\0=`\\
\lowercase{\gdef\comment@#1^^J{\comment@@#10endcomment\comment@@@}%
\gdef\comment@@#10endcomment{\FN@\comment@@@}%
\gdef\comment@@@#1\comment@@@{\ifx\next\comment@@@\let\next\comment@
 \else\def\next{\oldcodes@\endlinechar=`\^^M\relax}%
 \fi\next}}}
\def\pr@m@s{\ifx'\next\DN@##1{\prim@s}\else\let\next@\egroup\fi\next@}
\def\prime{{\null\prime@\null}}
\mathchardef\prime@="0230
\let\dsize\displaystyle

\let\ssize\scriptstyle

\message{math spacing,}
\def\,{\RIfM@\mskip\thinmuskip\relax\else\kern.16667em\fi}
\def\!{\RIfM@\mskip-\thinmuskip\relax\else\kern-.16667em\fi}
\let\thinspace\,
\let\negthinspace\!
\def\medspace{\RIfM@\mskip\medmuskip\relax\else\kern.222222em\fi}
\def\negmedspace{\RIfM@\mskip-\medmuskip\relax\else\kern-.222222em\fi}
\def\thickspace{\RIfM@\mskip\thickmuskip\relax\else\kern.27777em\fi}
\let\;\thickspace
\def\negthickspace{\RIfM@\mskip-\thickmuskip\relax\else
 \kern-.27777em\fi}
\atdef@,{\RIfM@\mskip.1\thinmuskip\else\leavevmode\null,\fi}
\atdef@!{\RIfM@\mskip-.1\thinmuskip\else\leavevmode\null!\fi}
\atdef@.{\RIfM@&&\else\leavevmode.\spacefactor3000 \fi}
\def\and{\DOTSB\;\mathchar"3026 \;}

\message{fractions,}
\def\frac#1#2{{#1\over#2}}

\newdimen\ex@
\ex@.2326ex
\Invalid@\thickness
\def\thickfrac{\relaxnext@
 \DN@{\ifx\next\thickness\let\next@\nextii@\else
 \DN@{\nextii@\thickness1}\fi\next@}%
 \DNii@\thickness##1##2##3{{##2\above##1\ex@##3}}%
 \FN@\next@}

\def\thickfracwithdelims#1#2{\relaxnext@\def\ldelim@{#1}\def\rdelim@{#2}%
 \DN@{\ifx\next\thickness\let\next@\nextii@\else
 \DN@{\nextii@\thickness1}\fi\next@}%
 \DNii@\thickness##1##2##3{{##2\abovewithdelims
 \ldelim@\rdelim@##1\ex@##3}}%
 \FN@\next@}

\def\:{\nobreak\hskip.1111em\mathpunct{}\nonscript\mkern-\thinmuskip{:}\hskip
 .3333emplus.0555em\relax}
\def\snug{\unskip\kern-\mathsurround}
\message{smash commands,}
\def\topsmash{\top@true\bot@false\smash@}
\def\botsmash{\top@false\bot@true\smash@}
\newif\iftop@
\newif\ifbot@
\def\smash{\top@true\bot@true\smash@}
\def\smash@{\RIfM@\expandafter\mathpalette\expandafter\mathsm@sh\else
 \expandafter\makesm@sh\fi}
\def\finsm@sh{\iftop@\ht\z@\z@\fi\ifbot@\dp\z@\z@\fi\leavevmode\boxz@}
\message{large operator symbols,}
\def\LimitsOnSums{\global\let\slimits@\displaylimits}
\def\NoLimitsOnSums{\global\let\slimits@\nolimits}
\LimitsOnSums
\mathchardef\coprod@="1360       \def\coprod{\DOTSB\coprod@\slimits@}
\mathchardef\bigvee@="1357       \def\bigvee{\DOTSB\bigvee@\slimits@}
\mathchardef\bigwedge@="1356     \def\bigwedge{\DOTSB\bigwedge@\slimits@}
\mathchardef\biguplus@="1355     \def\biguplus{\DOTSB\biguplus@\slimits@}
\mathchardef\bigcap@="1354       \def\bigcap{\DOTSB\bigcap@\slimits@}
\mathchardef\bigcup@="1353       \def\bigcup{\DOTSB\bigcup@\slimits@}
\mathchardef\prod@="1351         \def\prod{\DOTSB\prod@\slimits@}
\mathchardef\sum@="1350          \def\sum{\DOTSB\sum@\slimits@}
\mathchardef\bigotimes@="134E    \def\bigotimes{\DOTSB\bigotimes@\slimits@}
\mathchardef\bigoplus@="134C     \def\bigoplus{\DOTSB\bigoplus@\slimits@}
\mathchardef\bigodot@="134A      \def\bigodot{\DOTSB\bigodot@\slimits@}
\mathchardef\bigsqcup@="1346     \def\bigsqcup{\DOTSB\bigsqcup@\slimits@}
\message{integrals,}
\def\LimitsOnInts{\global\let\ilimits@\displaylimits}
\def\NoLimitsOnInts{\global\let\ilimits@\nolimits}
\NoLimitsOnInts
\def\int{\DOTSI\intop\ilimits@}
\def\oint{\DOTSI\ointop\ilimits@}
\def\intic@{\mathchoice{\hskip.5em}{\hskip.4em}{\hskip.4em}{\hskip.4em}}
\def\negintic@{\mathchoice
 {\hskip-.5em}{\hskip-.4em}{\hskip-.4em}{\hskip-.4em}}
\def\intkern@{\mathchoice{\!\!\!}{\!\!}{\!\!}{\!\!}}
\def\intdots@{\mathchoice{\plaincdots@}
 {{\cdotp}\mkern1.5mu{\cdotp}\mkern1.5mu{\cdotp}}
 {{\cdotp}\mkern1mu{\cdotp}\mkern1mu{\cdotp}}
 {{\cdotp}\mkern1mu{\cdotp}\mkern1mu{\cdotp}}}
\newcount\intno@
\def\iint{\DOTSI\intno@\tw@\FN@\ints@}
\def\iiint{\DOTSI\intno@\thr@@\FN@\ints@}
\def\iiiint{\DOTSI\intno@4 \FN@\ints@}
\def\idotsint{\DOTSI\intno@\z@\FN@\ints@}
\def\ints@{\findlimits@\ints@@}
\newif\iflimtoken@
\newif\iflimits@
\def\findlimits@{\limtoken@true\ifx\next\limits\limits@true
 \else\ifx\next\nolimits\limits@false\else
 \limtoken@false\ifx\ilimits@\nolimits\limits@false\else
 \ifinner\limits@false\else\limits@true\fi\fi\fi\fi}
\def\multint@{\int\ifnum\intno@=\z@\intdots@                                
 \else\intkern@\fi                                                          
 \ifnum\intno@>\tw@\int\intkern@\fi                                         
 \ifnum\intno@>\thr@@\int\intkern@\fi                                       
 \int}                                                                      
\def\multintlimits@{\intop\ifnum\intno@=\z@\intdots@\else\intkern@\fi
 \ifnum\intno@>\tw@\intop\intkern@\fi
 \ifnum\intno@>\thr@@\intop\intkern@\fi\intop}
\def\ints@@{\iflimtoken@                                                    
 \def\ints@@@{\iflimits@\negintic@\mathop{\intic@\multintlimits@}\limits    
  \else\multint@\nolimits\fi                                                
  \eat@}                                                                    
 \else                                                                      
 \def\ints@@@{\iflimits@\negintic@
  \mathop{\intic@\multintlimits@}\limits\else
  \multint@\nolimits\fi}\fi\ints@@@}
\def\LimitsOnNames{\global\let\nlimits@\displaylimits}
\def\NoLimitsOnNames{\global\let\nlimits@\nolimits@}
\LimitsOnNames
\def\nolimits@{\relaxnext@
 \DN@{\ifx\next\limits\DN@\limits{\nolimits}\else
  \let\next@\nolimits\fi\next@}%
 \FN@\next@}
\message{operator names,}
\def\newmcodes@{\mathcode`\'"27\mathcode`\*"2A\mathcode`\."613A%
 \mathcode`\-"2D\mathcode`\/"2F\mathcode`\:"603A }
\def\operatorname#1{\mathop{\newmcodes@\kern\z@\fam\z@#1}\nolimits@}
\def\operatornamewithlimits#1{\mathop{\newmcodes@\kern\z@\fam\z@#1}\nlimits@}
\def\qopname@#1{\mathop{\fam\z@#1}\nolimits@}
\def\qopnamewl@#1{\mathop{\fam\z@#1}\nlimits@}
\def\arccos{\qopname@{arccos}}
\def\arcsin{\qopname@{arcsin}}
\def\arctan{\qopname@{arctan}}
\def\arg{\qopname@{arg}}
\def\cos{\qopname@{cos}}
\def\cosh{\qopname@{cosh}}
\def\cot{\qopname@{cot}}
\def\coth{\qopname@{coth}}
\def\csc{\qopname@{csc}}
\def\deg{\qopname@{deg}}
\def\det{\qopnamewl@{det}}
\def\dim{\qopname@{dim}}
\def\exp{\qopname@{exp}}
\def\gcd{\qopnamewl@{gcd}}
\def\hom{\qopname@{hom}}
\def\inf{\qopnamewl@{inf}}
\def\injlim{\qopnamewl@{inj\,lim}}
\def\ker{\qopname@{ker}}
\def\lg{\qopname@{lg}}
\def\lim{\qopnamewl@{lim}}
\def\liminf{\qopnamewl@{lim\,inf}}
\def\limsup{\qopnamewl@{lim\,sup}}
\def\ln{\qopname@{ln}}
\def\log{\qopname@{log}}
\def\max{\qopnamewl@{max}}
\def\min{\qopnamewl@{min}}
\def\Pr{\qopnamewl@{Pr}}
\def\projlim{\qopnamewl@{proj\,lim}}
\def\sec{\qopname@{sec}}
\def\sin{\qopname@{sin}}
\def\sinh{\qopname@{sinh}}
\def\sup{\qopnamewl@{sup}}
\def\tan{\qopname@{tan}}
\def\tanh{\qopname@{tanh}}
\def\varinjlim{\mathop{\vtop{\ialign{##\crcr
 \hfil\rm lim\hfil\crcr\noalign{\nointerlineskip}\rightarrowfill\crcr
 \noalign{\nointerlineskip\kern-\ex@}\crcr}}}}
\def\varprojlim{\mathop{\vtop{\ialign{##\crcr
 \hfil\rm lim\hfil\crcr\noalign{\nointerlineskip}\leftarrowfill\crcr
 \noalign{\nointerlineskip\kern-\ex@}\crcr}}}}
\def\varliminf{\mathop{\underline{\vrule height\z@ depth.2exwidth\z@
 \hbox{\rm lim}}}}

\newdimen\buffer@
\buffer@\fontdimen13 \tenex
\newdimen\buffer
\buffer\buffer@

\def\ResetBuffer{\fontdimen13 \tenex\buffer@\global\buffer\buffer@}
\def\shave#1{\mathop{\hbox{$\m@th\fontdimen13 \tenex\z@                     
 \displaystyle{#1}$}}\fontdimen13 \tenex\buffer}

\message{multilevel sub/superscripts,}
\Invalid@\\
\def\Let@{\relax\iffalse{\fi\let\\=\cr\iffalse}\fi}
\Invalid@\vspace
\def\vspace@{\def\vspace##1{\crcr\noalign{\vskip##1\relax}}}
\def\multilimits@{\bgroup\vspace@\Let@
 \baselineskip\fontdimen10 \scriptfont\tw@
 \advance\baselineskip\fontdimen12 \scriptfont\tw@
 \lineskip\thr@@\fontdimen8 \scriptfont\thr@@
 \lineskiplimit\lineskip
 \vbox\bgroup\ialign\bgroup\hfil$\m@th\scriptstyle{##}$\hfil\crcr}
\def\Sb{_\multilimits@}
\def\endSb{\crcr\egroup\egroup\egroup}
\def\Sp{^\multilimits@}

\def\spreadlines#1{\RIfMIfI@\onlydmatherr@\spreadlines\else
 \openup#1\relax\fi\else\onlydmatherr@\spreadlines\fi}
\def\Mathstrut@{\copy\Mathstrutbox@}
\newbox\Mathstrutbox@
\setbox\Mathstrutbox@\null
\setboxz@h{$\m@th($}
\ht\Mathstrutbox@\ht\z@
\dp\Mathstrutbox@\dp\z@
\message{matrices,}
\newdimen\spreadmlines@
\def\spreadmatrixlines#1{\RIfMIfI@
 \onlydmatherr@\spreadmatrixlines\else
 \spreadmlines@#1\relax\fi\else\onlydmatherr@\spreadmatrixlines\fi}
\def\matrix{\null\,\vcenter\bgroup\Let@\vspace@
 \normalbaselines\openup\spreadmlines@\ialign
 \bgroup\hfil$\m@th##$\hfil&&\quad\hfil$\m@th##$\hfil\crcr
 \Mathstrut@\crcr\noalign{\kern-\baselineskip}}
\def\endmatrix{\crcr\Mathstrut@\crcr\noalign{\kern-\baselineskip}\egroup
 \egroup\,}
\def\format{\crcr\egroup\iffalse{\fi\ifnum`}=0 \fi\format@}
\newtoks\hashtoks@
\hashtoks@{#}
\def\format@#1\\{\def\preamble@{#1}%
 \def\l{$\m@th\the\hashtoks@$\hfil}%
 \def\c{\hfil$\m@th\the\hashtoks@$\hfil}%
 \def\r{\hfil$\m@th\the\hashtoks@$}%
 \edef\preamble@@{\preamble@}\ifnum`{=0 \fi\iffalse}\fi
 \ialign\bgroup\span\preamble@@\crcr}
\def\smallmatrix{\null\,\vcenter\bgroup\vspace@\Let@
 \baselineskip9\ex@\lineskip\ex@
 \ialign\bgroup\hfil$\m@th\scriptstyle{##}$\hfil&&\thickspace\hfil
 $\m@th\scriptstyle{##}$\hfil\crcr}
\def\endsmallmatrix{\crcr\egroup\egroup\,}

\newmuskip\dotsspace@
\dotsspace@1.5mu
\def\strip@#1 {#1}
\def\spacehdots#1\for#2{\multispan{#2}\xleaders
 \hbox{$\m@th\mkern\strip@#1 \dotsspace@.\mkern\strip@#1 \dotsspace@$}\hfill}
\def\hdotsfor#1{\spacehdots\@ne\for{#1}}
\def\multispan@#1{\omit\mscount#1\unskip\loop\ifnum\mscount>\@ne\sp@n\repeat}
\def\spaceinnerhdots#1\for#2\after#3{\multispan@{\strip@#2 }#3\xleaders
 \hbox{$\m@th\mkern\strip@#1 \dotsspace@.\mkern\strip@#1 \dotsspace@$}\hfill}
\def\innerhdotsfor#1\after#2{\spaceinnerhdots\@ne\for#1\after{#2}}
\def\cases{\bgroup\spreadmlines@\jot\left\{\,\matrix\format\l&\quad\l\\}
\def\endcases{\endmatrix\right.\egroup}
\message{multiline displays,}
\newif\ifinany@
\newif\ifinalign@
\newif\ifingather@
\def\strut@{\copy\strutbox@}
\newbox\strutbox@
\setbox\strutbox@\hbox{\vrule height8\p@ depth3\p@ width\z@}
\def\topaligned{\null\,\vtop\aligned@}
\def\botaligned{\null\,\vbox\aligned@}
\def\aligned{\null\,\vcenter\aligned@}
\def\aligned@{\bgroup\vspace@\Let@
 \ifinany@\else\openup\jot\fi\ialign
 \bgroup\hfil\strut@$\m@th\displaystyle{##}$&
 $\m@th\displaystyle{{}##}$\hfil\crcr}
\def\endaligned{\crcr\egroup\egroup}

\def\alignedat#1{\null\,\vcenter\bgroup\doat@{#1}\vspace@\Let@
 \ifinany@\else\openup\jot\fi\ialign\bgroup\span\preamble@@\crcr}
\newcount\atcount@
\def\doat@#1{\toks@{\hfil\strut@$\m@th
 \displaystyle{\the\hashtoks@}$&$\m@th\displaystyle
 {{}\the\hashtoks@}$\hfil}
 \atcount@#1\relax\advance\atcount@\m@ne                                    
 \loop\ifnum\atcount@>\z@\toks@=\expandafter{\the\toks@&\hfil$\m@th
 \displaystyle{\the\hashtoks@}$&$\m@th
 \displaystyle{{}\the\hashtoks@}$\hfil}\advance
  \atcount@\m@ne\repeat                                                     
 \xdef\preamble@{\the\toks@}\xdef\preamble@@{\preamble@}}

\def\gathered{\null\,\vcenter\bgroup\vspace@\Let@
 \ifinany@\else\openup\jot\fi\ialign
 \bgroup\hfil\strut@$\m@th\displaystyle{##}$\hfil\crcr}
\def\endgathered{\crcr\egroup\egroup}
\newif\iftagsleft@
\def\TagsOnLeft{\global\tagsleft@true}
\def\TagsOnRight{\global\tagsleft@false}
\TagsOnLeft
\newif\ifmathtags@
\def\TagsAsMath{\global\mathtags@true}
\def\TagsAsText{\global\mathtags@false}
\TagsAsText
\def\tagform@#1{\hbox{\rm(\ignorespaces#1\unskip)}}
\def\thetag{\leavevmode\tagform@}
\def\tag#1$${\iftagsleft@\leqno\else\eqno\fi                                
 \maketag@#1\maketag@                                                       
 $$}                                                                        
\def\maketag@{\FN@\maketag@@}
\def\maketag@@{\ifx\next"\expandafter\maketag@@@\else\expandafter\maketag@@@@
 \fi}
\def\maketag@@@"#1"#2\maketag@{\hbox{\rm#1}}                                
\def\maketag@@@@#1\maketag@{\ifmathtags@\tagform@{$\m@th#1$}\else
 \tagform@{#1}\fi}
\interdisplaylinepenalty\@M
\def\allowdisplaybreaks{\RIfMIfI@
 \onlydmatherr@\allowdisplaybreaks\else
 \interdisplaylinepenalty\z@\fi\else\onlydmatherr@\allowdisplaybreaks\fi}
\Invalid@\allowdisplaybreak
\Invalid@\displaybreak
\Invalid@\intertext
\def\allowdisplaybreak@{\def\allowdisplaybreak{\crcr\noalign{\allowbreak}}}
\def\displaybreak@{\def\displaybreak{\crcr\noalign{\break}}}
\def\intertext@{\def\intertext##1{\crcr\noalign{%
 \penalty\postdisplaypenalty \vskip\belowdisplayskip
 \vbox{\normalbaselines\noindent##1}%
 \penalty\predisplaypenalty \vskip\abovedisplayskip}}}
\newskip\centering@
\centering@\z@ plus\@m\p@
\def\align{\relax\ifingather@\DN@{\csname align (in
  \string\gather)\endcsname}\else
 \ifmmode\ifinner\DN@{\onlydmatherr@\align}\else
  \let\next@\align@\fi
 \else\DN@{\onlydmatherr@\align}\fi\fi\next@}
\newhelp\andhelp@
{An extra & here is so disastrous that you should probably exit^^J
and fix things up.}
\newif\iftag@
\newcount\and@
\def\align@{\inalign@true\inany@true
 \vspace@\allowdisplaybreak@\displaybreak@\intertext@
 \def\tag{\global\tag@true\ifnum\and@=\z@\DN@{&&}\else
          \DN@{&}\fi\next@}%
 \iftagsleft@\DN@{\csname align \endcsname}\else
  \DN@{\csname align \space\endcsname}\fi\next@}
\def\Tag@{\iftag@\else\errhelp\andhelp@\err@{Extra & on this line}\fi}
\newdimen\lwidth@
\newdimen\rwidth@
\newdimen\maxlwidth@
\newdimen\maxrwidth@
\newdimen\totwidth@
\def\measure@#1\endalign{\lwidth@\z@\rwidth@\z@\maxlwidth@\z@\maxrwidth@\z@
 \global\and@\z@                                                            
 \setbox@ne\vbox                                                            
  {\everycr{\noalign{\global\tag@false\global\and@\z@}}\Let@                
  \halign{\setboxz@h{$\m@th\displaystyle{\@lign##}$}
   \global\lwidth@\wdz@                                                     
   \ifdim\lwidth@>\maxlwidth@\global\maxlwidth@\lwidth@\fi                  
   \global\advance\and@\@ne                                                 
   &\setboxz@h{$\m@th\displaystyle{{}\@lign##}$}\global\rwidth@\wdz@        
   \ifdim\rwidth@>\maxrwidth@\global\maxrwidth@\rwidth@\fi                  
   \global\advance\and@\@ne                                                
   &\Tag@
   \eat@{##}\crcr#1\crcr}}
 \totwidth@\maxlwidth@\advance\totwidth@\maxrwidth@}                       
\def\displ@y@{\global\dt@ptrue\openup\jot
 \everycr{\noalign{\global\tag@false\global\and@\z@\ifdt@p\global\dt@pfalse
 \vskip-\lineskiplimit\vskip\normallineskiplimit\else
 \penalty\interdisplaylinepenalty\fi}}}
\def\black@#1{\noalign{\ifdim#1>\displaywidth
 \dimen@\prevdepth\nointerlineskip                                          
 \vskip-\ht\strutbox@\vskip-\dp\strutbox@                                   
 \vbox{\noindent\hbox to#1{\strut@\hfill}}
 \prevdepth\dimen@                                                          
 \fi}}
\expandafter\def\csname align \space\endcsname#1\endalign
 {\measure@#1\endalign\global\and@\z@                                       
 \ifingather@\everycr{\noalign{\global\and@\z@}}\else\displ@y@\fi           
 \Let@\tabskip\centering@                                                   
 \halign to\displaywidth
  {\hfil\strut@\setboxz@h{$\m@th\displaystyle{\@lign##}$}
  \global\lwidth@\wdz@\boxz@\global\advance\and@\@ne                        
  \tabskip\z@skip                                                           
  &\setboxz@h{$\m@th\displaystyle{{}\@lign##}$}
  \global\rwidth@\wdz@\boxz@\hfill\global\advance\and@\@ne                  
  \tabskip\centering@                                                       
  &\setboxz@h{\@lign\strut@\maketag@##\maketag@}
  \dimen@\displaywidth\advance\dimen@-\totwidth@
  \divide\dimen@\tw@\advance\dimen@\maxrwidth@\advance\dimen@-\rwidth@     
  \ifdim\dimen@<\tw@\wdz@\llap{\vtop{\normalbaselines\null\boxz@}}
  \else\llap{\boxz@}\fi                                                    
  \tabskip\z@skip                                                          
  \crcr#1\crcr                                                             
  \black@\totwidth@}}                                                      
\newdimen\lineht@
\expandafter\def\csname align \endcsname#1\endalign{\measure@#1\endalign
 \global\and@\z@
 \ifdim\totwidth@>\displaywidth\let\displaywidth@\totwidth@\else
  \let\displaywidth@\displaywidth\fi                                        
 \ifingather@\everycr{\noalign{\global\and@\z@}}\else\displ@y@\fi
 \Let@\tabskip\centering@\halign to\displaywidth
  {\hfil\strut@\setboxz@h{$\m@th\displaystyle{\@lign##}$}%
  \global\lwidth@\wdz@\global\lineht@\ht\z@                                 
  \boxz@\global\advance\and@\@ne
  \tabskip\z@skip&\setboxz@h{$\m@th\displaystyle{{}\@lign##}$}%
  \global\rwidth@\wdz@\ifdim\ht\z@>\lineht@\global\lineht@\ht\z@\fi         
  \boxz@\hfil\global\advance\and@\@ne
  \tabskip\centering@&\kern-\displaywidth@                                  
  \setboxz@h{\@lign\strut@\maketag@##\maketag@}%
  \dimen@\displaywidth\advance\dimen@-\totwidth@
  \divide\dimen@\tw@\advance\dimen@\maxlwidth@\advance\dimen@-\lwidth@
  \ifdim\dimen@<\tw@\wdz@
   \rlap{\vbox{\normalbaselines\boxz@\vbox to\lineht@{}}}\else
   \rlap{\boxz@}\fi
  \tabskip\displaywidth@\crcr#1\crcr\black@\totwidth@}}
\expandafter\def\csname align (in \string\gather)\endcsname
  #1\endalign{\vcenter{\align@#1\endalign}}
\Invalid@\endalign
\newif\ifxat@
\def\alignat{\RIfMIfI@\DN@{\onlydmatherr@\alignat}\else
 \DN@{\csname alignat \endcsname}\fi\else
 \DN@{\onlydmatherr@\alignat}\fi\next@}
\newif\ifmeasuring@
\newbox\savealignat@
\expandafter\def\csname alignat \endcsname#1#2\endalignat                   
 {\inany@true\xat@false
 \def\tag{\global\tag@true\count@#1\relax\multiply\count@\tw@
  \xdef\tag@{}\loop\ifnum\count@>\and@\xdef\tag@{&\tag@}\advance\count@\m@ne
  \repeat\tag@}%
 \vspace@\allowdisplaybreak@\displaybreak@\intertext@
 \displ@y@\measuring@true                                                   
 \setbox\savealignat@\hbox{$\m@th\displaystyle\Let@
  \attag@{#1}
  \vbox{\halign{\span\preamble@@\crcr#2\crcr}}$}%
 \measuring@false                                                           
 \Let@\attag@{#1}
 \tabskip\centering@\halign to\displaywidth
  {\span\preamble@@\crcr#2\crcr                                             
  \black@{\wd\savealignat@}}}                                               
\Invalid@\endalignat
\def\xalignat{\RIfMIfI@
 \DN@{\onlydmatherr@\xalignat}\else
 \DN@{\csname xalignat \endcsname}\fi\else
 \DN@{\onlydmatherr@\xalignat}\fi\next@}
\expandafter\def\csname xalignat \endcsname#1#2\endxalignat
 {\inany@true\xat@true
 \def\tag{\global\tag@true\def\tag@{}\count@#1\relax\multiply\count@\tw@
  \loop\ifnum\count@>\and@\xdef\tag@{&\tag@}\advance\count@\m@ne\repeat\tag@}%
 \vspace@\allowdisplaybreak@\displaybreak@\intertext@
 \displ@y@\measuring@true\setbox\savealignat@\hbox{$\m@th\displaystyle\Let@
 \attag@{#1}\vbox{\halign{\span\preamble@@\crcr#2\crcr}}$}%
 \measuring@false\Let@
 \attag@{#1}\tabskip\centering@\halign to\displaywidth
 {\span\preamble@@\crcr#2\crcr\black@{\wd\savealignat@}}}
\def\attag@#1{\let\Maketag@\maketag@\let\TAG@\Tag@                          
 \let\Tag@=0\let\maketag@=0
 \ifmeasuring@\def\llap@##1{\setboxz@h{##1}\hbox to\tw@\wdz@{}}%
  \def\rlap@##1{\setboxz@h{##1}\hbox to\tw@\wdz@{}}\else
  \let\llap@\llap\let\rlap@\rlap\fi                                         
 \toks@{\hfil\strut@$\m@th\displaystyle{\@lign\the\hashtoks@}$\tabskip\z@skip
  \global\advance\and@\@ne&$\m@th\displaystyle{{}\@lign\the\hashtoks@}$\hfil
  \ifxat@\tabskip\centering@\fi\global\advance\and@\@ne}
 \iftagsleft@
  \toks@@{\tabskip\centering@&\Tag@\kern-\displaywidth
   \rlap@{\@lign\maketag@\the\hashtoks@\maketag@}%
   \global\advance\and@\@ne\tabskip\displaywidth}\else
  \toks@@{\tabskip\centering@&\Tag@\llap@{\@lign\maketag@
   \the\hashtoks@\maketag@}\global\advance\and@\@ne\tabskip\z@skip}\fi      
 \atcount@#1\relax\advance\atcount@\m@ne
 \loop\ifnum\atcount@>\z@
 \toks@=\expandafter{\the\toks@&\hfil$\m@th\displaystyle{\@lign
  \the\hashtoks@}$\global\advance\and@\@ne
  \tabskip\z@skip&$\m@th\displaystyle{{}\@lign\the\hashtoks@}$\hfil\ifxat@
  \tabskip\centering@\fi\global\advance\and@\@ne}\advance\atcount@\m@ne
 \repeat                                                                    
 \xdef\preamble@{\the\toks@\the\toks@@}
 \xdef\preamble@@{\preamble@}
 \let\maketag@\Maketag@\let\Tag@\TAG@}                                      
\Invalid@\endxalignat
\def\xxalignat{\RIfMIfI@
 \DN@{\onlydmatherr@\xxalignat}\else\DN@{\csname xxalignat
  \endcsname}\fi\else
 \DN@{\onlydmatherr@\xxalignat}\fi\next@}
\expandafter\def\csname xxalignat \endcsname#1#2\endxxalignat{\inany@true
 \vspace@\allowdisplaybreak@\displaybreak@\intertext@
 \displ@y\setbox\savealignat@\hbox{$\m@th\displaystyle\Let@
 \xxattag@{#1}\vbox{\halign{\span\preamble@@\crcr#2\crcr}}$}%
 \Let@\xxattag@{#1}\tabskip\z@skip\halign to\displaywidth
 {\span\preamble@@\crcr#2\crcr\black@{\wd\savealignat@}}}
\def\xxattag@#1{\toks@{\tabskip\z@skip\hfil\strut@
 $\m@th\displaystyle{\the\hashtoks@}$&%
 $\m@th\displaystyle{{}\the\hashtoks@}$\hfil\tabskip\centering@&}%
 \atcount@#1\relax\advance\atcount@\m@ne\loop\ifnum\atcount@>\z@
 \toks@=\expandafter{\the\toks@&\hfil$\m@th\displaystyle{\the\hashtoks@}$%
  \tabskip\z@skip&$\m@th\displaystyle{{}\the\hashtoks@}$\hfil
  \tabskip\centering@}\advance\atcount@\m@ne\repeat
 \xdef\preamble@{\the\toks@\tabskip\z@skip}\xdef\preamble@@{\preamble@}}
\Invalid@\endxxalignat
\newdimen\gwidth@
\newdimen\gmaxwidth@
\def\gmeasure@#1\endgather{\gwidth@\z@\gmaxwidth@\z@\setbox@ne\vbox{\Let@
 \halign{\setboxz@h{$\m@th\displaystyle{##}$}\global\gwidth@\wdz@
 \ifdim\gwidth@>\gmaxwidth@\global\gmaxwidth@\gwidth@\fi
 &\eat@{##}\crcr#1\crcr}}}
\def\gather{\RIfMIfI@\DN@{\onlydmatherr@\gather}\else
 \ingather@true\inany@true\def\tag{&}%
 \vspace@\allowdisplaybreak@\displaybreak@\intertext@
 \displ@y\Let@
 \iftagsleft@\DN@{\csname gather \endcsname}\else
  \DN@{\csname gather \space\endcsname}\fi\fi
 \else\DN@{\onlydmatherr@\gather}\fi\next@}
\expandafter\def\csname gather \space\endcsname#1\endgather
 {\gmeasure@#1\endgather\tabskip\centering@
 \halign to\displaywidth{\hfil\strut@\setboxz@h{$\m@th\displaystyle{##}$}%
 \global\gwidth@\wdz@\boxz@\hfil&
 \setboxz@h{\strut@{\maketag@##\maketag@}}%
 \dimen@\displaywidth\advance\dimen@-\gwidth@
 \ifdim\dimen@>\tw@\wdz@\llap{\boxz@}\else
 \llap{\vtop{\normalbaselines\null\boxz@}}\fi
 \tabskip\z@skip\crcr#1\crcr\black@\gmaxwidth@}}
\newdimen\glineht@
\expandafter\def\csname gather \endcsname#1\endgather{\gmeasure@#1\endgather
 \ifdim\gmaxwidth@>\displaywidth\let\gdisplaywidth@\gmaxwidth@\else
 \let\gdisplaywidth@\displaywidth\fi\tabskip\centering@\halign to\displaywidth
 {\hfil\strut@\setboxz@h{$\m@th\displaystyle{##}$}%
 \global\gwidth@\wdz@\global\glineht@\ht\z@\boxz@\hfil&\kern-\gdisplaywidth@
 \setboxz@h{\strut@{\maketag@##\maketag@}}%
 \dimen@\displaywidth\advance\dimen@-\gwidth@
 \ifdim\dimen@>\tw@\wdz@\rlap{\boxz@}\else
 \rlap{\vbox{\normalbaselines\boxz@\vbox to\glineht@{}}}\fi
 \tabskip\gdisplaywidth@\crcr#1\crcr\black@\gmaxwidth@}}
\newif\ifctagsplit@
\def\CenteredTagsOnSplits{\global\ctagsplit@true}
\def\TopOrBottomTagsOnSplits{\global\ctagsplit@false}
\TopOrBottomTagsOnSplits
\def\split{\relax\ifinany@\let\next@\insplit@\else
 \ifmmode\ifinner\def\next@{\onlydmatherr@\split}\else
 \let\next@\outsplit@\fi\else
 \def\next@{\onlydmatherr@\split}\fi\fi\next@}
\def\insplit@{\global\setbox\z@\vbox\bgroup\vspace@\Let@\ialign\bgroup
 \hfil\strut@$\m@th\displaystyle{##}$&$\m@th\displaystyle{{}##}$\hfill\crcr}
\def\endsplit{\crcr\egroup\egroup\iftagsleft@\expandafter\lendsplit@\else
 \expandafter\rendsplit@\fi}
\def\rendsplit@{\global\setbox9 \vbox
 {\unvcopy\z@\global\setbox8 \lastbox\unskip}
 \setbox@ne\hbox{\unhcopy8 \unskip\global\setbox\tw@\lastbox
 \unskip\global\setbox\thr@@\lastbox}
 \global\setbox7 \hbox{\unhbox\tw@\unskip}
 \ifinalign@\ifctagsplit@                                                   
  \gdef\split@{\hbox to\wd\thr@@{}&
   \vcenter{\vbox{\moveleft\wd\thr@@\boxz@}}}
 \else\gdef\split@{&\vbox{\moveleft\wd\thr@@\box9}\crcr
  \box\thr@@&\box7}\fi                                                      
 \else                                                                      
  \ifctagsplit@\gdef\split@{\vcenter{\boxz@}}\else
  \gdef\split@{\box9\crcr\hbox{\box\thr@@\box7}}\fi
 \fi
 \split@}                                                                   
\def\lendsplit@{\global\setbox9\vtop{\unvcopy\z@}
 \setbox@ne\vbox{\unvcopy\z@\global\setbox8\lastbox}
 \setbox@ne\hbox{\unhcopy8\unskip\setbox\tw@\lastbox
  \unskip\global\setbox\thr@@\lastbox}
 \ifinalign@\ifctagsplit@                                                   
  \gdef\split@{\hbox to\wd\thr@@{}&
  \vcenter{\vbox{\moveleft\wd\thr@@\box9}}}
  \else                                                                     
  \gdef\split@{\hbox to\wd\thr@@{}&\vbox{\moveleft\wd\thr@@\box9}}\fi
 \else
  \ifctagsplit@\gdef\split@{\vcenter{\box9}}\else
  \gdef\split@{\box9}\fi
 \fi\split@}
\def\outsplit@#1$${\align\insplit@#1\endalign$$}
\newdimen\multlinegap@
\multlinegap@1em
\newdimen\multlinetaggap@
\multlinetaggap@1em
\def\MultlineGap#1{\global\multlinegap@#1\relax}
\def\multlinegap#1{\RIfMIfI@\onlydmatherr@\multlinegap\else
 \multlinegap@#1\relax\fi\else\onlydmatherr@\multlinegap\fi}
\def\nomultlinegap{\multlinegap{\z@}}
\def\multline{\RIfMIfI@
 \DN@{\onlydmatherr@\multline}\else
 \DN@{\multline@}\fi\else
 \DN@{\onlydmatherr@\multline}\fi\next@}
\newif\iftagin@
\def\tagin@#1{\tagin@false\in@\tag{#1}\ifin@\tagin@true\fi}
\def\multline@#1$${\inany@true\vspace@\allowdisplaybreak@\displaybreak@
 \tagin@{#1}\iftagsleft@\DN@{\multline@l#1$$}\else
 \DN@{\multline@r#1$$}\fi\next@}
\newdimen\mwidth@
\def\rmmeasure@#1\endmultline{%
 \def\shoveleft##1{##1}\def\shoveright##1{##1}
 \setbox@ne\vbox{\Let@\halign{\setboxz@h
  {$\m@th\@lign\displaystyle{}##$}\global\mwidth@\wdz@
  \crcr#1\crcr}}}
\newdimen\mlineht@
\newif\ifzerocr@
\newif\ifonecr@
\def\lmmeasure@#1\endmultline{\global\zerocr@true\global\onecr@false
 \everycr{\noalign{\ifonecr@\global\onecr@false\fi
  \ifzerocr@\global\zerocr@false\global\onecr@true\fi}}
  \def\shoveleft##1{##1}\def\shoveright##1{##1}%
 \setbox@ne\vbox{\Let@\halign{\setboxz@h
  {$\m@th\@lign\displaystyle{}##$}\ifonecr@\global\mwidth@\wdz@
  \global\mlineht@\ht\z@\fi\crcr#1\crcr}}}
\newbox\mtagbox@
\newdimen\ltwidth@
\newdimen\rtwidth@
\def\multline@l#1$${\iftagin@\DN@{\lmultline@@#1$$}\else
 \DN@{\setbox\mtagbox@\null\ltwidth@\z@\rtwidth@\z@
  \lmultline@@@#1$$}\fi\next@}
\def\lmultline@@#1\endmultline\tag#2$${%
 \setbox\mtagbox@\hbox{\maketag@#2\maketag@}
 \lmmeasure@#1\endmultline\dimen@\mwidth@\advance\dimen@\wd\mtagbox@
 \advance\dimen@\multlinetaggap@                                            
 \ifdim\dimen@>\displaywidth\ltwidth@\z@\else\ltwidth@\wd\mtagbox@\fi       
 \lmultline@@@#1\endmultline$$}
\def\lmultline@@@{\displ@y
 \def\shoveright##1{##1\hfilneg\hskip\multlinegap@}%
 \def\shoveleft##1{\setboxz@h{$\m@th\displaystyle{}##1$}%
  \setbox@ne\hbox{$\m@th\displaystyle##1$}%
  \hfilneg
  \iftagin@
   \ifdim\ltwidth@>\z@\hskip\ltwidth@\hskip\multlinetaggap@\fi
  \else\hskip\multlinegap@\fi\hskip.5\wd@ne\hskip-.5\wdz@##1}
  \halign\bgroup\Let@\hbox to\displaywidth
   {\strut@$\m@th\displaystyle\hfil{}##\hfil$}\crcr
   \hfilneg                                                                 
   \iftagin@                                                                
    \ifdim\ltwidth@>\z@                                                     
     \box\mtagbox@\hskip\multlinetaggap@                                    
    \else
     \rlap{\vbox{\normalbaselines\hbox{\strut@\box\mtagbox@}%
     \vbox to\mlineht@{}}}\fi                                               
   \else\hskip\multlinegap@\fi}                                             
\def\multline@r#1$${\iftagin@\DN@{\rmultline@@#1$$}\else
 \DN@{\setbox\mtagbox@\null\ltwidth@\z@\rtwidth@\z@
  \rmultline@@@#1$$}\fi\next@}
\def\rmultline@@#1\endmultline\tag#2$${\ltwidth@\z@
 \setbox\mtagbox@\hbox{\maketag@#2\maketag@}%
 \rmmeasure@#1\endmultline\dimen@\mwidth@\advance\dimen@\wd\mtagbox@
 \advance\dimen@\multlinetaggap@
 \ifdim\dimen@>\displaywidth\rtwidth@\z@\else\rtwidth@\wd\mtagbox@\fi
 \rmultline@@@#1\endmultline$$}
\def\rmultline@@@{\displ@y
 \def\shoveright##1{##1\hfilneg\iftagin@\ifdim\rtwidth@>\z@
  \hskip\rtwidth@\hskip\multlinetaggap@\fi\else\hskip\multlinegap@\fi}%
 \def\shoveleft##1{\setboxz@h{$\m@th\displaystyle{}##1$}%
  \setbox@ne\hbox{$\m@th\displaystyle##1$}%
  \hfilneg\hskip\multlinegap@\hskip.5\wd@ne\hskip-.5\wdz@##1}%
 \halign\bgroup\Let@\hbox to\displaywidth
  {\strut@$\m@th\displaystyle\hfil{}##\hfil$}\crcr
 \hfilneg\hskip\multlinegap@}
\def\endmultline{\iftagsleft@\expandafter\lendmultline@\else
 \expandafter\rendmultline@\fi}
\def\lendmultline@{\hfilneg\hskip\multlinegap@\crcr\egroup}
\def\rendmultline@{\iftagin@                                                
 \ifdim\rtwidth@>\z@                                                        
  \hskip\multlinetaggap@\box\mtagbox@                                       
 \else\llap{\vtop{\normalbaselines\null\hbox{\strut@\box\mtagbox@}}}\fi     
 \else\hskip\multlinegap@\fi                                                
 \hfilneg\crcr\egroup}
\def\bmod{\mskip-\medmuskip\mkern5mu\mathbin{\fam\z@ mod}\penalty900
 \mkern5mu\mskip-\medmuskip}
\def\pmod#1{\allowbreak\ifinner\mkern8mu\else\mkern18mu\fi
 ({\fam\z@ mod}\,\,#1)}
\def\pod#1{\allowbreak\ifinner\mkern8mu\else\mkern18mu\fi(#1)}
\def\mod#1{\allowbreak\ifinner\mkern12mu\else\mkern18mu\fi{\fam\z@ mod}\,\,#1}
\message{continued fractions,}
\newcount\cfraccount@
\def\cfrac{\bgroup\bgroup\advance\cfraccount@\@ne\strut
 \iffalse{\fi\def\\{\over\displaystyle}\iffalse}\fi}
\def\lcfrac{\bgroup\bgroup\advance\cfraccount@\@ne\strut
 \iffalse{\fi\def\\{\hfill\over\displaystyle}\iffalse}\fi}
\def\rcfrac{\bgroup\bgroup\advance\cfraccount@\@ne\strut\hfill
 \iffalse{\fi\def\\{\over\displaystyle}\iffalse}\fi}
\def\gloop@#1\repeat{\gdef\body{#1}\iterate}
\def\endcfrac{\gloop@\ifnum\cfraccount@>\z@\global\advance\cfraccount@\m@ne
 \egroup\hskip-\nulldelimiterspace\egroup\repeat}
\message{compound symbols,}
\def\binrel@#1{\setboxz@h{\thinmuskip0mu
  \medmuskip\m@ne mu\thickmuskip\@ne mu$#1\m@th$}%
 \setbox@ne\hbox{\thinmuskip0mu\medmuskip\m@ne mu\thickmuskip
  \@ne mu${}#1{}\m@th$}%
 \setbox\tw@\hbox{\hskip\wd@ne\hskip-\wdz@}}
\def\overset#1\to#2{\binrel@{#2}\ifdim\wd\tw@<\z@
 \mathbin{\mathop{\kern\z@#2}\limits^{#1}}\else\ifdim\wd\tw@>\z@
 \mathrel{\mathop{\kern\z@#2}\limits^{#1}}\else
 {\mathop{\kern\z@#2}\limits^{#1}}{}\fi\fi}
\def\underset#1\to#2{\binrel@{#2}\ifdim\wd\tw@<\z@
 \mathbin{\mathop{\kern\z@#2}\limits_{#1}}\else\ifdim\wd\tw@>\z@
 \mathrel{\mathop{\kern\z@#2}\limits_{#1}}\else
 {\mathop{\kern\z@#2}\limits_{#1}}{}\fi\fi}
\def\oversetbrace#1\to#2{\overbrace{#2}^{#1}}
\def\undersetbrace#1\to#2{\underbrace{#2}_{#1}}
\def\sideset#1\and#2\to#3{%
 \setbox@ne\hbox{$\dsize{\vphantom{#3}}#1{#3}\m@th$}%
 \setbox\tw@\hbox{$\dsize{#3}#2\m@th$}%
 \hskip\wd@ne\hskip-\wd\tw@\mathop{\hskip\wd\tw@\hskip-\wd@ne
  {\vphantom{#3}}#1{#3}#2}}
\def\rightarrowfill@#1{\setboxz@h{$#1-\m@th$}\ht\z@\z@
  $#1\m@th\copy\z@\mkern-6mu\cleaders
  \hbox{$#1\mkern-2mu\box\z@\mkern-2mu$}\hfill
  \mkern-6mu\mathord\rightarrow$}
\def\leftarrowfill@#1{\setboxz@h{$#1-\m@th$}\ht\z@\z@
  $#1\m@th\mathord\leftarrow\mkern-6mu\cleaders
  \hbox{$#1\mkern-2mu\copy\z@\mkern-2mu$}\hfill
  \mkern-6mu\box\z@$}
\def\leftrightarrowfill@#1{\setboxz@h{$#1-\m@th$}\ht\z@\z@
  $#1\m@th\mathord\leftarrow\mkern-6mu\cleaders
  \hbox{$#1\mkern-2mu\box\z@\mkern-2mu$}\hfill
  \mkern-6mu\mathord\rightarrow$}
\def\overrightarrow{\mathpalette\overrightarrow@}
\def\overrightarrow@#1#2{\vbox{\ialign{##\crcr\rightarrowfill@#1\crcr
 \noalign{\kern-\ex@\nointerlineskip}$\m@th\hfil#1#2\hfil$\crcr}}}

\def\overleftarrow{\mathpalette\overleftarrow@}
\def\overleftarrow@#1#2{\vbox{\ialign{##\crcr\leftarrowfill@#1\crcr
 \noalign{\kern-\ex@\nointerlineskip}$\m@th\hfil#1#2\hfil$\crcr}}}
\def\overleftrightarrow{\mathpalette\overleftrightarrow@}
\def\overleftrightarrow@#1#2{\vbox{\ialign{##\crcr\leftrightarrowfill@#1\crcr
 \noalign{\kern-\ex@\nointerlineskip}$\m@th\hfil#1#2\hfil$\crcr}}}
\def\underrightarrow{\mathpalette\underrightarrow@}
\def\underrightarrow@#1#2{\vtop{\ialign{##\crcr$\m@th\hfil#1#2\hfil$\crcr
 \noalign{\nointerlineskip}\rightarrowfill@#1\crcr}}}

\def\underleftarrow{\mathpalette\underleftarrow@}
\def\underleftarrow@#1#2{\vtop{\ialign{##\crcr$\m@th\hfil#1#2\hfil$\crcr
 \noalign{\nointerlineskip}\leftarrowfill@#1\crcr}}}
\def\underleftrightarrow{\mathpalette\underleftrightarrow@}
\def\underleftrightarrow@#1#2{\vtop{\ialign{##\crcr$\m@th\hfil#1#2\hfil$\crcr
 \noalign{\nointerlineskip}\leftrightarrowfill@#1\crcr}}}
\message{various kinds of dots,}
\let\DOTSI\relax
\let\DOTSB\relax

\newif\ifmath@
{\uccode`7=`\\ \uccode`8=`m \uccode`9=`a \uccode`0=`t \uccode`!=`h
 \uppercase{\gdef\math@#1#2#3#4#5#6\math@{\global\math@false\ifx 7#1\ifx 8#2%
 \ifx 9#3\ifx 0#4\ifx !#5\xdef\meaning@{#6}\global\math@true\fi\fi\fi\fi\fi}}}
\newif\ifmathch@
{\uccode`7=`c \uccode`8=`h \uccode`9=`\"
 \uppercase{\gdef\mathch@#1#2#3#4#5#6\mathch@{\global\mathch@false
  \ifx 7#1\ifx 8#2\ifx 9#5\global\mathch@true\xdef\meaning@{9#6}\fi\fi\fi}}}
\newcount\classnum@
\def\getmathch@#1.#2\getmathch@{\classnum@#1 \divide\classnum@4096
 \ifcase\number\classnum@\or\or\gdef\thedots@{\dotsb@}\or
 \gdef\thedots@{\dotsb@}\fi}
\newif\ifmathbin@
{\uccode`4=`b \uccode`5=`i \uccode`6=`n
 \uppercase{\gdef\mathbin@#1#2#3{\relaxnext@
  \DNii@##1\mathbin@{\ifx\space@\next\global\mathbin@true\fi}%
 \global\mathbin@false\DN@##1\mathbin@{}%
 \ifx 4#1\ifx 5#2\ifx 6#3\DN@{\FN@\nextii@}\fi\fi\fi\next@}}}
\newif\ifmathrel@
{\uccode`4=`r \uccode`5=`e \uccode`6=`l
 \uppercase{\gdef\mathrel@#1#2#3{\relaxnext@
  \DNii@##1\mathrel@{\ifx\space@\next\global\mathrel@true\fi}%
 \global\mathrel@false\DN@##1\mathrel@{}%
 \ifx 4#1\ifx 5#2\ifx 6#3\DN@{\FN@\nextii@}\fi\fi\fi\next@}}}
\newif\ifmacro@
{\uccode`5=`m \uccode`6=`a \uccode`7=`c
 \uppercase{\gdef\macro@#1#2#3#4\macro@{\global\macro@false
  \ifx 5#1\ifx 6#2\ifx 7#3\global\macro@true
  \xdef\meaning@{\macro@@#4\macro@@}\fi\fi\fi}}}
\def\macro@@#1->#2\macro@@{#2}
\newif\ifDOTS@
\newcount\DOTSCASE@
{\uccode`6=`\\ \uccode`7=`D \uccode`8=`O \uccode`9=`T \uccode`0=`S
 \uppercase{\gdef\DOTS@#1#2#3#4#5{\global\DOTS@false\DN@##1\DOTS@{}%
  \ifx 6#1\ifx 7#2\ifx 8#3\ifx 9#4\ifx 0#5\let\next@\DOTS@@\fi\fi\fi\fi\fi
  \next@}}}
{\uccode`3=`B \uccode`4=`I \uccode`5=`X
 \uppercase{\gdef\DOTS@@#1{\relaxnext@
  \DNii@##1\DOTS@{\ifx\space@\next\global\DOTS@true\fi}%
  \DN@{\FN@\nextii@}%
  \ifx 3#1\global\DOTSCASE@\z@\else
  \ifx 4#1\global\DOTSCASE@\@ne\else
  \ifx 5#1\global\DOTSCASE@\tw@\else\DN@##1\DOTS@{}%
  \fi\fi\fi\next@}}}
\newif\ifnot@
{\uccode`5=`\\ \uccode`6=`n \uccode`7=`o \uccode`8=`t
 \uppercase{\gdef\not@#1#2#3#4{\relaxnext@
  \DNii@##1\not@{\ifx\space@\next\global\not@true\fi}%
 \global\not@false\DN@##1\not@{}%
 \ifx 5#1\ifx 6#2\ifx 7#3\ifx 8#4\DN@{\FN@\nextii@}\fi\fi\fi
 \fi\next@}}}
\newif\ifkeybin@
\def\keybin@{\keybin@true
 \ifx\next+\else\ifx\next=\else\ifx\next<\else\ifx\next>\else\ifx\next-\else
 \ifx\next*\else\ifx\next:\else\keybin@false\fi\fi\fi\fi\fi\fi\fi}
\def\dots{\RIfM@\expandafter\mdots@\else\expandafter\tdots@\fi}
\def\tdots@{\unskip\relaxnext@
 \DN@{$\m@th\mathinner{\ldotp\ldotp\ldotp}\,
   \ifx\next,\,$\else\ifx\next.\,$\else\ifx\next;\,$\else\ifx\next:\,$\else
   \ifx\next?\,$\else\ifx\next!\,$\else$ \fi\fi\fi\fi\fi\fi}%
 \ \FN@\next@}
\def\mdots@{\FN@\mdots@@}
\def\mdots@@{\gdef\thedots@{\dotso@}
 \ifx\next\boldkey\gdef\thedots@\boldkey{\boldkeydots@}\else                
 \ifx\next\boldsymbol\gdef\thedots@\boldsymbol{\boldsymboldots@}\else       
 \ifx,\next\gdef\thedots@{\dotsc}
 \else\ifx\not\next\gdef\thedots@{\dotsb@}
 \else\keybin@
 \ifkeybin@\gdef\thedots@{\dotsb@}
 \else\xdef\meaning@{\meaning\next..........}\xdef\meaning@@{\meaning@}
  \expandafter\math@\meaning@\math@
  \ifmath@
   \expandafter\mathch@\meaning@\mathch@
   \ifmathch@\expandafter\getmathch@\meaning@\getmathch@\fi                 
  \else\expandafter\macro@\meaning@@\macro@                                 
  \ifmacro@                                                                
   \expandafter\not@\meaning@\not@\ifnot@\gdef\thedots@{\dotsb@}
  \else\expandafter\DOTS@\meaning@\DOTS@
  \ifDOTS@
   \ifcase\number\DOTSCASE@\gdef\thedots@{\dotsb@}%
    \or\gdef\thedots@{\dotsi}\else\fi                                      
  \else\expandafter\math@\meaning@\math@                                   
  \ifmath@\expandafter\mathbin@\meaning@\mathbin@
  \ifmathbin@\gdef\thedots@{\dotsb@}
  \else\expandafter\mathrel@\meaning@\mathrel@
  \ifmathrel@\gdef\thedots@{\dotsb@}
  \fi\fi\fi\fi\fi\fi\fi\fi\fi\fi\fi\fi
 \thedots@}
\def\plainldots@{\mathinner{\ldotp\ldotp\ldotp}}
\def\plaincdots@{\mathinner{\cdotp\cdotp\cdotp}}
\def\dotsi{\!\plaincdots@}
\let\dotsb@\plaincdots@
\newif\ifextra@
\newif\ifrightdelim@
\def\rightdelim@{\global\rightdelim@true                                    
 \ifx\next)\else                                                            
 \ifx\next]\else
 \ifx\next\rbrack\else
 \ifx\next\}\else
 \ifx\next\rbrace\else
 \ifx\next\rangle\else
 \ifx\next\rceil\else
 \ifx\next\rfloor\else
 \ifx\next\rgroup\else
 \ifx\next\rmoustache\else
 \ifx\next\right\else
 \ifx\next\bigr\else
 \ifx\next\biggr\else
 \ifx\next\Bigr\else                                                        
 \ifx\next\Biggr\else\global\rightdelim@false
 \fi\fi\fi\fi\fi\fi\fi\fi\fi\fi\fi\fi\fi\fi\fi}
\def\extra@{%
 \global\extra@false\rightdelim@\ifrightdelim@\global\extra@true            
 \else\ifx\next$\global\extra@true                                          
 \else\xdef\meaning@{\meaning\next..........}
 \expandafter\macro@\meaning@\macro@\ifmacro@                               
 \expandafter\DOTS@\meaning@\DOTS@
 \ifDOTS@
 \ifnum\DOTSCASE@=\tw@\global\extra@true                                    
 \fi\fi\fi\fi\fi}
\newif\ifbold@
\def\dotso@{\relaxnext@
 \ifbold@
  \let\next\delayed@
  \DNii@{\extra@\plainldots@\ifextra@\,\fi}%
 \else
  \DNii@{\DN@{\extra@\plainldots@\ifextra@\,\fi}\FN@\next@}%
 \fi
 \nextii@}
\def\extrap@#1{%
 \ifx\next,\DN@{#1\,}\else
 \ifx\next;\DN@{#1\,}\else
 \ifx\next.\DN@{#1\,}\else\extra@
 \ifextra@\DN@{#1\,}\else
 \let\next@#1\fi\fi\fi\fi\next@}
\def\ldots{\DN@{\extrap@\plainldots@}%
 \FN@\next@}
\def\cdots{\DN@{\extrap@\plaincdots@}%
 \FN@\next@}

\def\dotsc{\relaxnext@
 \DN@{\ifx\next;\plainldots@\,\else
  \ifx\next.\plainldots@\,\else\extra@\plainldots@
  \ifextra@\,\fi\fi\fi}%
 \FN@\next@}
\def\cdot{\mathchar"2201 }

\def\mapsto{\DOTSB\mapstochar\rightarrow}

\message{special superscripts,}
\def\dddot#1{{\mathop{#1}\limits^{\vbox to-1.4\ex@{\kern-\tw@\ex@
 \hbox{\rm...}\vss}}}}
\def\ddddot#1{{\mathop{#1}\limits^{\vbox to-1.4\ex@{\kern-\tw@\ex@
 \hbox{\rm....}\vss}}}}
\def\sphat{^{\mathchoice{}{}%
 {\,\,\botsmash{\hbox{\lower4\ex@\hbox{$\m@th\widehat{\null}$}}}}%
 {\,\botsmash{\hbox{\lower3\ex@\hbox{$\m@th\hat{\null}$}}}}}}

\def\spacute{^{\!\botsmash{\hbox{\lower\@ne ex\hbox{\'{}}}}}}
\def\spgrave{^{\mathchoice{}{}{}{\!}%
 \botsmash{\hbox{\lower\@ne ex\hbox{\`{}}}}}}
\def\spdot{^{\hbox{\raise\ex@\hbox{\rm.}}}}
\def\spddot{^{\hbox{\raise\ex@\hbox{\rm..}}}}
\def\spdddot{^{\hbox{\raise\ex@\hbox{\rm...}}}}
\def\spddddot{^{\hbox{\raise\ex@\hbox{\rm....}}}}
\def\spbreve{^{\!\botsmash{\hbox{\lower4\ex@\hbox{\u{}}}}}}

\message{\string\text,}
\def\textonlyfont@#1#2{\def#1{\RIfM@
 \Err@{Use \string#1\space only in text}\else#2\fi}}
\textonlyfont@\rm\tenrm
\textonlyfont@\it\tenit
\textonlyfont@\sl\tensl
\textonlyfont@\bf\tenbf
\def\oldnos#1{\RIfM@{\mathcode`\,="013B \fam\@ne#1}\else
 \leavevmode\hbox{$\m@th\mathcode`\,="013B \fam\@ne#1$}\fi}
\def\text{\RIfM@\expandafter\text@\else\expandafter\text@@\fi}
\def\text@@#1{\leavevmode\hbox{#1}}
\def\mathhexbox@#1#2#3{\text{$\m@th\mathchar"#1#2#3$}}
\def\dag{{\mathhexbox@279}}
\def\ddag{{\mathhexbox@27A}}
\def\S{{\mathhexbox@278}}
\def\P{{\mathhexbox@27B}}
\newif\iffirstchoice@
\firstchoice@true
\def\text@#1{\mathchoice
 {\hbox{\everymath{\displaystyle}\def\textfonti{\the\textfont\@ne}%
  \def\textfontii{\the\textfont\tw@}\textdef@@ T#1}}
 {\hbox{\firstchoice@false
  \everymath{\textstyle}\def\textfonti{\the\textfont\@ne}%
  \def\textfontii{\the\textfont\tw@}\textdef@@ T#1}}
 {\hbox{\firstchoice@false
  \everymath{\scriptstyle}\def\textfonti{\the\scriptfont\@ne}%
  \def\textfontii{\the\scriptfont\tw@}\textdef@@ S\rm#1}}
 {\hbox{\firstchoice@false
  \everymath{\scriptscriptstyle}\def\textfonti
  {\the\scriptscriptfont\@ne}%
  \def\textfontii{\the\scriptscriptfont\tw@}\textdef@@ s\rm#1}}}
\def\textdef@@#1{\textdef@#1\rm\textdef@#1\bf\textdef@#1\sl\textdef@#1\it}
\def\rmfam{0}
\def\textdef@#1#2{%
 \DN@{\csname\expandafter\eat@\string#2fam\endcsname}%
 \if S#1\edef#2{\the\scriptfont\next@\relax}%
 \else\if s#1\edef#2{\the\scriptscriptfont\next@\relax}%
 \else\edef#2{\the\textfont\next@\relax}\fi\fi}
\scriptfont\itfam\tenit \scriptscriptfont\itfam\tenit
\scriptfont\slfam\tensl \scriptscriptfont\slfam\tensl
\newif\iftopfolded@
\newif\ifbotfolded@
\def\topfoldedtext{\topfolded@true\botfolded@false\foldedtext@}
\def\botfoldedtext{\botfolded@true\topfolded@false\foldedtext@}
\def\foldedtext{\topfolded@false\botfolded@false\foldedtext@}
\Invalid@\foldedwidth
\def\foldedtext@{\relaxnext@
 \DN@{\ifx\next\foldedwidth\let\next@\nextii@\else
  \DN@{\nextii@\foldedwidth{.3\hsize}}\fi\next@}%
 \DNii@\foldedwidth##1##2{\setbox\z@\vbox
  {\normalbaselines\hsize##1\relax
  \tolerance1600 \noindent\ignorespaces##2}\ifbotfolded@\boxz@\else
  \iftopfolded@\vtop{\unvbox\z@}\else\vcenter{\boxz@}\fi\fi}%
 \FN@\next@}
\message{math font commands,}
\def\bold{\RIfM@\expandafter\bold@\else
 \expandafter\nonmatherr@\expandafter\bold\fi}
\def\bold@#1{{\bold@@{#1}}}
\def\bold@@#1{\fam\bffam\relax#1}
\def\slanted{\RIfM@\expandafter\slanted@\else
 \expandafter\nonmatherr@\expandafter\slanted\fi}
\def\slanted@#1{{\slanted@@{#1}}}
\def\slanted@@#1{\fam\slfam\relax#1}
\def\roman{\RIfM@\expandafter\roman@\else
 \expandafter\nonmatherr@\expandafter\roman\fi}
\def\roman@#1{{\roman@@{#1}}}
\def\roman@@#1{\fam\rmfam\relax#1}
\def\italic{\RIfM@\expandafter\italic@\else
 \expandafter\nonmatherr@\expandafter\italic\fi}
\def\italic@#1{{\italic@@{#1}}}
\def\italic@@#1{\fam\itfam\relax#1}
\def\Cal{\RIfM@\expandafter\Cal@\else
 \expandafter\nonmatherr@\expandafter\Cal\fi}
\def\Cal@#1{{\Cal@@{#1}}}
\def\Cal@@#1{\noaccents@\fam\tw@#1}
\mathchardef\Gamma="0000
\mathchardef\Delta="0001
\mathchardef\Theta="0002
\mathchardef\Lambda="0003
\mathchardef\Xi="0004
\mathchardef\Pi="0005
\mathchardef\Sigma="0006
\mathchardef\Upsilon="0007
\mathchardef\Phi="0008
\mathchardef\Psi="0009
\mathchardef\Omega="000A
\mathchardef\varGamma="0100
\mathchardef\varDelta="0101
\mathchardef\varTheta="0102
\mathchardef\varLambda="0103
\mathchardef\varXi="0104
\mathchardef\varPi="0105
\mathchardef\varSigma="0106
\mathchardef\varUpsilon="0107
\mathchardef\varPhi="0108
\mathchardef\varPsi="0109
\mathchardef\varOmega="010A
\let\alloc@@\alloc@
\def\hexnumber@#1{\ifcase#1 0\or 1\or 2\or 3\or 4\or 5\or 6\or 7\or 8\or
 9\or A\or B\or C\or D\or E\or F\fi}
\def\loadmsam{%
 \font@\tenmsa=msam10
 \font@\sevenmsa=msam7
 \font@\fivemsa=msam5
 \alloc@@8\fam\chardef\sixt@@n\msafam
 \textfont\msafam=\tenmsa
 \scriptfont\msafam=\sevenmsa
 \scriptscriptfont\msafam=\fivemsa
 \edef\next{\hexnumber@\msafam}%
 \mathchardef\dabar@"0\next39
 \edef\dashrightarrow{\mathrel{\dabar@\dabar@\mathchar"0\next4B}}%
 \edef\dashleftarrow{\mathrel{\mathchar"0\next4C\dabar@\dabar@}}%
 \let\dasharrow\dashrightarrow
 \edef\ulcorner{\delimiter"4\next70\next70 }%
 \edef\urcorner{\delimiter"5\next71\next71 }%
 \edef\llcorner{\delimiter"4\next78\next78 }%
 \edef\lrcorner{\delimiter"5\next79\next79 }%
 \edef\yen{{\noexpand\mathhexbox@\next55}}%
 \edef\checkmark{{\noexpand\mathhexbox@\next58}}%
 \edef\circledR{{\noexpand\mathhexbox@\next72}}%
 \edef\maltese{{\noexpand\mathhexbox@\next7A}}%
 \global\let\loadmsam\empty}%
\def\loadmsbm{%
 \font@\tenmsb=msbm10 \font@\sevenmsb=msbm7 \font@\fivemsb=msbm5
 \alloc@@8\fam\chardef\sixt@@n\msbfam
 \textfont\msbfam=\tenmsb
 \scriptfont\msbfam=\sevenmsb \scriptscriptfont\msbfam=\fivemsb
 \global\let\loadmsbm\empty
 }
\def\widehat#1{\ifx\undefined\msbfam \DN@{362}%
  \else \setboxz@h{$\m@th#1$}%
    \edef\next@{\ifdim\wdz@>\tw@ em%
        \hexnumber@\msbfam 5B%
      \else 362\fi}\fi
  \mathaccent"0\next@{#1}}
\def\widetilde#1{\ifx\undefined\msbfam \DN@{365}%
  \else \setboxz@h{$\m@th#1$}%
    \edef\next@{\ifdim\wdz@>\tw@ em%
        \hexnumber@\msbfam 5D%
      \else 365\fi}\fi
  \mathaccent"0\next@{#1}}
\message{\string\newsymbol,}
\def\newsymbol#1#2#3#4#5{\define#1{}%
  \count@#2\relax \advance\count@\m@ne 
 \ifcase\count@
   \ifx\undefined\msafam\loadmsam\fi \let\next@\msafam
 \or \ifx\undefined\msbfam\loadmsbm\fi \let\next@\msbfam
 \else  \Err@{\Invalid@@\string\newsymbol}\let\next@\tw@\fi
 \mathchardef#1="#3\hexnumber@\next@#4#5\space}

\def\Bbb{\RIfM@\expandafter\Bbb@\else
 \expandafter\nonmatherr@\expandafter\Bbb\fi}
\def\Bbb@#1{{\Bbb@@{#1}}}
\def\Bbb@@#1{\noaccents@\fam\msbfam\relax#1}
\message{bold Greek and bold symbols,}
\def\loadbold{%
 \font@\tencmmib=cmmib10 \font@\sevencmmib=cmmib7 \font@\fivecmmib=cmmib5
 \skewchar\tencmmib'177 \skewchar\sevencmmib'177 \skewchar\fivecmmib'177
 \alloc@@8\fam\chardef\sixt@@n\cmmibfam
 \textfont\cmmibfam\tencmmib
 \scriptfont\cmmibfam\sevencmmib \scriptscriptfont\cmmibfam\fivecmmib
 \font@\tencmbsy=cmbsy10 \font@\sevencmbsy=cmbsy7 \font@\fivecmbsy=cmbsy5
 \skewchar\tencmbsy'60 \skewchar\sevencmbsy'60 \skewchar\fivecmbsy'60
 \alloc@@8\fam\chardef\sixt@@n\cmbsyfam
 \textfont\cmbsyfam\tencmbsy
 \scriptfont\cmbsyfam\sevencmbsy \scriptscriptfont\cmbsyfam\fivecmbsy
 \let\loadbold\empty
}
\def\boldnotloaded#1{\Err@{\ifcase#1\or First\else Second\fi
       bold symbol font not loaded}}
\def\mathchari@#1#2#3{\ifx\undefined\cmmibfam
    \boldnotloaded@\@ne
  \else\mathchar"#1\hexnumber@\cmmibfam#2#3\space \fi}
\def\mathcharii@#1#2#3{\ifx\undefined\cmbsyfam
    \boldnotloaded\tw@
  \else \mathchar"#1\hexnumber@\cmbsyfam#2#3\space\fi}
\edef\bffam@{\hexnumber@\bffam}
\def\boldkey#1{\ifcat\noexpand#1A%
  \ifx\undefined\cmmibfam \boldnotloaded\@ne
  \else {\fam\cmmibfam#1}\fi
 \else
 \ifx#1!\mathchar"5\bffam@21 \else
 \ifx#1(\mathchar"4\bffam@28 \else\ifx#1)\mathchar"5\bffam@29 \else
 \ifx#1+\mathchar"2\bffam@2B \else\ifx#1:\mathchar"3\bffam@3A \else
 \ifx#1;\mathchar"6\bffam@3B \else\ifx#1=\mathchar"3\bffam@3D \else
 \ifx#1?\mathchar"5\bffam@3F \else\ifx#1[\mathchar"4\bffam@5B \else
 \ifx#1]\mathchar"5\bffam@5D \else
 \ifx#1,\mathchari@63B \else
 \ifx#1-\mathcharii@200 \else
 \ifx#1.\mathchari@03A \else
 \ifx#1/\mathchari@03D \else
 \ifx#1<\mathchari@33C \else
 \ifx#1>\mathchari@33E \else
 \ifx#1*\mathcharii@203 \else
 \ifx#1|\mathcharii@06A \else
 \ifx#10\bold0\else\ifx#11\bold1\else\ifx#12\bold2\else\ifx#13\bold3\else
 \ifx#14\bold4\else\ifx#15\bold5\else\ifx#16\bold6\else\ifx#17\bold7\else
 \ifx#18\bold8\else\ifx#19\bold9\else
  \Err@{\string\boldkey\space can't be used with #1}%
 \fi\fi\fi\fi\fi\fi\fi\fi\fi\fi\fi\fi\fi\fi\fi
 \fi\fi\fi\fi\fi\fi\fi\fi\fi\fi\fi\fi\fi\fi}
\def\boldsymbol#1{%
 \DN@{\Err@{You can't use \string\boldsymbol\space with \string#1}#1}%
 \ifcat\noexpand#1A%
   \let\next@\relax
   \ifx\undefined\cmmibfam \boldnotloaded\@ne
   \else {\fam\cmmibfam#1}\fi
 \else
  \xdef\meaning@{\meaning#1.........}%
  \expandafter\math@\meaning@\math@
  \ifmath@
   \expandafter\mathch@\meaning@\mathch@
   \ifmathch@
    \expandafter\boldsymbol@@\meaning@\boldsymbol@@
   \fi
  \else
   \expandafter\macro@\meaning@\macro@
   \expandafter\delim@\meaning@\delim@
   \ifdelim@
    \expandafter\delim@@\meaning@\delim@@
   \else
    \boldsymbol@{#1}%
   \fi
  \fi
 \fi
 \next@}
\def\mathhexboxii@#1#2{\ifx\undefined\cmbsyfam
    \boldnotloaded\tw@
  \else \mathhexbox@{\hexnumber@\cmbsyfam}{#1}{#2}\fi}
\def\boldsymbol@#1{\let\next@\relax\let\next#1%
 \ifx\next\cdot\mathcharii@201 \else
 \ifx\next\prime{{\null\mathcharii@030 \null}}\else
 \ifx\next\lbrack\mathchar"4\bffam@5B \else
 \ifx\next\rbrack\mathchar"5\bffam@5D \else
 \ifx\next\{\mathcharii@466 \else
 \ifx\next\lbrace\mathcharii@466 \else
 \ifx\next\}\mathcharii@567 \else
 \ifx\next\rbrace\mathcharii@567 \else
 \ifx\next\surd{{\mathcharii@170}}\else
 \ifx\next\S{{\mathhexboxii@78}}\else
 \ifx\next\P{{\mathhexboxii@7B}}\else
 \ifx\next\dag{{\mathhexboxii@79}}\else
 \ifx\next\ddag{{\mathhexboxii@7A}}\else
 \DN@{\Err@{You can't use \string\boldsymbol\space with \string#1}#1}%
 \fi\fi\fi\fi\fi\fi\fi\fi\fi\fi\fi\fi\fi}
\def\boldsymbol@@#1.#2\boldsymbol@@{\classnum@#1 \count@@@\classnum@        
 \divide\classnum@4096 \count@\classnum@                                    
 \multiply\count@4096 \advance\count@@@-\count@ \count@@\count@@@           
 \divide\count@@@\@cclvi \count@\count@@                                    
 \multiply\count@@@\@cclvi \advance\count@@-\count@@@                       
 \divide\count@@@\@cclvi                                                    
 \multiply\classnum@4096 \advance\classnum@\count@@                         
 \ifnum\count@@@=\z@                                                        
  \count@"\bffam@ \multiply\count@\@cclvi
  \advance\classnum@\count@
  \DN@{\mathchar\number\classnum@}%
 \else
  \ifnum\count@@@=\@ne                                                      
   \ifx\undefined\cmmibfam \DN@{\boldnotloaded\@ne}%
   \else \count@\cmmibfam \multiply\count@\@cclvi
     \advance\classnum@\count@
     \DN@{\mathchar\number\classnum@}\fi
  \else
   \ifnum\count@@@=\tw@                                                    
     \ifx\undefined\cmbsyfam
       \DN@{\boldnotloaded\tw@}%
     \else
       \count@\cmbsyfam \multiply\count@\@cclvi
       \advance\classnum@\count@
       \DN@{\mathchar\number\classnum@}%
     \fi
  \fi
 \fi
\fi}
\newif\ifdelim@
\newcount\delimcount@
{\uccode`6=`\\ \uccode`7=`d \uccode`8=`e \uccode`9=`l
 \uppercase{\gdef\delim@#1#2#3#4#5\delim@
  {\delim@false\ifx 6#1\ifx 7#2\ifx 8#3\ifx 9#4\delim@true
   \xdef\meaning@{#5}\fi\fi\fi\fi}}}
\def\delim@@#1"#2#3#4#5#6\delim@@{\if#32%
\let\next@\relax
 \ifx\undefined\cmbsyfam \boldnotloaded\@ne
 \else \mathcharii@#2#4#5\space \fi\fi}
\def\vert{\delimiter"026A30C }
\def\Vert{\delimiter"026B30D }
\let\|\Vert
\def\backslash{\delimiter"026E30F }
\def\boldkeydots@#1{\bold@true\let\next=#1\let\delayed@=#1\mdots@@
 \boldkey#1\bold@false}  
\def\boldsymboldots@#1{\bold@true\let\next#1\let\delayed@#1\mdots@@
 \boldsymbol#1\bold@false}
\message{Euler fonts,}

\def\frak{\mathfont@\frak}

\def\loadmathfont#1{%
   \expandafter\font@\csname ten#1\endcsname=#110
   \expandafter\font@\csname seven#1\endcsname=#17
   \expandafter\font@\csname five#1\endcsname=#15
   \edef\next{\noexpand\alloc@@8\fam\chardef\sixt@@n
     \expandafter\noexpand\csname#1fam\endcsname}%
   \next
   \textfont\csname#1fam\endcsname \csname ten#1\endcsname
   \scriptfont\csname#1fam\endcsname \csname seven#1\endcsname
   \scriptscriptfont\csname#1fam\endcsname \csname five#1\endcsname
   \expandafter\def\csname #1\expandafter\endcsname\expandafter{%
      \expandafter\mathfont@\csname#1\endcsname}%
 \expandafter\gdef\csname load#1\endcsname{}%
}
\def\mathfont@#1{\RIfM@\expandafter\mathfont@@\expandafter#1\else
  \expandafter\nonmatherr@\expandafter#1\fi}
\def\mathfont@@#1#2{{\mathfont@@@#1{#2}}}
\def\mathfont@@@#1#2{\noaccents@
   \fam\csname\expandafter\eat@\string#1fam\endcsname
   \relax#2}
\message{math accents,}
\def\accentclass@{7}
\def\noaccents@{\def\accentclass@{0}}
\def\makeacc@#1#2{\def#1{\mathaccent"\accentclass@#2 }}
\makeacc@\hat{05E}
\makeacc@\check{014}
\makeacc@\tilde{07E}
\makeacc@\acute{013}
\makeacc@\grave{012}
\makeacc@\dot{05F}
\makeacc@\ddot{07F}
\makeacc@\breve{015}
\makeacc@\bar{016}

\newcount\skewcharcount@
\newcount\familycount@
\def\theskewchar@{\familycount@\@ne
 \global\skewcharcount@\the\skewchar\textfont\@ne                           
 \ifnum\fam>\m@ne\ifnum\fam<16
  \global\familycount@\the\fam\relax
  \global\skewcharcount@\the\skewchar\textfont\the\fam\relax\fi\fi          
 \ifnum\skewcharcount@>\m@ne
  \ifnum\skewcharcount@<128
  \multiply\familycount@256
  \global\advance\skewcharcount@\familycount@
  \global\advance\skewcharcount@28672
  \mathchar\skewcharcount@\else
  \global\skewcharcount@\m@ne\fi\else
 \global\skewcharcount@\m@ne\fi}                                            
\newcount\pointcount@
\def\getpoints@#1.#2\getpoints@{\pointcount@#1 }
\newdimen\accentdimen@
\newcount\accentmu@
\def\dimentomu@{\multiply\accentdimen@ 100
 \expandafter\getpoints@\the\accentdimen@\getpoints@
 \multiply\pointcount@18
 \divide\pointcount@\@m
 \global\accentmu@\pointcount@}
\def\Makeacc@#1#2{\def#1{\RIfM@\DN@{\mathaccent@
 {"\accentclass@#2 }}\else\DN@{\nonmatherr@{#1}}\fi\next@}}
\def\unbracefonts@{\let\Cal@\Cal@@\let\roman@\roman@@\let\bold@\bold@@
 \let\slanted@\slanted@@}
\def\mathaccent@#1#2{\ifnum\fam=\m@ne\xdef\thefam@{1}\else
 \xdef\thefam@{\the\fam}\fi                                                 
 \accentdimen@\z@                                                           
 \setboxz@h{\unbracefonts@$\m@th\fam\thefam@\relax#2$}
 \ifdim\accentdimen@=\z@\DN@{\mathaccent#1{#2}}
  \setbox@ne\hbox{\unbracefonts@$\m@th\fam\thefam@\relax#2\theskewchar@$}
  \setbox\tw@\hbox{$\m@th\ifnum\skewcharcount@=\m@ne\else
   \mathchar\skewcharcount@\fi$}
  \global\accentdimen@\wd@ne\global\advance\accentdimen@-\wdz@
  \global\advance\accentdimen@-\wd\tw@                                     
  \global\multiply\accentdimen@\tw@
  \dimentomu@\global\advance\accentmu@\@ne                                 
 \else\DN@{{\mathaccent#1{#2\mkern\accentmu@ mu}%
    \mkern-\accentmu@ mu}{}}\fi                                             
 \next@}\Makeacc@\Hat{05E}
\Makeacc@\Check{014}
\Makeacc@\Tilde{07E}
\Makeacc@\Acute{013}
\Makeacc@\Grave{012}
\Makeacc@\Dot{05F}
\Makeacc@\Ddot{07F}
\Makeacc@\Breve{015}
\Makeacc@\Bar{016}
\def\Vec{\RIfM@\DN@{\mathaccent@{"017E }}\else
 \DN@{\nonmatherr@\Vec}\fi\next@}
\def\accentedsymbol#1#2{\csname newbox\expandafter\endcsname
  \csname\expandafter\eat@\string#1@box\endcsname
 \expandafter\setbox\csname\expandafter\eat@
  \string#1@box\endcsname\hbox{$\m@th#2$}\define
  #1{\copy\csname\expandafter\eat@\string#1@box\endcsname{}}}
\message{roots,}
\def\sqrt#1{\radical"270370 {#1}}
\let\underline@\underline
\let\overline@\overline
\def\underline#1{\underline@{#1}}
\def\overline#1{\overline@{#1}}
\Invalid@\leftroot
\Invalid@\uproot
\newcount\uproot@
\newcount\leftroot@
\def\root{\relaxnext@
  \DN@{\ifx\next\uproot\let\next@\nextii@\else
   \ifx\next\leftroot\let\next@\nextiii@\else
   \let\next@\plainroot@\fi\fi\next@}%
  \DNii@\uproot##1{\uproot@##1\relax\FN@\nextiv@}%
  \def\nextiv@{\ifx\next\space@\DN@. {\FN@\nextv@}\else
   \DN@.{\FN@\nextv@}\fi\next@.}%
  \def\nextv@{\ifx\next\leftroot\let\next@\nextvi@\else
   \let\next@\plainroot@\fi\next@}%
  \def\nextvi@\leftroot##1{\leftroot@##1\relax\plainroot@}%
   \def\nextiii@\leftroot##1{\leftroot@##1\relax\FN@\nextvii@}%
  \def\nextvii@{\ifx\next\space@
   \DN@. {\FN@\nextviii@}\else
   \DN@.{\FN@\nextviii@}\fi\next@.}%
  \def\nextviii@{\ifx\next\uproot\let\next@\nextix@\else
   \let\next@\plainroot@\fi\next@}%
  \def\nextix@\uproot##1{\uproot@##1\relax\plainroot@}%
  \bgroup\uproot@\z@\leftroot@\z@\FN@\next@}
\def\plainroot@#1\of#2{\setbox\rootbox\hbox{$\m@th\scriptscriptstyle{#1}$}%
 \mathchoice{\r@@t\displaystyle{#2}}{\r@@t\textstyle{#2}}
 {\r@@t\scriptstyle{#2}}{\r@@t\scriptscriptstyle{#2}}\egroup}
\def\r@@t#1#2{\setboxz@h{$\m@th#1\sqrt{#2}$}%
 \dimen@\ht\z@\advance\dimen@-\dp\z@
 \setbox@ne\hbox{$\m@th#1\mskip\uproot@ mu$}\advance\dimen@ 1.667\wd@ne
 \mkern-\leftroot@ mu\mkern5mu\raise.6\dimen@\copy\rootbox
 \mkern-10mu\mkern\leftroot@ mu\boxz@}
\def\boxed#1{\setboxz@h{$\m@th\displaystyle{#1}$}\dimen@.4\ex@
 \advance\dimen@3\ex@\advance\dimen@\dp\z@
 \hbox{\lower\dimen@\hbox{%
 \vbox{\hrule height.4\ex@
 \hbox{\vrule width.4\ex@\hskip3\ex@\vbox{\vskip3\ex@\boxz@\vskip3\ex@}%
 \hskip3\ex@\vrule width.4\ex@}\hrule height.4\ex@}%
 }}}
\message{commutative diagrams,}
\let\ampersand@\relax
\newdimen\minaw@
\minaw@11.11128\ex@
\newdimen\minCDaw@
\minCDaw@2.5pc
\def\minCDarrowwidth#1{\RIfMIfI@\onlydmatherr@\minCDarrowwidth
 \else\minCDaw@#1\relax\fi\else\onlydmatherr@\minCDarrowwidth\fi}
\newif\ifCD@
\def\CD{\bgroup\vspace@\relax\let\ampersand@&\iffalse}\fi
 \CD@true\vcenter\bgroup\Let@\tabskip\z@skip\baselineskip20\ex@
 \lineskip3\ex@\lineskiplimit3\ex@\halign\bgroup
 &\hfill$\m@th##$\hfill\crcr}
\def\endCD{\crcr\egroup\egroup\egroup}
\newdimen\bigaw@
\atdef@>#1>#2>{\ampersand@                                                  
 \setboxz@h{$\m@th\ssize\;{#1}\;\;$}
 \setbox@ne\hbox{$\m@th\ssize\;{#2}\;\;$}
 \setbox\tw@\hbox{$\m@th#2$}
 \ifCD@\global\bigaw@\minCDaw@\else\global\bigaw@\minaw@\fi                 
 \ifdim\wdz@>\bigaw@\global\bigaw@\wdz@\fi
 \ifdim\wd@ne>\bigaw@\global\bigaw@\wd@ne\fi                                
 \ifCD@\enskip\fi                                                           
 \ifdim\wd\tw@>\z@
  \mathrel{\mathop{\hbox to\bigaw@{\rightarrowfill@\displaystyle}}%
    \limits^{#1}_{#2}}
 \else\mathrel{\mathop{\hbox to\bigaw@{\rightarrowfill@\displaystyle}}%
    \limits^{#1}}\fi                                                        
 \ifCD@\enskip\fi                                                          
 \ampersand@}                                                              
\atdef@<#1<#2<{\ampersand@\setboxz@h{$\m@th\ssize\;\;{#1}\;$}%
 \setbox@ne\hbox{$\m@th\ssize\;\;{#2}\;$}\setbox\tw@\hbox{$\m@th#2$}%
 \ifCD@\global\bigaw@\minCDaw@\else\global\bigaw@\minaw@\fi
 \ifdim\wdz@>\bigaw@\global\bigaw@\wdz@\fi
 \ifdim\wd@ne>\bigaw@\global\bigaw@\wd@ne\fi
 \ifCD@\enskip\fi
 \ifdim\wd\tw@>\z@
  \mathrel{\mathop{\hbox to\bigaw@{\leftarrowfill@\displaystyle}}%
       \limits^{#1}_{#2}}\else
  \mathrel{\mathop{\hbox to\bigaw@{\leftarrowfill@\displaystyle}}%
       \limits^{#1}}\fi
 \ifCD@\enskip\fi\ampersand@}
\begingroup
 \catcode`\~=\active \lccode`\~=`\@
 \lowercase{%
  \global\atdef@)#1)#2){~>#1>#2>}
  \global\atdef@(#1(#2({~<#1<#2<}}
\endgroup
\atdef@ A#1A#2A{\llap{$\m@th\vcenter{\hbox
 {$\ssize#1$}}$}\Big\uparrow\rlap{$\m@th\vcenter{\hbox{$\ssize#2$}}$}&&}
\atdef@ V#1V#2V{\llap{$\m@th\vcenter{\hbox
 {$\ssize#1$}}$}\Big\downarrow\rlap{$\m@th\vcenter{\hbox{$\ssize#2$}}$}&&}
\atdef@={&\enskip\mathrel
 {\vbox{\hrule width\minCDaw@\vskip3\ex@\hrule width
 \minCDaw@}}\enskip&}
\atdef@|{\Big\Vert&&}
\atdef@\vert{\Big\Vert&&}
\def\pretend#1\haswidth#2{\setboxz@h{$\m@th\scriptstyle{#2}$}\hbox
 to\wdz@{\hfill$\m@th\scriptstyle{#1}$\hfill}}
\message{poor man's bold,}
\def\pmb{\RIfM@\expandafter\mathpalette\expandafter\pmb@\else
 \expandafter\pmb@@\fi}
\def\pmb@@#1{\leavevmode\setboxz@h{#1}%
   \dimen@-\wdz@
   \kern-.5\ex@\copy\z@
   \kern\dimen@\kern.25\ex@\raise.4\ex@\copy\z@
   \kern\dimen@\kern.25\ex@\box\z@
}
\def\binrel@@#1{\ifdim\wd2<\z@\mathbin{#1}\else\ifdim\wd\tw@>\z@
 \mathrel{#1}\else{#1}\fi\fi}
\newdimen\pmbraise@
\def\pmb@#1#2{\setbox\thr@@\hbox{$\m@th#1{#2}$}%
 \setbox4\hbox{$\m@th#1\mkern.5mu$}\pmbraise@\wd4\relax
 \binrel@{#2}%
 \dimen@-\wd\thr@@
   \binrel@@{%
   \mkern-.8mu\copy\thr@@
   \kern\dimen@\mkern.4mu\raise\pmbraise@\copy\thr@@
   \kern\dimen@\mkern.4mu\box\thr@@
}}
\def\documentstyle#1{\W@{}\input #1.sty\relax}
\message{syntax check,}
\font\dummyft@=dummy
\fontdimen1 \dummyft@=\z@
\fontdimen2 \dummyft@=\z@
\fontdimen3 \dummyft@=\z@
\fontdimen4 \dummyft@=\z@
\fontdimen5 \dummyft@=\z@
\fontdimen6 \dummyft@=\z@
\fontdimen7 \dummyft@=\z@
\fontdimen8 \dummyft@=\z@
\fontdimen9 \dummyft@=\z@
\fontdimen10 \dummyft@=\z@
\fontdimen11 \dummyft@=\z@
\fontdimen12 \dummyft@=\z@
\fontdimen13 \dummyft@=\z@
\fontdimen14 \dummyft@=\z@
\fontdimen15 \dummyft@=\z@
\fontdimen16 \dummyft@=\z@
\fontdimen17 \dummyft@=\z@
\fontdimen18 \dummyft@=\z@
\fontdimen19 \dummyft@=\z@
\fontdimen20 \dummyft@=\z@
\fontdimen21 \dummyft@=\z@
\fontdimen22 \dummyft@=\z@
\def\fontlist@{\\{\tenrm}\\{\sevenrm}\\{\fiverm}\\{\teni}\\{\seveni}%
 \\{\fivei}\\{\tensy}\\{\sevensy}\\{\fivesy}\\{\tenex}\\{\tenbf}\\{\sevenbf}%
 \\{\fivebf}\\{\tensl}\\{\tenit}}
\def\font@#1=#2 {\rightappend@#1\to\fontlist@\font#1=#2 }
\def\dodummy@{{\def\\##1{\global\let##1\dummyft@}\fontlist@}}
\def\nopages@{\output{\setbox\z@\box\@cclv \deadcycles\z@}%
 \alloc@5\toks\toksdef\@cclvi\output}
\let\galleys\nopages@
\newif\ifsyntax@
\newcount\countxviii@
\def\syntax{\syntax@true\dodummy@\countxviii@\count18
 \loop\ifnum\countxviii@>\m@ne\textfont\countxviii@=\dummyft@
 \scriptfont\countxviii@=\dummyft@\scriptscriptfont\countxviii@=\dummyft@
 \advance\countxviii@\m@ne\repeat                                           
 \dummyft@\tracinglostchars\z@\nopages@\frenchspacing\hbadness\@M}
\def\first@#1#2\end{#1}
\def\printoptions{\W@{Do you want S(yntax check),
  G(alleys) or P(ages)?}%
 \message{Type S, G or P, followed by <return>: }%
 \begingroup 
 \endlinechar\m@ne 
 \read\m@ne to\ans@
 \edef\ans@{\uppercase{\def\noexpand\ans@{%
   \expandafter\first@\ans@ P\end}}}%
 \expandafter\endgroup\ans@
 \if\ans@ P
 \else \if\ans@ S\syntax
 \else \if\ans@ G\galleys
 \else\message{? Unknown option: \ans@; using the `pages' option.}%
 \fi\fi\fi}
\def\alloc@#1#2#3#4#5{\global\advance\count1#1by\@ne
 \ch@ck#1#4#2\allocationnumber=\count1#1
 \global#3#5=\allocationnumber
 \ifalloc@\wlog{\string#5=\string#2\the\allocationnumber}\fi}
\def\document{\def\alloclist@{}\def\fontlist@{}}
\let\enddocument\bye

\let\proclaim\undefined
\let\footnote\undefined
\let\=\undefined
\let\>\undefined

\catcode`\@=\active
\message{... finished}

\expandafter\ifx\csname mathdefs.tex\endcsname\relax
  \expandafter\gdef\csname mathdefs.tex\endcsname{}
\else \message{Hey!  Apparently you were trying to
  \string\input{mathdefs.tex} twice.   This does not make sense.} 
\errmessage{Please edit your file (probably \jobname.tex) and remove
any duplicate ``\string\input'' lines}\endinput\fi




\catcode`\X=12\catcode`\@=11

\def\n@wcount{\alloc@0\count\countdef\insc@unt}
\def\n@wwrite{\alloc@7\write\chardef\sixt@@n}
\def\n@wread{\alloc@6\read\chardef\sixt@@n}
\def\r@s@t{\relax}\def\v@idline{\par}\def\@mputate#1/{#1}
\def\l@c@l#1X{\firstpart.#1}\def\gl@b@l#1X{#1}\def\t@d@l#1X{{}}

\def\crossrefs#1{\ifx\all#1\let\tr@ce=\all\else\def\tr@ce{#1,}\fi
   \n@wwrite\cit@tionsout\openout\cit@tionsout=\jobname.cit 
   \write\cit@tionsout{\tr@ce}\expandafter\setfl@gs\tr@ce,}
\def\setfl@gs#1,{\def\@{#1}\ifx\@\empty\let\next=\relax
   \else\let\next=\setfl@gs\expandafter\xdef
   \csname#1tr@cetrue\endcsname{}\fi\next}
\def\m@ketag#1#2{\expandafter\n@wcount\csname#2tagno\endcsname
     \csname#2tagno\endcsname=0\let\tail=\all\xdef\all{\tail#2,}
   \ifx#1\l@c@l\let\tail=\r@s@t\xdef\r@s@t{\csname#2tagno\endcsname=0\tail}\fi
   \expandafter\gdef\csname#2cite\endcsname##1{\expandafter
     \ifx\csname#2tag##1\endcsname\relax?\else\csname#2tag##1\endcsname\fi
     \expandafter\ifx\csname#2tr@cetrue\endcsname\relax\else
     \write\cit@tionsout{#2tag ##1 cited on page \folio.}\fi}
   \expandafter\gdef\csname#2page\endcsname##1{\expandafter
     \ifx\csname#2page##1\endcsname\relax?\else\csname#2page##1\endcsname\fi
     \expandafter\ifx\csname#2tr@cetrue\endcsname\relax\else
     \write\cit@tionsout{#2tag ##1 cited on page \folio.}\fi}
   \expandafter\gdef\csname#2tag\endcsname##1{\expandafter
      \ifx\csname#2check##1\endcsname\relax
      \expandafter\xdef\csname#2check##1\endcsname{}%
      \else\immediate\write16{Warning: #2tag ##1 used more than once.}\fi
      \multit@g{#1}{#2}##1/X%
      \write\t@gsout{#2tag ##1 assigned number \csname#2tag##1\endcsname\space
      on page \number\count0.}%
   \csname#2tag##1\endcsname}}

\def\multit@g#1#2#3/#4X{\def\t@mp{#4}\ifx\t@mp\empty%
      \global\advance\csname#2tagno\endcsname by 1 
      \expandafter\xdef\csname#2tag#3\endcsname
      {#1\number\csname#2tagno\endcsnameX}%
   \else\expandafter\ifx\csname#2last#3\endcsname\relax
      \expandafter\n@wcount\csname#2last#3\endcsname
      \global\advance\csname#2tagno\endcsname by 1 
      \expandafter\xdef\csname#2tag#3\endcsname
      {#1\number\csname#2tagno\endcsnameX}
      \write\t@gsout{#2tag #3 assigned number \csname#2tag#3\endcsname\space
      on page \number\count0.}\fi
   \global\advance\csname#2last#3\endcsname by 1
   \def\t@mp{\expandafter\xdef\csname#2tag#3/}%
   \expandafter\t@mp\@mputate#4\endcsname
   {\csname#2tag#3\endcsname\lastpart{\csname#2last#3\endcsname}}\fi}
\def\t@gs#1{\def\all{}\m@ketag#1e\m@ketag#1s\m@ketag\t@d@l p
\let\realscite\scite
\let\realstag\stag
   \m@ketag\gl@b@l r \n@wread\t@gsin
   \openin\t@gsin=\jobname.tgs \re@der \closein\t@gsin
   \n@wwrite\t@gsout\openout\t@gsout=\jobname.tgs }
\outer\def\localtags{\t@gs\l@c@l}
\outer\def\globaltags{\t@gs\gl@b@l}
\outer\def\newlocaltag#1{\m@ketag\l@c@l{#1}}
\outer\def\newglobaltag#1{\m@ketag\gl@b@l{#1}}

\newif\ifpr@ 
\def\m@kecs #1tag #2 assigned number #3 on page #4.%
   {\expandafter\gdef\csname#1tag#2\endcsname{#3}
   \expandafter\gdef\csname#1page#2\endcsname{#4}
   \ifpr@\expandafter\xdef\csname#1check#2\endcsname{}\fi}
\def\re@der{\ifeof\t@gsin\let\next=\relax\else
   \read\t@gsin to\t@gline\ifx\t@gline\v@idline\else
   \expandafter\m@kecs \t@gline\fi\let \next=\re@der\fi\next}
\def\pretags#1{\pr@true\pret@gs#1,,}
\def\pret@gs#1,{\def\@{#1}\ifx\@\empty\let\n@xtfile=\relax
   \else\let\n@xtfile=\pret@gs \openin\t@gsin=#1.tgs \message{#1} \re@der 
   \closein\t@gsin\fi \n@xtfile}

\newcount\sectno\sectno=0\newcount\subsectno\subsectno=0
\newif\ifultr@local \def\ultralocal{\ultr@localtrue}
\def\firstpart{\number\sectno}
\def\lastpart#1{\ifcase#1 \or a\or b\or c\or d\or e\or f\or g\or h\or 
   i\or k\or l\or m\or n\or o\or p\or q\or r\or s\or t\or u\or v\or w\or 
   x\or y\or z \fi}

\def\resetall{\global\advance\sectno by 1\subsectno=0
   \gdef\firstpart{\number\sectno}\r@s@t}
\def\resetsub{\global\advance\subsectno by 1
   \gdef\firstpart{\number\sectno.\number\subsectno}\r@s@t}
\def\newsection#1\par{\resetall\vskip0pt plus.3\vsize\penalty-250
   \vskip0pt plus-.3\vsize\bigskip\bigskip
   \message{#1}\leftline{\bf#1}\nobreak\bigskip}
\def\subsection#1\par{\ifultr@local\resetsub\fi
   \vskip0pt plus.2\vsize\penalty-250\vskip0pt plus-.2\vsize
   \bigskip\smallskip\message{#1}\leftline{\bf#1}\nobreak\medskip}


\newdimen\marginshift

\newdimen\margindelta
\newdimen\marginmax
\newdimen\marginmin

\def\margininit{       
\marginmax=3 true cm                  
				      
\margindelta=0.1 true cm              
\marginmin=0.1true cm                 
\marginshift=\marginmin
}    

\def\t@gsjj#1,{\def\@{#1}\ifx\@\empty\let\next=\relax\else\let\next=\t@gsjj
   \def\@@{p}\ifx\@\@@\else
   \expandafter\gdef\csname#1cite\endcsname##1{\citejj{##1}}
   \expandafter\gdef\csname#1page\endcsname##1{?}
   \expandafter\gdef\csname#1tag\endcsname##1{\tagjj{##1}}\fi\fi\next}
\newif\ifshowstuffinmargin
\showstuffinmarginfalse
\def\jjtags{\ifx\shlhetal\relax 
  \else
\ifx\shlhetal\undefinedcontrolseq
\else
\showstuffinmargintrue
\ifx\all\relax\else\expandafter\t@gsjj\all,\fi\fi \fi
}

\def\tagjj#1{\realstag{#1}\oldmginpar{\zeigen{#1}}}
\def\citejj#1{\rechnen{#1}\mginpar{\zeigen{#1}}}     

\def\rechnen#1{\expandafter\ifx\csname stag#1\endcsname\relax ??\else
                           \csname stag#1\endcsname\fi}

\newdimen\theight

\def\marginfont{\sevenrm}

\def\trymarginbox#1{\setbox0=\hbox{\marginfont\hskip\marginshift #1}%
		\global\marginshift\wd0 
		\global\advance\marginshift\margindelta}

\def \oldmginpar#1{%
\ifvmode\setbox0\hbox to \hsize{\hfill\rlap{\marginfont\quad#1}}%
\ht0 0cm
\dp0 0cm
\box0\vskip-\baselineskip
\else 
             \vadjust{\trymarginbox{#1}%
		\ifdim\marginshift>\marginmax \global\marginshift\marginmin
			\trymarginbox{#1}%
                \fi
             \theight=\ht0
             \advance\theight by \dp0    \advance\theight by \lineskip
             \kern -\theight \vbox to \theight{\rightline{\rlap{\box0}}%
\vss}}\fi}

\newdimen\upordown
\global\upordown=8pt
\font\tinyfont=cmtt8 
\def\mginpar#1{\smash{\hbox to 0cm{\kern-10pt\raise7pt\hbox{\tinyfont #1}\hss}}}
\def\mginpar#1{{\hbox to 0cm{\kern-10pt\raise\upordown\hbox{\tinyfont #1}\hss}}\global\upordown-\upordown}


\def\t@gsoff#1,{\def\@{#1}\ifx\@\empty\let\next=\relax\else\let\next=\t@gsoff
   \def\@@{p}\ifx\@\@@\else
   \expandafter\gdef\csname#1cite\endcsname##1{\zeigen{##1}}
   \expandafter\gdef\csname#1page\endcsname##1{?}
   \expandafter\gdef\csname#1tag\endcsname##1{\zeigen{##1}}\fi\fi\next}
\def\verbatimtags{\showstuffinmarginfalse
\ifx\all\relax\else\expandafter\t@gsoff\all,\fi}
\def\zeigen#1{\hbox{$\scriptstyle\langle$}#1\hbox{$\scriptstyle\rangle$}}


\def\margintag#1{\ifshowstuffinmargin\oldmginpar{\zeigen{#1}}\fi}

\def\(#1){\edef\dot@g{\ifmmode\ifinner(\hbox{\noexpand\etag{#1}})
   \else\noexpand\eqno(\hbox{\noexpand\etag{#1}})\fi
   \else(\noexpand\ecite{#1})\fi}\dot@g}

\newif\ifbr@ck
\def\eat#1{}
\def\[#1]{\br@cktrue[\br@cket#1'X]}
\def\br@cket#1'#2X{\def\temp{#2}\ifx\temp\empty\let\next\eat
   \else\let\next\br@cket\fi
   \ifbr@ck\br@ckfalse\br@ck@t#1,X\else\br@cktrue#1\fi\next#2X}
\def\br@ck@t#1,#2X{\def\temp{#2}\ifx\temp\empty\let\neext\eat
   \else\let\neext\br@ck@t\def\temp{,}\fi
   \def\teemp{#1}\ifx\teemp\empty\else\rcite{#1}\fi\temp\neext#2X}
\def\resetbr@cket{\gdef\[##1]{[\rtag{##1}]}}
\def\references{\resetbr@cket\newsection References\par}

\newtoks\symb@ls\newtoks\s@mb@ls\newtoks\p@gelist\n@wcount\ftn@mber
    \ftn@mber=1\newif\ifftn@mbers\ftn@mbersfalse\newif\ifbyp@ge\byp@gefalse
\def\defm@rk{\ifftn@mbers\n@mberm@rk\else\symb@lm@rk\fi}
\def\n@mberm@rk{\xdef\m@rk{{\the\ftn@mber}}%
    \global\advance\ftn@mber by 1 }
\def\rot@te#1{\let\temp=#1\global#1=\expandafter\r@t@te\the\temp,X}
\def\r@t@te#1,#2X{{#2#1}\xdef\m@rk{{#1}}}
\def\b@@st#1{{$^{#1}$}}\def\str@p#1{#1}
\def\symb@lm@rk{\ifbyp@ge\rot@te\p@gelist\ifnum\expandafter\str@p\m@rk=1 
    \s@mb@ls=\symb@ls\fi\write\f@nsout{\number\count0}\fi \rot@te\s@mb@ls}
\def\byp@ge{\byp@getrue\n@wwrite\f@nsin\openin\f@nsin=\jobname.fns 
    \n@wcount\currentp@ge\currentp@ge=0\p@gelist={0}
    \re@dfns\closein\f@nsin\rot@te\p@gelist
    \n@wread\f@nsout\openout\f@nsout=\jobname.fns }
\def\m@kelist#1X#2{{#1,#2}}
\def\re@dfns{\ifeof\f@nsin\let\next=\relax\else\read\f@nsin to \f@nline
    \ifx\f@nline\v@idline\else\let\t@mplist=\p@gelist
    \ifnum\currentp@ge=\f@nline
    \global\p@gelist=\expandafter\m@kelist\the\t@mplistX0
    \else\currentp@ge=\f@nline
    \global\p@gelist=\expandafter\m@kelist\the\t@mplistX1\fi\fi
    \let\next=\re@dfns\fi\next}
\def\symbols#1{\symb@ls={#1}\s@mb@ls=\symb@ls} 
\def\bigsymbol{\textstyle}
\symbols{\bigsymbol\ast,\dagger,\ddagger,\sharp,\flat,\natural,\star}
\def\ftnumbers{\ftn@mberstrue} \def\ftsymbols{\ftn@mbersfalse}
\def\paginal{\byp@ge} \def\resetftnumbers{\ftn@mber=1}
\def\ftnote#1{\defm@rk\expandafter\expandafter\expandafter\footnote
    \expandafter\b@@st\m@rk{#1}}

\long\def\jump#1\endjump{}
\def\ssum{\mathop{\lower .1em\hbox{$\textstyle\Sigma$}}\nolimits}

\def\qed{\nobreak\kern 1em \vrule height .5em width .5em depth 0em}
\def\newneq{\hbox{\rlap{\hbox to 1\wd9{\hss$=$\hss}}\raise .1em 
   \hbox to 1\wd9{\hss$\scriptscriptstyle/$\hss}}}
\def\subsetne{\setbox9 = \hbox{$\subset$}\mathrel{\hbox{\rlap
   {\lower .4em \newneq}\raise .13em \hbox{$\subset$}}}}
\def\supsetne{\setbox9 = \hbox{$\subset$}\mathrel{\hbox{\rlap
   {\lower .4em \newneq}\raise .13em \hbox{$\supset$}}}}

\def\vbar{\mathchoice{\vrule height6.3ptdepth-.5ptwidth.8pt\kern-.8pt}
   {\vrule height6.3ptdepth-.5ptwidth.8pt\kern-.8pt}
   {\vrule height4.1ptdepth-.35ptwidth.6pt\kern-.6pt}
   {\vrule height3.1ptdepth-.25ptwidth.5pt\kern-.5pt}}
\def\f@dge{\mathchoice{}{}{\mkern.5mu}{\mkern.8mu}}
\def\b@c#1#2{{\rm \mkern#2mu\vbar\mkern-#2mu#1}}
\def\b@b#1{{\rm I\mkern-3.5mu #1}}
\def\b@a#1#2{{\rm #1\mkern-#2mu\f@dge #1}}
\def\bb#1{{\count4=`#1 \advance\count4by-64 \ifcase\count4\or\b@a A{11.5}\or
   \b@b B\or\b@c C{5}\or\b@b D\or\b@b E\or\b@b F \or\b@c G{5}\or\b@b H\or
   \b@b I\or\b@c J{3}\or\b@b K\or\b@b L \or\b@b M\or\b@b N\or\b@c O{5} \or
   \b@b P\or\b@c Q{5}\or\b@b R\or\b@a S{8}\or\b@a T{10.5}\or\b@c U{5}\or
   \b@a V{12}\or\b@a W{16.5}\or\b@a X{11}\or\b@a Y{11.7}\or\b@a Z{7.5}\fi}}

\catcode`\X=11 \catcode`\@=12




\let\thischap\jobname

\def\partof#1{\csname returnthe#1part\endcsname}
\def\chapof#1{\csname returnthe#1chap\endcsname}

\def\setchapter#1,#2,#3;{%
  \expandafter\def\csname returnthe#1part\endcsname{#2}%
  \expandafter\def\csname returnthe#1chap\endcsname{#3}%
}

\setchapter 300a,A,II.A;
\setchapter 300b,A,II.B;
\setchapter 300c,A,II.C;
\setchapter 300d,A,II.D;
\setchapter 300e,A,II.E;
\setchapter 300f,A,II.F;
\setchapter 300g,A,II.G;
\setchapter  E53,B,N;
\setchapter  88r,B,I;
\setchapter  600,B,III;
\setchapter  705,B,IV;
\setchapter  734,B,V;

\def\cprefix#1{
\edef\theotherpart{\partof{#1}}\edef\theotherchap{\chapof{#1}}%
\ifx\theotherpart\thispart
   \ifx\theotherchap\thischap 
    \else 
     \theotherchap%
    \fi
   \else 
     \theotherchap\fi}

\def\sectioncite[#1]#2{%
     \cprefix{#2}#1}

\edef\thispart{\partof{\thischap}}
\edef\thischap{\chapof{\thischap}}

\def\lastpage of '#1' is #2.{\expandafter\def\csname lastpage#1\endcsname{#2}}


\def\spuriousreset{}


\expandafter\ifx\csname citeadd.tex\endcsname\relax
\expandafter\gdef\csname citeadd.tex\endcsname{}
\else \message{Hey!  Apparently you were trying to
\string\input{citeadd.tex} twice.   This does not make sense.} 
\errmessage{Please edit your file (probably \jobname.tex) and remove
any duplicate ``\string\input'' lines}\endinput\fi

\sectno=-1   
\localtags
\jjtags
\NoBlackBoxes
\define\mr{\medskip\roster}
\define\sn{\smallskip\noindent}
\define\mn{\medskip\noindent}
\define\bn{\bigskip\noindent}
\define\ub{\underbar}
\define\wilog{\text{without loss of generality}}
\define\ermn{\endroster\medskip\noindent}
\define\dbca{\dsize\bigcap}

\define \nl{\newline}
\magnification=\magstep 1
\documentstyle{amsppt}

{    
\catcode`@11

\ifx\alicetwothousandloaded@\relax
  \endinput\else\global\let\alicetwothousandloaded@\relax\fi

\gdef\subjclass{\let\savedef@\subjclass
 \def\subjclass##1\endsubjclass{\let\subjclass\savedef@
   \toks@{\def\usualspace{{\rm\enspace}}\eightpoint}%
   \toks@@{##1\unskip.}%
   \edef\thesubjclass@{\the\toks@
     \frills@{{\noexpand\rm2000 {\noexpand\it Mathematics Subject
       Classification}.\noexpand\enspace}}%
     \the\toks@@}}%
  \nofrillscheck\subjclass}
} 


\expandafter\ifx\csname alice2jlem.tex\endcsname\relax
  \expandafter\xdef\csname alice2jlem.tex\endcsname{\the\catcode`@}
\else \message{Hey!  Apparently you were trying to
\string\input{alice2jlem.tex}  twice.   This does not make sense.}
\errmessage{Please edit your file (probably \jobname.tex) and remove
any duplicate ``\string\input'' lines}\endinput\fi

\expandafter\ifx\csname bib4plain.tex\endcsname\relax
  \expandafter\gdef\csname bib4plain.tex\endcsname{}
\else \message{Hey!  Apparently you were trying to \string\input
  bib4plain.tex twice.   This does not make sense.}
\errmessage{Please edit your file (probably \jobname.tex) and remove
any duplicate ``\string\input'' lines}\endinput\fi

\def\renewcommand{\newcommand}	       
\edef\cite{\the\catcode`@}%
\catcode`@ = 11
\let\@oldatcatcode = \cite
\chardef\@letter = 11
\chardef\@other = 12
%
%
%
%
\def\@innerdef#1#2{\edef#1{\expandafter\noexpand\csname #2\endcsname}}%
%
%
\@innerdef\@innernewcount{newcount}%
\@innerdef\@innernewdimen{newdimen}%
\@innerdef\@innernewif{newif}%
\@innerdef\@innernewwrite{newwrite}%
%
%
%
\def\@gobble#1{}%
%
%
%
\ifx\inputlineno\@undefined
   \let\@linenumber = \empty 
\else
   \def\@linenumber{\the\inputlineno:\space}%
\fi
%
%
%
\def\@futurenonspacelet#1{\def\cs{#1}%
   \afterassignment\@stepone\let\@nexttoken=
}%
\begingroup 
\def\\{\global\let\@stoken= }%
\\ 
\endgroup
\def\@stepone{\expandafter\futurelet\cs\@steptwo}%
\def\@steptwo{\expandafter\ifx\cs\@stoken\let\@@next=\@stepthree
   \else\let\@@next=\@nexttoken\fi \@@next}%
\def\@stepthree{\afterassignment\@stepone\let\@@next= }%
%
%
%
\def\@getoptionalarg#1{%
   \let\@optionaltemp = #1%
   \let\@optionalnext = \relax
   \@futurenonspacelet\@optionalnext\@bracketcheck
}%
%
%
\def\@bracketcheck{%
   \ifx [\@optionalnext
      \expandafter\@@getoptionalarg
   \else
      \let\@optionalarg = \empty
      \expandafter\@optionaltemp
   \fi
}%
\def\@@getoptionalarg[#1]{%
   \def\@optionalarg{#1}%
   \@optionaltemp
}%
%
%
%
\def\@nnil{\@nil}%
\def\@fornoop#1\@@#2#3{}%
\def\@for#1:=#2\do#3{%
   \edef\@fortmp{#2}%
   \ifx\@fortmp\empty \else
      \expandafter\@forloop#2,\@nil,\@nil\@@#1{#3}%
   \fi
}%
\def\@forloop#1,#2,#3\@@#4#5{\def#4{#1}\ifx #4\@nnil \else
       #5\def#4{#2}\ifx #4\@nnil \else#5\@iforloop #3\@@#4{#5}\fi\fi
}%
\def\@iforloop#1,#2\@@#3#4{\def#3{#1}\ifx #3\@nnil
       \let\@nextwhile=\@fornoop \else
      #4\relax\let\@nextwhile=\@iforloop\fi\@nextwhile#2\@@#3{#4}%
}%
%
%
%
\@innernewif\if@fileexists
\def\@testfileexistence{\@getoptionalarg\@finishtestfileexistence}%
\def\@finishtestfileexistence#1{%
   \begingroup
      \def\extension{#1}%
      \immediate\openin0 =
         \ifx\@optionalarg\empty\jobname\else\@optionalarg\fi
         \ifx\extension\empty \else .#1\fi
         \space
      \ifeof 0
         \global\@fileexistsfalse
      \else
         \global\@fileexiststrue
      \fi
      \immediate\closein0
   \endgroup
}%
%
%
%
%
\def\bibliographystyle#1{%
   \@readauxfile
   \@writeaux{\string\bibstyle{#1}}%
}%
\let\bibstyle = \@gobble
%
%
\let\bblfilebasename = \jobname
\def\bibliography#1{%
   \@readauxfile
   \@writeaux{\string\bibdata{#1}}%
   \@testfileexistence[\bblfilebasename]{bbl}%
   \if@fileexists
      \nobreak
      \@readbblfile
   \fi
}%
\let\bibdata = \@gobble
%
%
\def\nocite#1{%
   \@readauxfile
   \@writeaux{\string\citation{#1}}%
}%
\@innernewif\if@notfirstcitation
%
%
\def\cite{\@getoptionalarg\@cite}%
%
%
\def\@cite#1{%
   \let\@citenotetext = \@optionalarg
   \printcitestart
   \nocite{#1}%
   \@notfirstcitationfalse
   \@for \@citation :=#1\do
   {%
      \expandafter\@onecitation\@citation\@@
   }%
   \ifx\empty\@citenotetext\else
      \printcitenote{\@citenotetext}%
   \fi
   \printcitefinish
}%
\newif\ifweareinprivate
\weareinprivatetrue
\ifx\shlhetal\undefinedcontrolseq\weareinprivatefalse\fi
\ifx\shlhetal\relax\weareinprivatefalse\fi
\def\@onecitation#1\@@{%
   \if@notfirstcitation
      \printbetweencitations
   \fi
   \expandafter \ifx \csname\@citelabel{#1}\endcsname \relax
      \if@citewarning
         \message{\@linenumber Undefined citation `#1'.}%
      \fi
     \ifweareinprivate
      \expandafter\gdef\csname\@citelabel{#1}\endcsname{%
\strut 
\vadjust{\vskip-\dp\strutbox
\vbox to 0pt{\vss\parindent0cm \leftskip=\hsize 
\advance\leftskip3mm
\advance\hsize 4cm\strut\openup-4pt 
\rightskip 0cm plus 1cm minus 0.5cm ?  #1 ?\strut}}
         {\tt
            \escapechar = -1
            \nobreak\hskip0pt\pfeilsw
            \expandafter\string\csname#1\endcsname
             \pfeilso
            \nobreak\hskip0pt
         }%
      }%
     \else  
      \expandafter\gdef\csname\@citelabel{#1}\endcsname{%
            {\tt\expandafter\string\csname#1\endcsname}
      }%
     \fi  
   \fi
   \csname\@citelabel{#1}\endcsname
   \@notfirstcitationtrue
}%
%
%
\def\@citelabel#1{b@#1}%
%
%
\def\@citedef#1#2{\expandafter\gdef\csname\@citelabel{#1}\endcsname{#2}}%
%
%
%
\def\@readbblfile{%
   \ifx\@itemnum\@undefined
      \@innernewcount\@itemnum
   \fi
   \begingroup
      \def\begin##1##2{%
         \setbox0 = \hbox{\biblabelcontents{##2}}%
         \biblabelwidth = \wd0
      }%
      \def\end##1{}
      %
      %
      \@itemnum = 0
      \def\bibitem{\@getoptionalarg\@bibitem}%
      \def\@bibitem{%
         \ifx\@optionalarg\empty
            \expandafter\@numberedbibitem
         \else
            \expandafter\@alphabibitem
         \fi
      }%
      \def\@alphabibitem##1{%
         \expandafter \xdef\csname\@citelabel{##1}\endcsname {\@optionalarg}%
         \ifx\biblabelprecontents\@undefined
            \let\biblabelprecontents = \relax
         \fi
         \ifx\biblabelpostcontents\@undefined
            \let\biblabelpostcontents = \hss
         \fi
         \@finishbibitem{##1}%
      }%
      \def\@numberedbibitem##1{%
         \advance\@itemnum by 1
         \expandafter \xdef\csname\@citelabel{##1}\endcsname{\number\@itemnum}%
         \ifx\biblabelprecontents\@undefined
            \let\biblabelprecontents = \hss
         \fi
         \ifx\biblabelpostcontents\@undefined
            \let\biblabelpostcontents = \relax
         \fi
         \@finishbibitem{##1}%
      }%
      \def\@finishbibitem##1{%
         \biblabelprint{\csname\@citelabel{##1}\endcsname}%
         \@writeaux{\string\@citedef{##1}{\csname\@citelabel{##1}\endcsname}}%
         \ignorespaces
      }%
      %
      %
      \let\em = \bblem
      \let\newblock = \bblnewblock
      \let\sc = \bblsc
      \frenchspacing
      \clubpenalty = 4000 \widowpenalty = 4000
      \tolerance = 10000 \hfuzz = .5pt
      \everypar = {\hangindent = \biblabelwidth
                      \advance\hangindent by \biblabelextraspace}%
      \bblrm
      \parskip = 1.5ex plus .5ex minus .5ex
      \biblabelextraspace = .5em
      \bblhook
      \input \bblfilebasename.bbl
   \endgroup
}%
%
%
\@innernewdimen\biblabelwidth
\@innernewdimen\biblabelextraspace
%
%
%
\def\biblabelprint#1{%
   \noindent
   \hbox to \biblabelwidth{%
      \biblabelprecontents
      \biblabelcontents{#1}%
      \biblabelpostcontents
   }%
   \kern\biblabelextraspace
}%
%
%
%
\def\biblabelcontents#1{{\bblrm [#1]}}%
%
%
\def\bblrm{\rm}%
%
%
\def\bblem{\it}%
%
%
\def\bblsc{\ifx\@scfont\@undefined
              \font\@scfont = cmcsc10
           \fi
           \@scfont
}%
%
%
\def\bblnewblock{\hskip .11em plus .33em minus .07em }%
%
%
\let\bblhook = \empty
%
%
%
\def\printcitestart{[}
\def\printcitefinish{]}
\def\printbetweencitations{, }
\def\printcitenote#1{, #1}
%
%
%
\let\citation = \@gobble
%
%
%
\@innernewcount\@numparams
%
%
\def\newcommand#1{%
   \def\@commandname{#1}%
   \@getoptionalarg\@continuenewcommand
}%
%
%
\def\@continuenewcommand{%
   \@numparams = \ifx\@optionalarg\empty 0\else\@optionalarg \fi \relax
   \@newcommand
}%
%
%
\def\@newcommand#1{%
   \def\@startdef{\expandafter\edef\@commandname}%
   \ifnum\@numparams=0
      \let\@paramdef = \empty
   \else
      \ifnum\@numparams>9
         \errmessage{\the\@numparams\space is too many parameters}%
      \else
         \ifnum\@numparams<0
            \errmessage{\the\@numparams\space is too few parameters}%
         \else
            \edef\@paramdef{%
               \ifcase\@numparams
                  \empty  No arguments.
               \or ####1%
               \or ####1####2%
               \or ####1####2####3%
               \or ####1####2####3####4%
               \or ####1####2####3####4####5%
               \or ####1####2####3####4####5####6%
               \or ####1####2####3####4####5####6####7%
               \or ####1####2####3####4####5####6####7####8%
               \or ####1####2####3####4####5####6####7####8####9%
               \fi
            }%
         \fi
      \fi
   \fi
   \expandafter\@startdef\@paramdef{#1}%
}%
%
%
%
%
\def\@readauxfile{%
   \if@auxfiledone \else 
      \global\@auxfiledonetrue
      \@testfileexistence{aux}%
      \if@fileexists
         \begingroup
            \endlinechar = -1
            \catcode`@ = 11
            \input \jobname.aux
         \endgroup
      \else
         \message{\@undefinedmessage}%
         \global\@citewarningfalse
      \fi
      \immediate\openout\@auxfile = \jobname.aux
   \fi
}%
%
%
\newif\if@auxfiledone
\ifx\noauxfile\@undefined \else \@auxfiledonetrue\fi
%
%
%
%
\@innernewwrite\@auxfile
\def\@writeaux#1{\ifx\noauxfile\@undefined \write\@auxfile{#1}\fi}%
%
%
%
\ifx\@undefinedmessage\@undefined
   \def\@undefinedmessage{No .aux file; I won't give you warnings about
                          undefined citations.}%
\fi
%
%
\@innernewif\if@citewarning
\ifx\noauxfile\@undefined \@citewarningtrue\fi
%
%
%
\catcode`@ = \@oldatcatcode

\def\pfeilso{\leavevmode
            \vrule width 1pt height9pt depth 0pt\relax
           \vrule width 1pt height8.7pt depth 0pt\relax
           \vrule width 1pt height8.3pt depth 0pt\relax
           \vrule width 1pt height8.0pt depth 0pt\relax
           \vrule width 1pt height7.7pt depth 0pt\relax
            \vrule width 1pt height7.3pt depth 0pt\relax
            \vrule width 1pt height7.0pt depth 0pt\relax
            \vrule width 1pt height6.7pt depth 0pt\relax
            \vrule width 1pt height6.3pt depth 0pt\relax
            \vrule width 1pt height6.0pt depth 0pt\relax
            \vrule width 1pt height5.7pt depth 0pt\relax
            \vrule width 1pt height5.3pt depth 0pt\relax
            \vrule width 1pt height5.0pt depth 0pt\relax
            \vrule width 1pt height4.7pt depth 0pt\relax
            \vrule width 1pt height4.3pt depth 0pt\relax
            \vrule width 1pt height4.0pt depth 0pt\relax
            \vrule width 1pt height3.7pt depth 0pt\relax
            \vrule width 1pt height3.3pt depth 0pt\relax
            \vrule width 1pt height3.0pt depth 0pt\relax
            \vrule width 1pt height2.7pt depth 0pt\relax
            \vrule width 1pt height2.3pt depth 0pt\relax
            \vrule width 1pt height2.0pt depth 0pt\relax
            \vrule width 1pt height1.7pt depth 0pt\relax
            \vrule width 1pt height1.3pt depth 0pt\relax
            \vrule width 1pt height1.0pt depth 0pt\relax
            \vrule width 1pt height0.7pt depth 0pt\relax
            \vrule width 1pt height0.3pt depth 0pt\relax}

\def\pfeilsw{ \leavevmode 
            \vrule width 1pt height0.3pt depth 0pt\relax
            \vrule width 1pt height0.7pt depth 0pt\relax
            \vrule width 1pt height1.0pt depth 0pt\relax
            \vrule width 1pt height1.3pt depth 0pt\relax
            \vrule width 1pt height1.7pt depth 0pt\relax
            \vrule width 1pt height2.0pt depth 0pt\relax
            \vrule width 1pt height2.3pt depth 0pt\relax
            \vrule width 1pt height2.7pt depth 0pt\relax
            \vrule width 1pt height3.0pt depth 0pt\relax
            \vrule width 1pt height3.3pt depth 0pt\relax
            \vrule width 1pt height3.7pt depth 0pt\relax
            \vrule width 1pt height4.0pt depth 0pt\relax
            \vrule width 1pt height4.3pt depth 0pt\relax
            \vrule width 1pt height4.7pt depth 0pt\relax
            \vrule width 1pt height5.0pt depth 0pt\relax
            \vrule width 1pt height5.3pt depth 0pt\relax
            \vrule width 1pt height5.7pt depth 0pt\relax
            \vrule width 1pt height6.0pt depth 0pt\relax
            \vrule width 1pt height6.3pt depth 0pt\relax
            \vrule width 1pt height6.7pt depth 0pt\relax
            \vrule width 1pt height7.0pt depth 0pt\relax
            \vrule width 1pt height7.3pt depth 0pt\relax
            \vrule width 1pt height7.7pt depth 0pt\relax
            \vrule width 1pt height8.0pt depth 0pt\relax
            \vrule width 1pt height8.3pt depth 0pt\relax
            \vrule width 1pt height8.7pt depth 0pt\relax
            \vrule width 1pt height9pt depth 0pt\relax
      }


\def\widestnumber#1#2{}

\def\citewarning#1{\ifx\shlhetal\relax 
    \else
    \par{#1}\par
    \fi
}

\def\rm{\fam0 \tenrm}

\def\fakesubhead#1\endsubhead{\bigskip\noindent{\bf#1}\par}



%
%
%

%

\font\textrsfs=rsfs10
\font\scriptrsfs=rsfs7
\font\scriptscriptrsfs=rsfs5

\newfam\rsfsfam
\textfont\rsfsfam=\textrsfs
\scriptfont\rsfsfam=\scriptrsfs
\scriptscriptfont\rsfsfam=\scriptscriptrsfs

\edef\oldcatcodeofat{\the\catcode`\@}
\catcode`\@11

\def\Cal@@#1{\noaccents@ \fam \rsfsfam #1}

\catcode`\@\oldcatcodeofat


\expandafter\ifx \csname margininit\endcsname \relax\else\margininit\fi

\long\def\red#1\endred{}
\long\def\green#1\endgreen{}
\long\def\blue#1\endblue{}
\long\def\private#1\endprivate{}

\def\endred{ \unmatched endred! }
\def\endgreen{ \unmatched endgreen! }
\def\endblue{ \unmatched endblue! }
\def\endprivate{ \unmatched endprivate! }

\ifx\latexcolors\undefinedcs\def\latexcolors{}\fi

\def\emptycs{}
\def\evaluatelatexcolors{%
        \ifx\latexcolors\emptycs\else
        \expandafter\xxevaluate\latexcolors\xxfertig\evaluatelatexcolors\fi}
\def\xxevaluate#1,#2\xxfertig{\setupthiscolor{#1}%
        \def\latexcolors{#2}}


\font\smallfont=cmsl7
\def\rutgerscolor{\ifmmode\else\endgraf\fi\smallfont
\advance\leftskip0.5cm\relax}
\def\setupthiscolor#1{\edef\tmptmpcs{\noexpand\bgroup\noexpand\rutgerscolor
\noexpand\def\noexpand\currentcolor{#1}%
\noexpand}%
\expandafter\let\csname#1\endcsname\tmptmpcs
\def\tmptmpcs{\checkColorUnmatched{#1}\popthecolor}
\expandafter\let\csname end#1\endcsname\tmptmpcs}

\def\checkColorUnmatched#1{\def\expectcolor{#1}%
    \ifx\expectcolor\currentcolor   
    \else \edef\failhere{\noexpand\tryingToClose '\currentcolor' with end\expectcolor}\failhere\fi}

\def\currentcolor{???}

\def\popthecolor{\ifmmode\else\endgraf\fi\egroup}

\expandafter\def\csname#1\endcsname{}

\evaluatelatexcolors

 \let\outerhead\head
 \def\head{\innerhead}
 \let\innerhead\outerhead

 \let\outersubhead\subhead
 \def\subhead{\innersubhead}
 \let\innersubhead\outersubhead

 \let\outersubsubhead\subsubhead
 \def\subsubhead{\innersubsubhead}
 \let\innersubsubhead\outersubsubhead

 \def\proclaim{\innerproclaim}
 \let\innerproclaim\outerproclaim

 %
 %
 %
 %

\def\demo#1{\medskip\noindent{\it #1.\/}}
\def\enddemo{\smallskip}

\def\remark#1{\medskip\noindent{\it #1.\/}}
\def\endremark{\smallskip}

\pageheight{8.5truein}
\topmatter
\title{$\aleph_n$-free abelian group with no non-zero homomorphism to
$\Bbb Z$ \\
 Sh883} \endtitle
\author {Saharon Shelah \thanks {\null\newline I would like to thank 
Alice Leonhardt for the beautiful typing. \null\newline
This research was supported by German-Israeli Foundation for
Scientific Research and Development \null\newline
 First Typed - 06/Mar/6 \null\newline
 Latest Revision - 06/Aug/7} \endthanks} \endauthor 


\affil{The Hebrew University of Jerusalem \\
Einstein Institute of Mathematics \\
Edmond J. Safra Campus, Givat Ram \\
Jerusalem 91904, Israel
 \medskip
 Department of Mathematics \\
 Hill Center-Busch Campus \\
  Rutgers, The State University of New Jersey \\
 110 Frelinghuysen Road \\
 Piscataway, NJ 08854-8019 USA} \endaffil

\abstract  We, for any natural $n$, construct an $\aleph_n$-free
abelian groups which have few homomorphisms to $\Bbb Z$.  For this we
use ``$\aleph_n$-free $(n+1)$-dimensional black boxes".
The method is relevant to
e.g. construction of $\aleph_n$-free abelian groups with a prescribed
endomorphism ring.
\endabstract
\endtopmatter
\document

\newpage

\head {Annotated Content} \endhead
 \spuriousreset
\bn
\S1 Constructing $\aleph_{k(*)+1}$-free Abelian group
\mr
\item "{${{}}$}"  [We introduce ``$\bold x$ is a combinatorial
$k(*)$-parameter.  We also give a short cut for getting only ``there
is a non-Whitehead $\aleph_{k(*)+1}$-free non-free abelian group'' only
(this is from \scite{6.1} on).  This is similar to \cite[\S5]{Sh:771},
so proofs are put in an appendix, except \scite{af.63}, note that
\scite{af.63}(3) really belongs to \S3.]
\endroster
\bn
\S2 Black boxes
\mr
\item "{${{}}$}"  [We prove that we have black boxes in this context,
see \scite{ya.7}; it is based on the silly black box.  Now
\scite{ya.14} belongs to the short cut.]
\endroster
\bn
\S3 Constructing abelian groups from combinatorial parameter
\mr
\item "{${{}}$}"  [For $\bold x \in K^{\text{cb}}_{k(*)+1}$ we define
a class ${\Cal G}_{\bold x}$ of abelian groups constructed from it and
a black box.  We prove they are all $\aleph_{k(*)+1}$-free of
cardinality $|\Gamma|^{\bold x} + \aleph_0$ and for some $G \in {\Cal
G}_{\bold x}$ satisfies Hom$(G,\Bbb Z) = \{0\}$.]
\endroster
\bn
\S4 Appendix 1
\mr
\item "{${{}}$}"  [We give the proofs from \cite{Sh:771} and with the
relevant changes.]
\endroster
\newpage

\head {\S0 Introduction} \endhead  \resetall \sectno=0
 \spuriousreset
\bigskip

For regular $\theta = \aleph_n$ we look for a $\theta$-free
abelian group $G$ with Hom$(G,\Bbb Z) = \{0\}$.  We first construct
$G$ and a subgroup $\Bbb Z z \subseteq G$ which
is not a direct summand.  If instead ``not direct product" we ask
``not free" so naturally of cardinality $\theta$, we know much: see
\cite{EM02}.

We can ask further questions on abelian groups, their endormorphism
rings, similarly on modules; naturally questions whose answer is
known when we demand $\aleph_1$-free instead $\aleph_n$-free; see
\cite{GbTl06} .  But we feel
those two cases can serve as a base
for them (or we can immitate the proofs).  Also this concentration is
reasonable for sorting out the set theoretical situation.  Why not $\theta
= \aleph_\omega$ and higher cardinals? (there are more reasonable
cardinals for which such results are not excluded), note that 
even in previous questions historically this was harder.

For $n=1$ we can use ${}^\omega \Bbb Z$ and $z = \langle
1,1,1,\ldots\rangle$.  But there is such an abelian group of
cardinality $\aleph_1$, by \cite[\S4]{Sh:98}.  However, if MA then
$\aleph_2 < 2^{\aleph_0} \Rightarrow$ any $\aleph_2$-free abelian group
of cardinality $< 2^{\aleph_0}$ fail the question.

The groups we construct are in a sense complete, like ${}^\omega 
\Bbb Z$.  They are essentially 
from \cite[\S5]{Sh:771} only there $S = \{0,1\}$ as there we are
interested in Borel abelian groups.  See earlier \cite{Sh:161}, see
representations of \cite{Sh:161} in \cite[\S3]{Sh:523}, \cite{EM02}.

However we still like to have $\theta = \aleph_\omega$,
i.e. $\aleph_\omega$-free abelian groups.   Concerning this we
continue in \cite{Sh:F691}.

We thank Ester Sternfield for corrections.
\bn
We shall use freely the well known theorem saying
\demo{\stag{hq.14} Theorem}   A subgroup of a free abelian group is a free
abelian group.
\enddemo
\bigskip

\definition{\stag{haq.7} Definition}  1) Pr$(\lambda,\kappa)$: means
that for some $\bar G$ we have:
\mr
\item "{$(a)$}"  $\bar G = \langle G_\alpha:\alpha \le \kappa +2
\rangle$
\sn
\item "{$(b)$}"  $\bar G$ is an increasing continuous sequence of free
abelian groups 
\sn
\item "{$(c)$}"  $|G_{\kappa +1}| \le \lambda$,
\sn
\item "{$(d)$}"  $G_{\kappa +1}/G_\alpha$ is free for 
$\alpha < \kappa$,
\sn
\item "{$(e)$}"  $G_0 = \{0\}$
\sn
\item "{$(f)$}"  $G_\beta/G_\alpha$ is free if $\alpha \le \beta \le
\kappa$
\sn
\item "{$(g)$}"   some $h \in \text{ Hom}(G_\kappa;\Bbb Z)$ cannot be
extended to $\hat h \in \text{ Hom}(G_\kappa,\Bbb Z)$.
\ermn
2) We write Pr$^-(\lambda,\theta,\kappa)$ be defined as above, only
   replacing ``$G_{\kappa +1}/G_\alpha$ is free for $\alpha < \kappa$"
   by ``$G_{\kappa +1}/G_\kappa$ is $\theta$-free.
\enddefinition
\newpage

\head {\S1 Constructing $\aleph_{k(*)+1}$-free abelian groups} \endhead  \resetall \sectno=1
 \spuriousreset
\bigskip

\definition{\stag{af.1} Definition}  1) We say $\bold x$ is a
combinatorial parameter if $\bold x =
(k,S,\Lambda) = (k^{\bold x},S^{\bold x},\Lambda^{\bold x})$ and 
they satisfy clauses (a)-(c)
\mr
\item"{$(a)$}"  $k < \omega$
\sn
\item "{$(b)$}"  $S$ is a set (in \cite{Sh:771}, $S = \{0,1\}$),
\sn
\item "{$(c)$}"  $\Lambda \subseteq {}^{k+1}({}^\omega S)$ and for
simplicity $|\Lambda| \ge \aleph_0$ if not said otherwise.
\ermn
1A) We say $\bold x$ is an abelian group $k$-parameter when 
$\bold x = (k,S,\Lambda,\bold a)$ such that (a),(b),(c) from part (1) and:
\mr
\item "{$(d)$}"  $\bold a$ is a function from $\Lambda \times \omega$ to
$\Bbb Z$.
\ermn
1B) Let $\bold x = (k^{\bold x},S^{\bold x},\Lambda^{\bold x})$ or
$\bold x = (k^{\bold x},S^{\bold x},\Lambda^{\bold x},\bold
a^{\bold x})$.  A parameter is a $k$-parameter for some $k$ and 
$K^{\text{cb}}_{k(*)}/\bold K^{\text{gr}}_{k(*)}$ is the class of
combinatorial/abelian group  $k(*)$-parameters.  
\nl
2) If $\bold x$ is an abelian group parameter and $\Lambda \subseteq
\Lambda^{\bold x}$ \ub{then} $\bold x \restriction \Lambda = (k(*)^{\bold x},
S^{\bold x},\Lambda,\bold a^{\bold x} \restriction (\Lambda \times \omega))$.
\nl
3) We may write $\bold a^{\bold x}_{\bar \eta,n}$ instead 
$\bold a^{\bold x}(\eta,n)$.  Let $w_{k,m} = w(k,m) = \{\ell \le k:\ell
\ne m\}$.
\nl
4) We say $\bold x$ is full when $\Lambda^{\bold x} =
   {}^{k(*)}({}^\omega S)$.
\nl
5) If $\Lambda \subseteq \Lambda^{\bold x}$ let $\bold x \restriction
\Lambda$ be $(k^{\bold x},S^{\bold x},\Lambda)$ or $(k^{\bold x},
S^{\bold x},\Lambda,\bold a \restriction \Lambda)$ as suitable.
We may write $\bold x = (\bold y,\bold a)$ if $\bold a = \bold
a^{\bold x},\bold y = (k^{\bold x},S^{\bold x},\Lambda^{\bold x})$.
\enddefinition
\bigskip

\demo{\stag{af.3} Convention}  If $\bold x$ is clear from the context
we may write $k$ or $k(*),S,\Lambda,\bold a$ instead of $k^{\bold x},
X^{\bold s},\Lambda^{\bold x},\bold a^{\bold x}$.
\enddemo
\bn
A variant of the above is
\definition{\stag{ya.77} Definition}  1) For $\bar S = \langle S_n:m
\le k\rangle$ we define when $\bold x$ is a
$\bar S$-parameter: $\bar \eta \in \Lambda^{\bold x} \wedge m \le
k^{\bold x} \Rightarrow \eta_m \in {}^\omega(S_m)$.
\nl
2) We say $\bar \alpha$ is a $(\bold x,\bar \chi)$-black box or
   Qr$(\bold x,\bar \chi)$ when:
\mr
\item "{$(a)$}"  $\bar \chi = \langle \chi_m:m \le k^{\bold x}\rangle$
\sn
\item "{$(b)$}"  $\bar \alpha = \langle \bar \alpha_{\bar \eta}:\bar
\eta \in \Lambda^{\bold x}\rangle$
\sn
\item "{$(c)$}"  $\bar \alpha_\eta = \langle \alpha_{\bar \eta,m,n}:m
\le k^{\bold x},n < \omega\rangle$ and $\alpha_{\eta,m,n} < \chi_m$
\sn
\item "{$(d)$}"  if $h_m:\Lambda^{\bold x}_m \rightarrow \chi_m$ for
$m \le k^{\bold x}$ then for some $\bar \eta \in \Lambda^{\bold x}$ we
have: $m \le k^{\bold x} \wedge n < w \Rightarrow h(\bar \eta
\upharpoonleft \langle m,n \rangle) = \alpha_{\bar \eta,m,n}$, see
Definition \scite{af.4} below on ``$\bar \eta \upharpoonleft \langle
m,n\rangle$. 
\ermn
2A) We may replace $\bar \chi$ by $\chi$ if $\bar \chi =\langle
\chi_\ell:\ell \le k^{\bold x}\rangle$.  We may replace $\bold x$ by
$\Lambda^{\bold x}$ (so say Qr$(\Lambda^{\bold x},\bar \chi)$ or say
$\bar \alpha$ is a $(\Lambda,\bar \chi)$-black box).
\nl 
3) We say a $\bar S$-parameter $\bold x$ is full when 
$\Lambda^{\bold x} = \dsize \prod_{m \le k} {}^\omega(S_m)$.
\enddefinition
\bigskip

\definition{\stag{af.4} Definition}  For an $k(*)$-parameter $\bold x$
and for $m \le k(*)$ let
\mr
\item "{$(a)$}"  $\Lambda^{\bold x}_m = \Lambda_{\bold x,m} = 
\{\bar \eta:\bar \eta = \langle \eta_\ell:\ell \le k(*)\rangle$ 
and $\eta_m \in {}^{\omega >}S \text{ and } \ell \le k(*) 
\wedge \ell \ne m \Rightarrow \eta_\ell \in {}^\omega S$ and
for some $\bar \eta' \in \Lambda$ we have $n < \omega,\bar \eta = \bar \eta'
\upharpoonleft (m,n)\}$ where
\nl
$\bar \eta = \bar \eta' \upharpoonleft \langle m,n\rangle$ means 
$\eta_m = \eta'_m \restriction n$
and $\ell \le k(*) \wedge \ell \ne m \Rightarrow \eta_\ell =
\eta'_\ell\}$
\sn
\item "{$(b)$}"  $\Lambda^{\bold x}_{\le k(*)}$ is 
$\cup\{\Lambda^{\bold x}_m:m \le k(*)\}$
\sn
\item "{$(c)$}"  $\bold m(\bar \eta) = m$ if $\bar \eta \in
\Lambda^{\bold x}_m$.
\endroster
\enddefinition
\bigskip

\definition{\stag{af.5} Definition}  1) We say a combinatorial
$k(*)$-parameter $\bold x$ is free when there is a list $\langle \bar
\eta^\alpha:\alpha < \alpha(*)\rangle$ of $\Lambda^{\bold x}$ such
that for every $\alpha$ for some $m \le k(*)$ and $n < \omega$ we have
\mr
\item "{$(*)$}"  $\bar \eta^\alpha_m \upharpoonleft \langle m,n\rangle \notin
\{\eta^\beta_m \upharpoonleft \langle m,n\rangle:\beta < \alpha\}$.
\ermn
2) We say a combinatorial $k$-parameter $\bold x$ is $\theta$-free
when $(k,S^{\bold x},\Lambda)$ is free for every $\Lambda \subseteq
\Lambda^{\bold x}$ of cardinality $< \theta$.
\enddefinition
\bigskip

\remark{Remark}  1) We can require in $(*)$ even $(\exists^\infty
n)[\eta^\alpha_m(n) \notin \cup\{\eta^\beta_\ell(n'):\ell \le k,\beta
< \alpha,n < \omega\}]$.

At present this seems an immaterial change.
\endremark
\bigskip

\definition{\stag{6.1} Definition}  For $k(*) < \omega$ and
$k(*)$-parameter $\bold x$ we define
an abelian group $G = G_{\bold x}$ as follows: it is generated
by $\{x_{\bar \eta}:m \le k(*) \text{ and } \bar \eta \in
\Lambda^{\bold x}_m\} \cup \{y_{\bar \eta,n}:n < \omega$ and
$\bar \eta \in \Lambda^{\bold x}\} \cup \{z\}$ freely except the equations:
\mr
\item "{$\boxtimes_{\bar \eta,n}$}"  $(n!)y_{\bar \eta,n+1} = y_{\bar \eta,n}
+ \bold a^{\bold x}_{\bar \eta,n} z + 
\sum\{x_{\bar \eta \upharpoonleft <m,n>}:m \le k(*)\}$.
\ermn
(Note that if $m_1 < m_2 \le k(*)$ then $\bar\eta_{m_1} \ne \bar
\eta_{m_2}$ having different index sets).
\enddefinition
\bigskip

\demo{\stag{af.21} Explanation}   A canonical example of a non-free group is
$(\Bbb Q,+)$.  Other examples are related to it after we divide by
something.  The $y$'s here play that role of provided (hidden) copies
of $\Bbb Q$.  What about $x$'s?  We use
$m \le k(*)$ to give $\langle y_{\bar \eta,n}:n < \omega \rangle,k(*)$
``chances", ``opportunities" to avoid having $(\Bbb Q,+)$ as 
a quotient, one for each cardinal
$\le \aleph_{k(*)}$.  More specifically, for each $m(*) \le k(*)$ if 
$H \subseteq G$ is the
subgroup which is generated by $X = \{x_{\bar \eta \upharpoonleft <m,n>}:m \ne
m(*)$ and $\bar\eta \in {}^{k(*)+1}({}^\omega S)$ and $m \le k(*)\}$, 
still in $G/H$ the set $\{y_{\bar \eta,n}:n < \omega\}$ does not
generate a copy of $\Bbb Q$, as witnessed by 
$\{x_{\bar \eta \upharpoonleft <m(*),n>}:n < \omega\}$.
\enddemo
\bn
As a warm up we note:
\proclaim{\stag{af.28} Claim}  For $k(*) < \omega$ and $k(*)$-parameter
$\bold x$ the abelian group $G_{\bold x}$ is an $\aleph_1$-free abelian group.
\endproclaim
\bn
Now systematically
\definition{\stag{af.35} Definition}  Let $\bold x$ be a
$k(*)$-parameter.
\nl
1) For $U \subseteq {}^\omega S$ let $G_U = G^{\bold x}_U$ be
the subgroup of $G$ generated by $Y_U = Y^{\bold x}_U = \{z\}
= \{y_{\bar \eta,n}:\bar \eta \in \Lambda \cap {}^{k(*)+1}(U)$ and 
$n < \omega\} \cup
\{x_{\eta \upharpoonleft <m,n>}:m \le k(*)$ and $\bar \eta \in
{}^{(k(*)+1)}(U)$ and $n < \omega\}$.  Let
$G^+_U = G^{\bold x,+}_U$ be the divisible hull of $G_U$ and 
$G^+ = G^+_{({}^\omega S)}$.
\nl
2) For $U \subseteq {}^\omega S$ and finite $u \subseteq {}^\omega S$
let $G_{U,u}$ be the subgroup \footnote{note that if $u=\{\eta\}$ then
$G_{U,u} = G_U$}  of 
$G$ generated by $\cup\{G_{U \cup (u \backslash \{\eta\})}:
\eta \in u\}$; and for $\bar \eta \in {}^{k(*)
\ge} U$ let $G_{U,\bar \eta}$ be the
subgroup of $G$ generated by $\cup \{G_{U \cup \{\eta_k:k < \ell
g(\bar \eta) \text{ and } k\ne \ell\}}:\ell < \ell g(\bar \eta)\}$.
\nl
3) For $U \subseteq {}^\omega S$ let $\Xi_U = \Xi^{\bold x}_U = 
\{\text{the equation }
\boxtimes_{\bar \eta,n}:\bar \eta \in \Lambda \cap {}^{k(*)+1} U$ and $n <
\omega\}$.  Let $\Xi_{U,u} = \Xi^{\bold x}_{U,u} = 
\cup\{\Xi_{U \cup u \backslash
\{\beta\}}:\beta \in u\}$.
\enddefinition
\bigskip

\proclaim{\stag{af.42} Claim}   Let $\bold x \in \bold K_{k(*)}$.
\nl
0) If $U_1 \subseteq U_2 \subseteq
{}^\omega S$ \ub{then} $G^+_{U_1} \subseteq G^+_{U_2}
\subseteq G^+$.
\nl
1) For any $n(*) < \omega$, the abelian group
$G^+_U$ (which is a vector space over $\Bbb Q$), has the basis
$Y^{n(*)}_{U_i} := \{z\} \cup
\{y_{\bar \eta,n(*)}:\bar \eta \in \Lambda \cap {}^{k(*)+1}(U)\} 
\cup \{x_{\bar \eta \upharpoonleft <m,n>}:m \le k(*),\bar \eta \in
{}^{k(*)+1}(U)$ and $n < \omega\}$. 
\nl
2) For $U \subseteq {}^\omega S$ the abelian group
$G_U$ is generated by $Y_U$ freely (as an abelian group) except the
set $\Xi_U$ of equations.
\nl
3) If $m(*) < \omega$ and 
$U_m \subseteq {}^\omega S$ for $m < m(*)$ \ub{then} the
subgroup $G_{U_0} + \ldots + G_{U_{m(*)-1}}$ of $G$ is 
generated by $Y_{U_0} \cup Y_{U_1} \cup
\ldots \cup Y_{U_{m(*)-1}}$ freely (as an abelian group) except the
equations in $\Xi_{U_0} \cup \Xi_{U_1} \cup \ldots \cup \Xi_{U_{m(*)-1}}$
provided that
\mr
\item "{$\circledast$}"  if $\eta_0,\dotsc,\eta_{k(*)} \in
\cup\{U_m:m < m(*)\}$ are such that 
\nl

$(\forall \ell \le k(*))
(\exists m < m(*))[\{\eta_0,\dotsc,\eta_{k(*)}\}
\backslash \{\eta_\ell\} \subseteq U_m)$ 
\nl

\ub{then} for some $m < m(*)$
we have $\{\eta_0,\dotsc,\eta_{k(*)}\} \subseteq U_m$.
\ermn
4) If $m(*) \le k(*)$ and $U_\ell 
= U \backslash U'_\ell$ for $\ell < m(*)$ and
$\langle U'_\ell:\ell < m(*)\rangle$ are pairwise disjoint \ub{then}
$\circledast$ holds.
\nl
5) $G_{U,u} \subseteq G_{U \cup u}$ if $U \subseteq {}^\omega S$ and
$u \subseteq {}^\omega S \backslash U$ is finite; moreover $G_{U,u}
\subseteq_{\text{pr}} G_{U \cup u} \subseteq_{\text{pr}} G$.
\nl
6) If $\langle U_\alpha:\alpha < \alpha(*)\rangle$ is
 $\subseteq$-increasing continuous \ub{then} also $\langle
 G_{U_\alpha}:\alpha < \alpha(*)\rangle$ is $\subseteq$-increasing
   continuous.
\nl
7) If $U_1 \subseteq U_2 \subseteq U \subseteq {}^\omega S$ and $u \subseteq
{}^\omega S \backslash U$ is finite, $|u| < k(*)$ and $U_2
\backslash U_1 = \{\eta\}$ and $v=u \cup\{\eta\}$ \ub{then}
$(G_{U,u} + G_{U_2 \cup u})/(G_{U,u} + G_{U_1 \cup u})$ is isomorphic to
$G_{U_1 \cup v}/G_{U_1,v}$.
\nl
8) If $U \subseteq {}^\omega S$ and $u \subseteq {}^\omega S
   \backslash U$ has $\le k(*)$ members \ub{then} $(G_{U,u} + G_u)/G_{U,u}$
   is isomorphic to $G_u/G_{\emptyset,u}$.
\endproclaim
\bigskip
\bn
\margintag{af.49}\ub{\stag{af.49} Discussion}:  For the reader 
we write what the group $G_{\bold x}$ is
for the case $k(*)=0$.  So, omitting constant indexes and replacing
sequences of length one by the unique entry we get that it is generated by
$y_{\eta,n}$ (for $\eta \in{}^\omega S,n < \omega$) and $x_\nu$ (for $\nu \in
{}^{\omega >} S$) freely as an abelian group except the equations
$(n!)y_{\eta,n+1} = y_{\eta,n} + x_{\eta \restriction n}$. \nl
Note that if $K$ is the countable subgroup generated by $\{x_\nu:\nu
\in {}^{\omega >} 2\}$ then $G/K$ is a divisible group of cardinality
continuum hence $G$ is not free.  So $G$ is $\aleph_1$-free but not free.
\bn
Now we have the abelian group version of freeness, see generally \scite{af.59}.
\proclaim{\stag{af.56} The Freeness Claim}  Let $\bold x \in \bold K_{k(*)}$.
\nl
1) The abelian group $G_{U \cup
u}/G_{U,u}$ is free \ub{if} $U \subseteq {}^\omega S,u \subseteq
{}^\omega S \backslash U$ and $|u| \le k \le k(*)$ and $|U| \le
\aleph_{k(*)-k}$. \nl
2) If $U \subseteq {}^\omega S$ and $|U| \le \aleph_{k(*)}$, 
\ub{then} $G_U$ is free.
\endproclaim
\bigskip

\proclaim{\stag{af.59} Claim}  1) If $\bold x$ is a combinatorial
$k(*)$-parameter \ub{then} $\bold x$ is $\aleph_{k(*)+1}$-free.
\nl
2) If $\bold x$ is an abelian group parameter and $(k^{\bold x},
S^{\bold x},\Lambda^{\bold x})$ is free, \ub{then} $G_{\bold x}$ is free.
\endproclaim
\bigskip

\demo{Proof}  1) Easily follows by (2).
\nl
2) Similar and follows from \scite{k.14} as easily $G$ belongs to
$({\Cal G}_{(k(*)},S^{\bold x},\Lambda^{\bold x})$.
\enddemo
\bigskip

\proclaim{\stag{af.63} Claim}  Assume $\bold x \in K^{\text{cb}}_{k(*)}$ is full
(i.e. $\Lambda^{\bold x} = {}^{k(*)+1}({}^\omega(S^{\bold x})))$.
\nl
1) If $U \subseteq {}^\omega S$ and $|U| \ge (|S| +
\aleph_0)^{+k(*)+1}$ \ub{Then} $G^{\bold x}_U$ is not free.
\nl
2) If $|S^{\bold x}| \ge \aleph_{k(*)+1}$ \ub{then} $G_{\bold x}$ is
   not free.
\nl
3) Assume $\bold x \in K^{\text{cb}}_{k(*)},|S^{\bold x}_\ell| +
\lambda_\ell < \lambda_{\ell +1}$ for $\ell < k(*)$ and
$|\Lambda^{\bold x}| \ge \lambda_{k_(*)}$ and $G \in 
{\Cal G}_{\bold x}$ (see \S3) \ub{then} $G$ is not free.
\endproclaim
\bigskip

\demo{Proof}  1) Assume toward contradiction that $G_U$ is free and let
$\chi$ be large enough; for notational simplicity assume $|U| =
\aleph_{\alpha +1,k(*)+1}$, this is O.K. as a subgroup of a free abelian group is a
free abelian group.  Let $\aleph_\alpha = |S|$.  We
choose $N_\ell$ by downward induction on $\ell \le k(*)$ such that
\mr
\item "{$(a)$}"  $N_\ell$ is an elementary submodel 
\footnote{${\Cal H}(\chi)$ is $\{x$: the transitive closure of $x$ has
cardinality $< \chi\}$ and $<^*_\chi$ is a well ordering of ${\Cal H}(\chi)$}
of $({\Cal H}(\chi),\in,<^*_\chi)$
\sn
\item "{$(b)$}"  $\|N_\ell\| = |N_\ell \cap \aleph_{\alpha +k(*)}| =
\aleph_\ell$ and $\{\zeta:\zeta \le \aleph_{\alpha + \ell +1}\} 
\subseteq N_\ell$
\sn
\item "{$(c)$}"   $\langle x_{\bar \eta}:\bar \eta \in \Lambda^{\bold
x}_{\le k(*)}\rangle,\langle y_{\bar \eta,n}:\bar \eta \in
\Lambda^{\bold x}$ and $n < \omega\rangle,{\Cal U}$ and $G_{\Cal U}$
belong to $G_U \in N_\ell$ and $N_{\ell +1},\dotsc,N_{k(*)}
\in N_\ell$.
\ermn
Let $G_\ell = G_U \cap N_\ell$, a subgroup of $G_U$.
Now
\mr
\item "{$(*)_0$}"  $G_U/(\Sigma\{G_\ell:\ell \le k(*)\})$ 
is a free (abelian) group 
\nl
[easy or see \cite{Sh:52}, that is: \nl
as $G_U$ is free we can prove by induction on $k \le k(*) + 1$ then
$G/(\Sigma\{G_{k(*)+1-\ell}:\ell < k\})$ is free, for $k = 0$ this is the
assumption toward contradiction, the induction step is by Ax VI in
\cite{Sh:52} for abelian groups and for $k=k(*)+1$ we get the desired 
conclusion.] 
\ermn
Now
\mr
\item "{$(*)_1$}"  letting $U^1_\ell$ be $U$ for $\ell = k(*)+1$ and
$\dbca^{k(*)}_{m=\ell} (N_m \cap U)$ for $\ell \le k(*)$; we have: $U^1_\ell$
has cardinality $\aleph_{\alpha + \ell}$ for $\ell \le k(*) + 1$ \nl
[Why?  By downward induction on $\ell$.
For $\ell = k(*)+1$ this holds by an assumption.  For 
$\ell = k(*)$ this holds by clause (b).  For $\ell < k(*)$ this
holds by the choice of $N_\ell$ as the
set $\dbca^{k(*)}_{m=\ell +1} (N_m \cap U)$ has cardinality
$\aleph_{\alpha + \ell +1} \ge \aleph_\ell$ and belong to $N_\ell$ and clause
(b) above.]
\sn
\item "{$(*)_2$}"  $U^2_\ell =: U^1_{\ell +1} \backslash (N_\ell \cap U)$
has cardinality $\aleph_{\ell +1}$ for $\ell \le k(*)$ \nl
[Why?  As $|U^1_{\ell +1}| = \aleph_{\ell +1} > \aleph_\ell =
\|N_\ell\| \ge |N_\ell \cap U|$.]
\sn
\item "{$(*)_3$}"  for $m < \ell \le k(*)$ the set 
$U^3_{m,\ell} =: U^2_\ell \cap
\dbca^{\ell-1}_{r = m} N_r$ has cardinality $\aleph_{\alpha +m}$   \nl
[Why?  By downward induction on $m$.  For $m=\ell -1$ as $U^2_\ell \in
N_m$ and $|U^2_\ell| = \aleph_{\alpha +\ell +1}$ and clause (b).  
For $m < \ell$ similarly.]
\ermn
Now for $\ell=0$ choose $\eta^*_\ell \in U^2_\ell$, possible by
$(*)_2$ above.  Then for $\ell >0,\ell \le k(*)$ choose 
$\eta^*_\ell \in U^3_{0,\ell}$.  This is possible by $(*)_3$.  So
clearly
\mr
\item "{$(*)_4$}"  $\eta^*_\ell \in U$ and $\eta^*_\ell \in N_m \cap U
\Leftrightarrow \ell \ne m$ for $\ell,m \le k(*)$. \nl
[Why?  If $\ell=0$, then by its choice, $\eta^*_\ell \in U^2_\ell$,
hence by the definition of $U^2_\ell$ in $(*)_2$ we
have $\eta^*_\ell \notin N_\ell$, and 
$\eta^*_\ell \in U^1_{\ell +1}$ hence $\eta^*_\ell \in N_{\ell +1}
\cap \ldots \cap N_{k(*)}$ by $(*)_1$ so $(*)_4$ holds for $\ell = 0$.
If $\ell > 0$ then by its choice, $\eta^*_\ell \in U^3_{0,\ell}$ but
$U^3_{m,\ell} \subseteq U^2_\ell$ by $(*)_3$ so 
$\eta^*_\ell \in U^2_\ell$ hence as
before $\eta^*_\ell \in N_{\ell +1} \cap \ldots \cap N_{k(*)}$ and
$\eta^*_\ell \notin N_\ell$.  Also by $(*)_3$ we have
$\eta^*_\ell \in \dbca^{\ell -1}_{r=0}
N_\ell$ so $(*)_4$ really holds.]
\ermn
Let $\bar \eta^* = \langle \eta^*_\ell:\ell \le k(*) \rangle$ and let $G'$
be the subgroup of $G_U$ generated by $\{x_{\bar \eta 
\upharpoonleft <m,n>}:m \le
k(*)$ and $\bar \eta \in {}^{k(*)+1}U$ and $n < \omega\} 
\cup \{y_{\bar \eta,n}:\bar \eta \in {}^{k(*)+1} U 
\text{ but } \bar \eta \ne \bar \eta^*$ and $n < \omega\}$.
Easily $G_\ell \subseteq G'$ recalling $G_\ell = N_\ell \cap G_U$
hence $\Sigma\{G_\ell:\ell \le k(*)\}
\subseteq G'$, but $y_{\bar \eta^*,0} \notin G'$ hence
\mr
\item "{$(*)_5$}"  $y_{\bar \eta^*,0} \notin \sum\{G_\ell:\ell \le
k(*)\}$. 
\ermn
But for every $n$
\mr
\item "{$(*)_6$}"  $\bar n!y_{\bar \eta^*,n+1} - y_{\bar \eta^*,n} = 
\Sigma\{x_{\bar \eta^* \upharpoonleft <m,n>}:m \le k(*)\} \in 
\Sigma\{G_\ell:\ell \le k(*)\}$. \nl
[Why?  $x_{\bar \eta^* \upharpoonleft <m,n>} \in G_m$ as
$\bar \eta^* \restriction (k(*))+1 \backslash \{m\}) \in N_m$ by
$(*)_4$.]
\ermn
We can conclude that in $G_U/\sum\{G_\ell:\ell \le k(*)\}$, the
element $y_{\bar \eta^*,0} + \sum\{G_\ell:\ell \le k(*)\}$ is not zero
(by $(*)_5$) but is divisible by every natural number by $(*)_6$. \nl
This contradicts $(*)_0$ so we are done.
\nl
2),3)  Left to the reader.  \hfill$\square_{\scite{af.63}}$
\enddemo
\newpage

\head {\S2 Black Boxes} \endhead  \resetall \sectno=2
 \spuriousreset
\bigskip

\proclaim{\stag{ya.7} Claim}  1) Assume $k(*) < \omega,\chi =
\chi^{\aleph_0}$ and $\lambda = \beth_{k(*)}(\chi),S = \lambda,
\Lambda_{k(*)} = {}^{k(*)+1}({}^\omega S)$ or just $S_\ell = \chi_\ell
= \beth_\ell(\chi)$ for $\ell \le k(*)$ and $\Lambda_{k(*)} = \dsize
\prod_{\ell \le k(*)} {}^\omega(S_\ell)$  and $\bold x^{k(*)} =
(k(*),\lambda,\Lambda_{k(*)})$ so $\bold x$ is a full $\langle
S_\ell:\ell \le k(*)\rangle$-parameter.  
\ub{Then} $\Lambda$ has a $\chi$-black box,
i.e. {\rm Qr}$(\Lambda_{\bold x^{k(*)}},\chi)$.
\nl
2) Moreover, $\bold x$ has the $\langle \chi_\ell:\ell \le
k(*)\rangle$-black box, i.e. for every $\bar B = \langle B_{\bar
\eta}:\bar \eta \in \Lambda^{\bold x}_{\le k(*)}\rangle$ satisfying
clause (c) we can find $\langle h_{\bar \eta}:\bar \eta \in
   \Lambda\rangle$ such that:
\mr
\item "{$(a)$}"  $h_{\bar \eta}$ is a function with domain $\{\bar
\eta \upharpoonleft \langle m,n\rangle:m \le k(*),n < \omega\}$
\sn
\item "{$(b)$}"  $h_{\bar \eta}(\bar \eta \upharpoonleft \langle
m,n\rangle) \in B_{\bar \eta \upharpoonleft <m,n-1>}$
\sn 
\item "{$(c)$}"  $B_{\bar \eta \upharpoonleft \langle m,k(*)\rangle}$ 
is a set of cardinality $\beth_m(\chi)$
\sn
\item "{$(d)$}"  if $h$ is a function with domain $\Lambda^{\bold
x}_{\le k(*)}$ such that $h(\bar \eta \upharpoonleft \langle
m,n\rangle) \in B_{(\bar \eta \upharpoonleft <m,n>)}$ and $\alpha_\ell
< \beth_\ell(\chi)$ for $\ell \le k(*)$ then for some $\bar \eta \in
\Lambda^{\bold x},h_{\bar \eta} \subseteq h$ and $\eta_\ell(0) =
\alpha_\ell$ for $\ell \le k(*)$.
\ermn
3) Assume $\chi_\ell = \lambda^{\aleph_0}_\ell,\chi_{\ell +1} =
 \chi^{\chi_\ell}_{\ell +1}$ for $\ell \le k(*)$.  
If $|S_\ell| = \lambda_\ell$ for $\ell \le
k(*),\bold x$ is a full combinatorial $(\bar S,k(*))$-parameter,
and $|B_{\bar \eta \upharpoonleft <m,n>}| \le \chi_m$ for $\bar \eta
 \in \Lambda^{\bold x}$ \ub{then} we can find $\langle h_{\bar
 \eta}:\bar \eta \in \Lambda^{\bold x}\rangle$ as in part (2)
 replacing $\beth_\ell(\chi)$ by $\lambda_\ell$, moreover such that:
\mr
\item "{$(e)$}"  if $\bar \eta \in \Lambda$ then $\eta_\ell$ is
 increasing
\sn
\item "{$(f)$}"  if $\lambda_\ell$ is regular then we can in clause
(d) above add: if for $E_\ell$ is a club of $\lambda_\ell$ for $\ell
\le k(*)$ then we can demand 
$\eta_\ell \in {}^\omega(E_\ell \cup \{\alpha^*_\ell\})$
\sn
\item "{$(g)$}"  if $\lambda_\ell$ is singular $\lambda_\ell =
\Sigma\{\lambda_{\ell,i}:i < \text{\rm cf}(\lambda_\ell)\}$, 
{\rm cf}$(\lambda_{i,\ell}) = \lambda_{i,\ell}$ increasing with $i$ we can add:
if $u_\ell \in [\text{\rm cf}(\lambda_\ell)]$ is unbounded, $E_{\ell,i}$ a
club of $\lambda_{\ell,i}$ then $\eta_\ell \in {}^\omega(E_{i,\ell}
\cup \{\alpha^*_\ell\})$ for some $i$.
\endroster
\endproclaim
\bigskip

\demo{Proof}  Part (1) follows form part (2) which follows from part
(3), so let us prove part (3).
Note that \wilog \, $B_{\bar \nu} = |B_{\bar \nu}|$ and we use
$\alpha_{\bar \eta,m,n} = h_{\bar \eta}(\bar \eta \upharpoonleft
\langle m,n\rangle$ for $\bar \eta \in\Lambda_{\bold x},m \le k(*)$
and $n < \omega$.  We prove
by part (3) by induction on $k(*)$.  Let $\Lambda_k = \Lambda^{\bold x}$.
\mn
\ub{Case 1}:  $k(*)=0$.

By the silly black box, see \cite[III,\S4]{Sh:300}, or better
\cite[VI,\S2]{Sh:e}, see below for details on such a proof.
\mn
\ub{Case 2}:  $k(*) = k+1$.

Let $\langle \alpha^k_{\bar \eta,m,n}:\bar \eta \in \Lambda_k,n <
\omega,m \le k\rangle$ witness part (2) for $k$, i.e. for $\bold x^k$, so no
need to assume $\bold x^k$ is full.  So $\lambda = \lambda_{k(*)},\chi
= \chi_{k(*)}$ and let $\bold H =  \{h:h$ is a function from 
$\Lambda_k$ to $\chi\}$.  So $|\bold H| \le
(\lambda)^{\lambda^{\aleph_0}_k} = \chi$.  
By the silly black box, see below, we can find $\langle \bar
h_\eta:\eta \in {}^\omega \lambda\rangle$ such that
\mr
\item "{$\circledast_1$}"  $(a) \quad \bar h_\eta = \langle h_{\eta,n}:n
<\omega\rangle$ and $h_{\eta,n} \in \bold H$ for $\eta \in {}^\omega \lambda$
\sn
\item "{${{}}$}"  $(b) \quad$ if $\bar f = \langle f_\nu:
\nu \in {}^{\omega >}\lambda\rangle$ and $f_\nu \in \bold H$ 
for every such $\nu$ and $\alpha < \lambda$ \ub{then} 
\nl

\hskip25pt for some increasing 
$\eta \in {}^\omega \lambda$ we have $\alpha = \eta(0)$ and $n <
 \omega \Rightarrow h_{\eta,n}$
\nl

\hskip25pt $= f_{\eta \restriction n}$.
\ermn
[Why?  First assume $\chi =\lambda$.  Let $\langle g_\alpha:\alpha <
\lambda\rangle$ enumerate ${\bold H}$ such that for each $g \in 
{\bold H}$ the set $\{\alpha < \lambda:g_\alpha=g\}$ is unbounded in
$\lambda$.  Now for $\eta \in {}^\omega \lambda$ and $n< \omega$ let
$h_{\eta,n} =g_{\eta(n+1)}$.  So clause (a) holds and as for clause
(b), let $\bar f =  \langle f_\nu:\nu \in {}^{\omega >}\lambda\rangle$
be given, $f_\nu \in \bold H$.  

We choose $\alpha_n$ by induction on $n < \omega$ such that:
\mr
\item "{$(a)$}"  $\alpha_0 = \alpha$
\sn
\item "{$(b)$}"  $\alpha_n < \lambda$ and $\alpha_n > \alpha_m$ if $n
=m+1$
\sn
\item "{$(c)$}"  if $n=m+1 > 1$ then $\alpha_n$ satisfies
$g_{\alpha_n} = f_{\langle \alpha_\ell:\ell < m\rangle}$.
\ermn
Now $\eta =: \langle \alpha_n:n < \omega\rangle$ is as required.  
If $\chi > \lambda$ but still $\chi \le \lambda^{\aleph_0}$, let
$\langle g_\alpha:\alpha < \chi^{\aleph_0}\rangle$ list $\bold H$, and
let $\bold h_n:\chi \rightarrow \lambda$ for $n < \omega$ be such
that $\alpha < \beta < \chi \Rightarrow (\forall^* n)(\bold
h_n(\alpha) \ne \bold h_n(\beta))$ and let cd$:\lambda \rightarrow
{}^{\omega >} \lambda$ be one to one onto.  Now for $\eta \in
{}^\omega \lambda$ and $n < \omega$ let $h_{\eta,n}$ be $g_\alpha$
where $\alpha$ is the unique ordinal $\beta < \chi$ such that for
every $k < \omega$ large enough $(\text{cd}(\eta(k)))(n) = \bold h_n(\alpha)$. 

Next we shall define $\bar \alpha^{k(*)} = \langle \alpha^{k(*)}_{\bar
\eta,m,n}:\bar \eta \in \Lambda_{k+1},m \le k(*),n < \omega\rangle$ as
required; so let $\bar \eta = \langle \eta_\ell:\ell \le k(*)\rangle
\in \Lambda_{k(*)}$ we define $\bar \alpha^{k(*)}_{\bar \eta} = \langle
\alpha^{k(*)}_{\bar \eta,m,n}:m \le k(*),n < \omega\rangle$ as
follows:
\mr
\item "{$(*)$}"   if $\eta_{k(*)} \in {}^\omega \lambda$ and $\langle\eta_0,
\dotsc,\eta_{k(*)-1}\rangle \in \Lambda_k$ \ub{then} for $m \le k(*)$
and $n < \omega$
{\roster
\itemitem{ $(\alpha)$ }   if $m = k(*)$ then $\alpha^{k(*)}_{\bar
\eta,m,n} = h_{\eta_{k(*)},n}(\langle
\eta_0,\dotsc,\eta_{k(*)-1}\rangle) < \lambda_m$
\sn  
\itemitem{ $(\beta)$ }   if $m < k(*)$, i.e. $m \le k$ then 
$\alpha^{k(*)}_{\bar \eta,m,n} = \alpha^k_{\bar \eta \restriction
k(*),m,n} < \lambda_m$.
\endroster}
\ermn
Clearly $\alpha^{k(*)}_{\bar \eta,m,n} < \lambda_m$ we shall prove that
$\langle \bar \alpha^{k(*)}_{\bar \eta,m,n}:\bar \eta \in 
\Lambda^{k+1},m \le k(*),n < \omega\rangle$ witness Qr$(\bold
x^{k(*)},\chi)$, this suffices.
\nl
Why does this hold?  Let $h$ be a function with domain
$\Lambda^{\bold x^{k(*)}}_{\le k(*)}$ as in part (3) and
$\alpha^*_\ell < \lambda_\ell$ for $\ell \le k(*)$.

For $\nu \in {}^{\omega >} \lambda$ let $f_\nu:\Lambda_k
\rightarrow \lambda = \lambda_{k(*)}$ be defined by: 
$f_\nu(\langle \eta_\ell:\ell \le
k\rangle) =: h(\langle \eta_\ell:\ell \le k\rangle \char 94 \langle
\nu \rangle)$.  So by $\circledast_1$ above for some increasing 
$\eta^*_{k(*)} \in {}^\omega \lambda$ we
have $\eta^*_{k(*)}(0) = \alpha^*_{k(*)}$ and
\mr
\item "{$\odot$}"   $\langle \eta_0,\dotsc,\eta_k\rangle \in \Lambda_k
\wedge n < \omega \Rightarrow f_{\eta^*_{k(*)} \restriction n} =
h(\langle \eta_0,\dotsc,\eta^*_k,\eta^*_{\eta(*)}\rangle)$.  
\ermn
Now we define $h'$ with domain $\Lambda^{\bold x^k}_{\le k}$ by:  if 
$\bar \eta \in \Lambda^{\bold x^k}_{\le k}$ then 
$h'(\bar \eta) = h(\bar \eta \char 94 \langle \eta^*_{k(*)}\rangle)$.

So by the choice of $\bar \alpha^k$ we can find $\langle
\eta^*_0,\dotsc,\eta^*_k\rangle \in \Lambda_k$ with no repetitions
such that $\eta^*_\ell(0) = \alpha^*_\ell$ for $\ell \le k$ and

$$
m \le k \wedge n < \omega \Rightarrow \alpha_{\langle \bar
\eta^*_0,\dotsc,\eta^*_k\rangle,m,\ell} = h'(\langle
\eta^*_0,\dotsc,\eta^*_k\rangle \upharpoonleft (m,n)\rangle).
$$
\mn
Now we can check that $\langle
\eta^*_0,\dotsc,\eta^*_k,\eta^*_{k(*)}\rangle$ is as required. 
\hfill$\square_{\scite{ya.7}}$
\enddemo
\bigskip

\demo{\stag{ya.21} Conclusion}  For every $k < \omega$ there is an
$\aleph_{k+1}$-free abelian group $G$ of cardinality $\beth_{k+1}$ and
pure (non-zero) subgroup $\Bbb Z_z \subseteq G$ such that $\Bbb Z z$
is not a direct summand of $G$.
\enddemo
\bigskip

\demo{Proof}  Let $\chi = 2^{\aleph_0}$ and $\bold x$ be a combinatorial
$k$-parmeter as guaranteed by \scite{ya.7}.  Now by \scite{ya.14}(2) below we
can expand $\bold x$ to an abelian group $k$-parameter, so $G_{\bold
x}$ is as required.
\enddemo
\bigskip

\proclaim{\stag{ya.14} Claim}  1) If $\bold x$ is a combinatorial
$k$-parameter such that {\rm Qr}$(\bold x,2^{\aleph_0})$ \ub{then} for some
$\bold a,(\bold x,\bold a)$ is an abelian group $k$-parameter such
that $h \in \text{\rm Hom}(G_{\bold x},\Bbb Z z) \Rightarrow h(z) = 0$.
\nl
2) For every $k$ there is an $\aleph_n$-free abelian group $G$ of
cardinality $\beth_{k+1}$ and $z \in G$ a pure $z \in G$ as above. 
\endproclaim
\bigskip

\demo{Proof}  1) Let $\bar \alpha$ witness Qr$(\bold x,2^{\aleph_0})$.
For each $\bar \eta \in \Lambda^{\bold x}$ we shall choose a sequence
$\langle \bold a_{\bar \eta,n}:n < \omega\rangle$ of integers such
that for any $b \in \Bbb Z \backslash \{0\}$ for no 
$\bar c \in {}^\omega \Bbb Z$ do we have 
(letting $b_{\bar \eta,m,n}$ is: $\alpha_{\bar \eta,m,n}$
if $< \omega$, be $- (\alpha_{\bar \eta,m,n} - \omega)$ if $\alpha_{\bar
\eta,m,n} \in [\omega,\omega + \omega)$ and be $0$ if $\alpha_{\bar \eta,m,n}
\ge \omega + \omega)$:
\mr
\item "{$\boxtimes_{\bar \eta}$}"  for each $n < \omega$ we have
$$
n!c_{n+1} = c_n + \bold a_{\bar \eta,n} b + \Sigma\{b_{\bar \eta,m,n}:m \le
k(*)\}.
$$
\ermn
This is easy: for each pair $(b,c_0) \in \Bbb Z \times \Bbb Z$ there is
at most one sequence $\langle c_0,c_1,c_2,\ldots\rangle$ of integers such
that $\boxtimes_{\bar \eta}$ holds for them, so $\le \aleph_0$
sequences are excluded, so the choice of $\langle \bold a_{\bar
\eta,n}:n < \omega\rangle$ is possible.

Now toward contradiction assume that $h$ is a homomorphism from
$G_{\bold x}$ to $z \Bbb Z$ such that $h(z) = bz,b \in \Bbb Z
\backslash \{0\}$.  We define $h':\Lambda^{\bold x}_{\le k}
\rightarrow \chi$ by $h'(\bar \eta) = n$ if $n < \omega$ and
$h(x_{\bar \eta}) = nz$ and $h'(\bar \eta) = \omega +n$ if $n <
\omega$ and $h(x_{\bar \eta}) = (-n)z$.

By the choice of $\bar \alpha$, for some $\bar \eta \in \Lambda^{\bold
x}$ we have: $m \le k \wedge n < \omega \Rightarrow h'(\bar \eta
\upharpoonleft \langle m,n\rangle) = \alpha_{\bar \eta,m,n}$.  Hence 
$h(x_{\bar \eta \upharpoonleft (m,n)}) = b_{\bar \eta,m,n} z$ for $m
\le k,n < \omega$.

Let $c_n \in \Bbb Z$ be such that $h(y_{\bar \eta,n}) = c_n z$.  Now
the equation $\boxtimes_{\bar \eta,n}$ in Definition \scite{6.1} is
mapped to the $n$-th equation in $\boxtimes_{\bar \eta}$, so an
obvious contradiction.  \hfill$\square_{\scite{ya.14}}$
\nl
2) By part (1) and \scite{ya.21}.
\enddemo
\bigskip

\remark{\stag{ya.28} Remark}  1) We can replace $\chi$ by a set of
cardinality $\chi$ in Definition \scite{ya.77}.  Using $\Bbb Z z$
instead of $\chi$ simplify the notation in the proof of \scite{ya.14}.
\nl
2) We have not tried to save in the cardinality of $G$ in
\scite{ya.14}(2), using as basic of the induction the abelian group
of cardinality $\aleph_0$ or $\aleph_1$.
\endremark
\bigskip

\proclaim{\stag{ya.84} Claim}   1) If $\chi_0 =
\chi^{\aleph_0}_0,\chi_{m+1} = 2^{\chi_m}$ and $\lambda_m = \chi_m$
for $m \le k$ \ub{there} the $\bar \chi$-full has the $\bar
\chi$-black box.
\endproclaim
\bigskip

\demo{\stag{ya.91} Conclusion}  Assume $\mu_0 < \ldots < \mu_{k(*)}$ are strong
limit of cofinality $\aleph_0$ (or $\mu_0 = \aleph_0$), $\lambda_\ell
= \mu^+_\ell,\chi_\ell =2^{\mu_\ell}$.

Then in \scite{ya.7} for $\bar \eta \in \Lambda^{\bold x}$ we can let
$h_{\bar \eta,m}$ has domain $\{\bar \nu \in \Lambda^{\bold
x}_m:[\nu_\ell = \eta_\ell$ for $\ell = m+1,\dotsc,k(*)\}$.
\enddemo
\newpage

\head {\S3 Constructing abelian groups from combinatorial parameters}
\endhead  \resetall 
 \spuriousreset
\bigskip

\definition{\stag{k.7} Definition}  1) We say $F$ is a $\mu$-regressive
function on a combinatorial parameter $\bold x \in K^{\text{cb}}_{k(*)}$ 
when: $S^{\bold x}$
is a set of ordinals and:
\mr
\item "{$(a)$}"  Dom$(F)$ is $\Lambda^{\bold x}$
\sn
\item "{$(b)$}"  Rang$(F) \subseteq [\Lambda^{\bold x} \cup
\Lambda^{\bold x}_{\le k(*)}]^{\le \aleph_0}$
\sn
\item "{$(c)$}"  for every $\bar \eta \in \Lambda^{\bold x}$ and $\ell
\le k(*)$ we \footnote{actually, suffice to have it for $\ell=k(*)$}
have sup Rang$(\eta_\ell) >
\sup(\cup\{\text{Rang}(\nu_\ell):\bar \nu \in F(\bar \eta)\})$.
\ermn
1A) We say $F$ is finitary when $F(\bar \eta)$ is finite for every
$\bar \eta$.
\nl
1B) We say $F$ is simple if $\eta_{k(*)}(0)$ determined $F(\bar \eta)$
for $\bar \eta \in \Lambda^{\bold x}$.
\nl
2) For $\bold x,F$ as above and $\Lambda \subseteq \Lambda^{\bold x}$
 we say that $\Lambda$ is free for $(\bold x,F)$ when: there is a
 sequence $\langle \bar \eta^\alpha:\alpha < \alpha(*)\rangle$
 listing $\Lambda' = \Lambda \cup \bigcup\{F(\bar \eta):\bar \eta
 \in \Lambda\}$ and sequence $\langle \ell_\alpha:\alpha <
 \alpha(*)\rangle$ such that
\mr
\item "{$(a)$}"  $\ell_\alpha \le k(*)$
\sn   
\item "{$(b)$}"  if $\alpha < \alpha(*)$ and $\bar \eta^\alpha \in
 \Lambda$ then $F(\bar \eta^\alpha) \subseteq \{\bar \eta^\beta,\bar
 \eta^\beta \upharpoonleft \langle m,n\rangle:\beta < \alpha,n <
 \omega,m \le k(*)\}$ 
\sn
\item "{$(c)$}"  if $\alpha < \alpha(*),\bar \eta^\alpha \in \Lambda$
 then for some $n < \omega,\bar \eta^\alpha \upharpoonleft \langle
 \ell_\alpha,n\rangle \notin \{\bar \eta^\beta \upharpoonleft \langle
 \ell_\alpha,n\rangle:\beta < \alpha,\eta^\beta \in \Lambda\} \cup
 \{\bar \eta^\beta:\beta < \alpha\}$.
\ermn
3) We say $\bold x$ is $\theta$-free for $F$ is $(\bold x,F)$ is
 $\mu$-free \ub{when} $\bold x,F$ are as in part (1) and every
 $\Lambda \subseteq \Lambda^{\bold x}$ of cardinality $< \theta$ is free
 for $(\bold x,F)$.
\enddefinition
\bigskip

\proclaim{\stag{k.14} Claim}  1) If $\bold x \in K^{\text{cb}}_{k(*)}$
 and $F$ is a regressive function on $\bold x$ \ub{then} $(\bold x,F)$
 is $\aleph_{k(*)+1}$-free provided that $F$ is finitary or simple.
\nl
2) In addition: if $k \le k(*),\Lambda \subseteq \Lambda^{\bold x}$ 
has cardinality $\le \aleph_k$ and $\bar u = \langle
 u_{\bar \eta}:\bar \eta \in \Lambda\rangle$ satisfies 
$u_{\bar \eta} \subseteq \{0,\dotsc,k(*)\},|u_\eta| >k$, \ub{then} we can
 find $\langle \bar \eta_\alpha:\alpha < \aleph_k\rangle,\langle
\ell_\alpha:\alpha < \aleph_k\rangle,\langle n_\alpha:\alpha <
\aleph_k\rangle$ such that:
\mr
\item "{$(a)$}"   $\Lambda \subseteq \{\bar \eta^\alpha:\alpha <
\aleph_k$
\sn
\item "{$(b)$}"  if $\bar \eta_\alpha \in \Lambda^{\bold x}$ then
$\ell_\alpha \in u_{\bar \eta^\alpha},n_\alpha < \omega$
\sn
\item "{$(c)$}"   $\bar \eta^\alpha \upharpoonleft
\langle \ell_\alpha,n_\alpha\rangle \notin \{\bar \eta^\beta
\upharpoonleft \langle \ell_\alpha,n_\alpha\rangle:\beta < \alpha\}
\cup\{\bar \eta^\beta:\beta < \alpha\}$.
\endroster
\endproclaim
\bigskip

\demo{Proof}  1) Follows by part (2) for the case $k=k(*),u_{\bar
\eta} = \{0,\dotsc,k(*)\}$ for every $\bar \eta \in \Lambda$.
\nl
2) Without loss of generality $\Lambda$ is closed under $\bar \eta
\mapsto F(\bar \eta) \cap \Lambda^{\bold x}$.
We prove this by induction on $k$.  
\enddemo
\bn
\ub{Case 1}:  $k=0$.
\bn
\ub{Subcase 1A}:  Ignoring $F$.

Let $\langle \bar \eta^\alpha:\alpha < |\Lambda|\rangle$ list
$\Lambda$ with no repetitions (so $\alpha < |\Lambda| \Rightarrow
\alpha < \omega$).  Now $\alpha < |\Lambda| 
\Rightarrow u_{\bar \eta^\alpha} \ne \emptyset$ and let $\ell_\alpha = 
\text{ min}(u_{\bar \eta^\alpha}) \le k(*)$.
Now for each $\alpha < |\Lambda|$ we know that $\beta < \alpha
\Rightarrow \bar \eta^\beta \ne \bar\eta^\alpha$, hence for some $n =
n_{\alpha,\beta} < \omega$ we have $\bar \eta^\beta \upharpoonleft
\langle \ell_\alpha,n_{\alpha,\beta}\rangle \ne \bar \eta^\alpha
\upharpoonleft \langle \ell_\alpha,n_{\alpha,\beta}\rangle$.

Let $n_\alpha = \sup\{n_{\alpha,\beta}:\beta < \alpha\rangle$ it is $<
\omega$ as $\alpha < \omega$.  Now $\langle
(\ell_\alpha,n_\alpha):\alpha < |\Lambda|\rangle$ is as required.
\bn
\ub{Subcase 1B}:  $\bar \eta \in \Lambda \Rightarrow F(\bar \eta)$ is
finite.

Let $\langle \eta^1_\alpha:\alpha < |\Lambda|\rangle$ list $\Lambda$,
we choose $w_j$ by induction on $j \le j(*),j(*) < \omega$ such that:
\mr
\item "{$(a)$}"  $w_j \subseteq |\Lambda|$ is finite
\sn
\item "{$(b)$}"  $j \in w_{j+1}$
\sn
\item "{$(c)$}"  if $\alpha \in w_j$ then $F(\bar \eta^\alpha) \cap
\Lambda \subseteq \{\bar \eta^\alpha:\beta \in \Lambda_j\}$
\sn
\item "{$(d)$}"  $w_{j(*)} = (\Lambda)$ and $w_0 = \emptyset$
\sn
\item "{$(e)$}"  $w_j \subseteq w_{j+1}$.
\ermn
No problem to do this.

Now let $\langle \beta(j,i):i < i^*_j\rangle$ list $w_{j+1} \backslash
w_j$ such that: if $i_1,i_2 < i^*_j$ and $\bar \eta^{\beta(j,i_1)} \in
F(\bar \eta^{\beta(j,i_2))}$ then $i_1 < i_2$ existence by $F$ being
regressive.  Let $\langle \bar
\nu_{j,i}:i < i^{**}_j\rangle$ list $\cup\{F(\bar \eta^\alpha):\alpha
\in w_{j+1} \backslash w_j\} \backslash \Lambda^{\bold x} \backslash
\{F(\bar \eta^\alpha):\alpha \in w_j\}$.

Let $\alpha^*_j = \Sigma\{i^{**}_{j(1)} + i^*_{j(2)}:j(1) < j\}$.  Now
we define $\bar\rho_\varepsilon$ for $\varepsilon < \alpha^*_j$ for 
$j < j(*)$ as follows:
\mr
\item "{$(a)$}"  $\rho_{\alpha^*_j +i} = \nu_{j,i}$ if $i < i^{**}_j$
\sn
\item "{$(b)$}"  $\bar\rho_{\alpha^*_j+i^{**}_j+i} = \bar
\eta^{\beta(j,i)}$ if $i<i^*_j$.
\ermn
Lastly, we choose $\eta_{\alpha_j + i^{**}_{j+i}} < \omega$ as in case 1A.

Now check. 
\bn
\ub{Subcase 1C}:  $F$ is simple.

Note that $F(\bar \eta)$ when defined is determined by
$\eta_{k(*)}(0)$ and is included in $\{\bar \nu \in \Lambda^{\bold
x}_{\le k(*)} \cup \Lambda^{\bold x}:\text{ sup Rang}(\nu_{k(*)}) <
\eta_{k(*)}(0)\}$.  So let $u = \{\eta_{k(*)}(0):\bar \eta \in
\Lambda\}$ and $u^* = u \cup\{\sup(u)+1\}$ and 
for $\alpha \in u$ let $\Lambda_\alpha = \{\bar \eta
\in \Lambda:\eta_{k(*)}(0) = \alpha\}$ and let $\Lambda_{< \alpha} =
\cup\{\Lambda_\alpha:\alpha \in u\}$.  Now by induction on $\beta \in
u^+$ we choose $\langle(\bar
\eta^\varepsilon,\ell_\varepsilon):\varepsilon <
\varepsilon_\beta\rangle$ such that it is a required for $\Lambda_{<
\alpha}$.  For $\beta = \text{ min}(u)$ this is trivial and if otp$(u
\cap \beta)$ is a limit ordinal this is obvious.  So assume $\alpha =
\text{ max}(u \cap \beta)$, we use
Subcase 1A on $\Lambda_\alpha$, and combine them naturally promising
$\ell_\alpha = k(*) \Rightarrow n_\alpha > 1$. 
\bn
\ub{Case 2}:  $k=k_* +1$.

Let $\langle \Lambda_\varepsilon:\varepsilon < \aleph_k\rangle$ be
$\subseteq$-increasing continuous with union $\Lambda,|\Lambda_{1+
\varepsilon}|  =
\aleph_{k_*},\Lambda_0 = \emptyset$, each $\Lambda_\varepsilon$ closed
enough, mainly:
\mr
\item "{$\circledast_1$}"  if $\bar \eta^i \in \Lambda_\varepsilon$
for $i<i(*) < \omega,\bar\rho \in \Lambda$ and $\{\rho_\ell:\ell \le
k(*)\} \subseteq \{\eta^i_\ell:\ell \le k(*),i<i(*)\}$ \ub{then} $\bar
\rho \in \Lambda$
\sn
\item "{$\circledast_2$}"  $\Lambda_\varepsilon$ is closed under $\bar
\eta \mapsto F(\bar \eta) \cap \Lambda^{\bold x}$.
\ermn
Easily
\mr
\item "{$\odot$}"  if $\varepsilon < \aleph_k,\bar \eta \in
\Lambda_{\varepsilon +1} \backslash \Lambda_\varepsilon$ then
$u'_{\bar \eta} = \{\ell \in u_{\bar \eta}:\eta_\ell$ belongs to
$\{\nu_\ell:\bar \nu \in \Lambda_\varepsilon\}\}$ has at most one
member.
\ermn
 Apply the induction hypothesis to $\Lambda_{\varepsilon +1} \backslash
\Lambda_\varepsilon$ for each $\varepsilon$ and combine but for
$\Lambda_{\varepsilon +1} \backslash \Lambda_\varepsilon$ we use
$\langle u_{\bar \eta} \backslash u'_{\bar \eta}:\bar \eta \in
\Lambda_{\varepsilon +1} \backslash \Lambda_\varepsilon\rangle$, so
$|u_{\bar \eta} \backslash u'_{\bar \eta}| \ge k-1 = k_*$.
\hfill$\square_{\scite{k.14}}$ 
\bigskip

\definition{\stag{k.21} Definition}  For a combinatorial parameter
$\bold x$ we define ${\Cal G}_{\bold x}$, the class of abelian groups
derived from $\bold x$ as follows: $G \in {\Cal G}_{\bold x}$ is there
is a simple (or finitary) regressive $F$ on $\Lambda^{\bold x}$ 
and $G$ is generated by $\{y_{\bar \eta,n}:\eta \in \Lambda^{\bold x},
n < \omega\} \cup \{x_{\bar
\eta}:\bar \eta \in \Lambda^{\bold x}_{\le k(*)}\}$ freely except
\mr
\item "{$\boxtimes_{\bar \eta,n}$}"  $(n!)y_{\bar \eta,n+1} = y_{\bar \eta,n}
+ b^{\bold x}_{\bar \eta,n} z_{\bar \eta,n} + 
\sum\{x_{\bar \eta \upharpoonleft <m,n>}:m \le k(*)\}$
\ermn
where
\mr
\item "{$\odot$}"   $(a) \quad b_{\bar \eta,n} \in \Bbb Z$
\sn
\item "{${{}}$}"  $(b) \quad z_{\bar \eta,n}$ is a linear combination
of

$$
\align
\{x_{\bar \nu}:\bar \nu \in F(\bar \eta) \backslash \Lambda^{\bold
x}\} \cup \{y_{\bar \eta,n}:&\bar \eta \in F(\bar \eta) \cap
\Lambda^{\bold x} \text{ and} \\
  &(\forall m \le k(*))(\bar \eta \upharpoonleft \langle m,n\rangle)
  \in F(\bar \eta)\}.
\endalign
$$
\endroster
\enddefinition
\bigskip

\proclaim{\stag{k.28} Claim}  If $\bold x \in K^{\text{cb}}_{k(*)}$ and
  $G \in {\Cal G}_{\bold x}$ (i.e. $G$ is an abelian group derived
  from $\bold x$, \ub{then} $G$ is $\aleph_{k(*)+1}$-free.
\endproclaim
\bigskip

\demo{Proof}  We use claim \scite{k.14}.  So let $H$ be a subgroup of
$G$ of cardinality $\le \aleph_{k(*)}$.  We can find $\Lambda$ such
that
\mr
\item "{$(*)$}"  $(a) \quad \Lambda \subseteq \Lambda^{\bold x}$ has
cardinality $\le \aleph_{k(*)}$
\sn
\item "{${{}}$}"  $(b) \quad$ every equation which $X_\Lambda =
\{x_{\bar \eta \upharpoonleft<m,n>},y_{\bar \eta,n}:m \le k(*),n <
\omega,\bar \eta \in \Lambda\}$ satisfies in $G$, is implied by the
equations in $\Gamma_\Lambda = 
\{\boxtimes_{\bar \eta,n}:\bar \eta \in \Lambda\}$
\sn 
\item "{${{}}$}"  $(c) \quad H \subseteq G_\Lambda = \langle 
x_{\bar \eta \upharpoonleft <m,n>},y_{\bar \eta,n}:\bar \eta \in
\Lambda,m \le k(*),n < \omega\rangle$.
\ermn
So it sufices to prove that $G_\Lambda$ is a free (abelian) group.

Let $\langle(\bar \eta^\alpha,\ell_\alpha):\alpha < \alpha(*)\rangle$
be as proved to exist in \scite{k.14}.  Let ${\Cal U} = \{\alpha <
\alpha(*):\bar \eta^\alpha \in \Lambda\} \cup \{\alpha(*)\}$ 
and for $\alpha \in {\Cal U}$ let $X^0_\alpha = 
\{x_{\bar \eta^\beta \upharpoonleft <m,n>}:\beta \in \alpha \cap {\Cal U},
m \le k(*)$ and $n < \omega\}$ and $X^1_\alpha = X^0_\alpha \cup
\{\bar \eta_\beta:\beta \in \alpha \backslash {\Cal U}\}$.
 So for $\alpha \in {\Cal U}$ there is 
$\bar n_\alpha = \langle n_{\alpha,\ell}:\ell \in
v_\alpha\rangle$ such that: $\ell_\alpha \in v_\alpha \subseteq
\{0,\dotsc,v_\alpha\},n_{\alpha,\ell} < \omega$ and $X^1_{\alpha +1}
\backslash X^1_\alpha = \{x_{\bar \eta
\upharpoonleft <\ell,n>}:\ell \in v_\alpha$ and $n \in
[n_{\alpha,\ell},\omega)\}$.

For $\alpha \le \alpha(*)$ let $G_{\Lambda,\alpha} = 
\langle\{y_{\bar \eta^\beta,n}:x_{\bar \nu}:\beta \in {\Cal U} 
\cap \alpha,\bar \nu \in X^1_\beta\}\rangle_{G_\Lambda}$.  
Clearly $\langle G_{\Lambda,\alpha}:\alpha \le
\alpha(*)\rangle$ is purely increasing continuous with union
$G_\Lambda$, and $G_{\Lambda,0} = \{0\}$.  So it suffices to prove
that $G_{\Lambda,\alpha +1}/G_{\Lambda,\alpha}$ is free.  If $\alpha
\notin {\Cal U}$ the quotient is the trivial group, and if $\alpha \in
{\Cal U}$ we can use $\ell_\alpha \in v_\alpha$ to prove that is free
giving a basis.  \hfill$\square_{\scite{k.28}}$
\enddemo
\bigskip

\demo{\stag{k.35} Conclusion}  For every $k(*) < \omega$ there is an
$\aleph_{k(*)+1}$-free abelian group $G$ of cardinality $\lambda =
\beth_{k(*)+1}$ such that Hom$(G,\Bbb Z) = \{0\}$.
\enddemo
\bigskip

\demo{Proof}  We use $\bold x$ and $\langle h_{\bar \eta}:\bar \eta
\in \Lambda^{\bold x}\rangle$ from \scite{ya.7} (?), and we shall choose $G \in
{\Cal G}_{\bold x}$.  So $G$ is $\aleph_{k(*)+1}$-free by \scite{k.28}.  

Let ${\Cal S} = \{\langle(a_i,\bar \eta_i):i < i_1\rangle \char 94
\langle(b_j,\bar \nu_j,n_j):j < j_1\rangle:i_1 < \omega,a_i \in \Bbb
Z,\bar \eta_i \in \Lambda^{\bold x}_{\le k(*)}$ and $j_1 < \omega,b_j
\in \Bbb Z,\nu_j \in \Lambda^{\bold x},n_j < \omega\}$ (actually $S =
\Lambda^{\bold x}_{\le k(*)}$ suffice).

So $|S| = \lambda_{k(*)}$ and let $\bar p$ be such that:
\mr
\item "{$(a)$}"  $\bar p = \langle p^\alpha:\alpha < \lambda\rangle$
\sn
\item "{$(b)$}"  $\bar p$ lists ${\Cal S}$
\sn
\item "{$(c)$}"  $p^\alpha = \langle (a^\alpha_i,\bar \eta^\alpha_i):i
< i_\alpha\rangle \char 94 \langle
(b^\alpha_j,\nu^\alpha_j,\eta^\alpha_j):j < j_\alpha\rangle$
\sn
\item "{$(d)$}"  sup Rang$(\eta^\alpha_{i,k(*)}) < \alpha$, sup
Rang$(\nu6\alpha_{j',k(*)}) < \alpha$ if $i < i_\alpha,j < j_\alpha$.
\ermn
Now to apply definition \scite{k.21} we have to choose $z_\alpha$ (for
Definition \scite{k.21} as $\Sigma\{a^\alpha_i,x_{\eta_i}:i <
i_\alpha\} + \Sigma\{b^\alpha_j,y_{\bar \nu^\alpha_j},\eta^\alpha_j:j <
j_\alpha\}$ and $z_{\bar \eta} = z_{\eta_{k(*)}}(0)$ for $\bar \eta
\in \Lambda^{\bold x}$ then for $\bar \eta \in \Lambda^{\bold x}$ 
we choose $\langle b_{\bar \eta,n}:n < \omega\rangle
\in {}^\omega \Bbb Z$ such that: there is no function $h$ 
from $\{z_{\bar \eta}\} \cup
\{y_{\bar \eta,n}:n < \omega\} \cup \{x_{\bar \nu \restriction
<m,n>}:m \le k(*),n < \omega\}$ into $\Bbb Z$ satisfying
\mr
\item "{$\circledast$}"  $(a) \quad h(z_{\bar \eta}) \ne 0$ and
\sn
\item "{${{}}$}"  $(b) \quad h(x_{\bar \eta \upharpoonleft <m,n>}) =
h(\bar \eta \upharpoonleft \langle m,n \rangle)$ for $m \le k(*),n < \omega$
\sn
\item "{${{}}$}"  $(c) \quad$ for every $n$
{\roster
\itemitem{ $(*)_n$ }  $n!h(y_{\bar \eta,n+1}) = h(y_{\bar \eta,n}) +
b_{\bar \eta,n} h(z_{\bar \eta}) + 
\Sigma\{\{x_{\bar \eta \upharpoonleft <m,n>}):m \le k(*)\}$.
\endroster}
\ermn
E.g. for each $\rho \in {}^\omega 2$ we can try $b^\rho_n = \rho(n)$
and assume toward contradiction that for each $\rho \in {}^\omega 2$
there is $h_\rho$ as above.  Hence for some $c \in \Bbb Z \backslash
\{0\}$ the set $\{\rho \in {}^\omega 2:h_\rho(z_{\bar \eta})=c\}$ is
uncountable.  So we can find $\rho_1 \ne \rho_2$ such that $h_{\rho_1}
= c = h_{\rho_2}(x_\nu)$ and $\rho_1 \restriction (|c| +7) = \rho_2
\restriction (|c|+7)$.  So for some $n \ge |c| +7,\rho_1 \restriction
n = \rho_2 \restriction n$ and $\rho_1(n) \ne \rho_2(n)$.

Now consider the equation $(*)_n$ for $h_{\bar\rho_1}$ and $h_{\bar
\rho_2}$, subtract them and get $(\rho_1(n)-\rho_2(n)) c$ is divisible
by $n!$, clear contradiction.  So $G \in {\Cal G}_{\bold x}$ is well
defined and is $\aleph_{k(*)+1}$-free by \scite{k.28}.  Suppose $h \in
\text{ Hom}(G,\Bbb Z)$ is non-zero, so for some $\alpha <
\lambda_{k(*)},h(z_\alpha) \ne 0$ (actually as $G^1 = \langle \{x_{\bar
\nu}:\bar \nu \in \Lambda^{\bold x}_{\le k(*)}\}\rangle_G$ is a
subgroup such that $G/G^1$ is divisible necessarily $h \restriction
G^1$ is not zero hence in \scite{ya.7}(2) for 
some $\bar \nu \in \Lambda^{\bold x}_{\le
k(*)}$ we have $h(x_{\bar \nu}) \ne 0$.  Let $\bold y = \{\bar \nu\}$
and so by the choice of $\langle h_{\bar \eta}:\bar \eta \in
\Lambda\rangle$ for some $\bar \eta \in \Lambda^{\bold x}$ we have
$h_{\bar \eta} = h \restriction
\{x_{\bar \eta \upharpoonleft <m,n>}:m \le k(*),n < \omega\}$.  We
clearly get a contradiction.  \hfill$\square_{\scite{k.35}}$
\enddemo
\bigskip

\remark{Remark}  We can give more details as in the proof of \scite{ya.14}.
\endremark
\bigskip

\demo{\stag{ya.42} Conclusion}  For every $n \le m < \omega$ there is a
purely increasing sequence $\langle G_\alpha:\alpha \le
\omega_n+1\rangle$ of abelian groups, $G_\alpha,G_\beta/G_\alpha$ are
free for $\alpha < \beta \le \omega_n$ and
$G_{\omega_n+1}/G_{\omega_n}$ is $\aleph_n$-free and for some $h \in
\text{ Hom}(G_\kappa,\Bbb Z)$ has no extension in
Hom$(G_{\omega_n+1},\Bbb Z)$. 
\enddemo
\bigskip

\demo{Proof}  Let $G,z$ be as in \scite{ya.21}.  So also $G/\Bbb Z z$
is $\aleph_n$-free.  Let $G_\alpha = \langle \{z\}\rangle_G$ for
$\alpha \le \omega_2,G_{\omega_n+1} = G$.
\enddemo
\newpage

\head {\S4 Appendix 1} \endhead  \resetall \sectno=4
 \spuriousreset
\bigskip

\demo{\stag{af.19} Notation}  If $\bar \eta^* \in \Lambda^{\bold x}_m$ and $\bar
\eta = \bar \eta^* \restriction \{\ell \le k(*):\ell \ne m\}$ and $\nu
= \eta^*_m$ then let $x_{m,\bar \eta,\nu} := x_{\bar \eta^*}$.  (See
proof of \scite{af.56}).
\enddemo
\bigskip

\demo{Proof of \scite{af.28}}  Let $U \subseteq {}^\omega S$ be countable (and
infinite) and define $G'_U$ like $G$
restricting ourselves to $\eta_\ell \in U$; by the L\"owenheim-Skolem
argument it suffices to prove that $G'_U$ is a free abelian group.  
List $\Lambda \cap {}^{k(*)+1}U$ without
repetitions as $\langle \bar \eta_t:t < t^* \le \omega \rangle$, and 
choose $s_t < \omega$ by induction on $t < \omega$  such that 
$[r < t \and \bar \eta_r \restriction k(*) = \bar \eta_t
\restriction k(*) \Rightarrow \emptyset = \{\eta_{t,k(*)} 
\restriction \ell:\ell \in [s_t,\omega)\} \cap 
\{\eta_{r,k(*)} \restriction \ell:\ell \in [s_r,\omega)\}]$.
\sn
Let
\medskip

$\qquad \quad Y_1 = \{x_{m,\bar \eta,\nu}:m < k(*),\bar \eta \in
{}^{k(*)+1 \backslash \{m\}}U
\text{ and } \nu \in {}^{\omega >}2\}$

$$
\align
Y_2 = \biggl\{ x_{m,\bar \eta,\nu}:&m = k(*),\bar \eta \in {}^{k(*)}U
\text{ and for no } t < t^* \text{ do we have} \\
  &\bar \eta = \bar \eta_t \restriction k(*) \and 
\nu \in \{\eta_{t,k(*)} \restriction \ell:s_t \le \ell < \omega\} \biggr\}
\endalign
$$

$\qquad \quad Y_3 = \{y_{\bar \eta_t,n}:t < t^* \text{ and } n \in
[s_t,\omega)\}$.
\mn
Now
\mr
\item "{$(*)_1$}"  $Y_1 \cup Y_2 \cup Y_3 \cup \{z\}$ generates
$G'_U$.
\ermn
[Why?  Let $G'$ be the subgroup of $G'_U$ which $Y_1 \cup Y_2 \cup Y_3$
generates.  First we prove by induction on $n < \omega$ that for $\bar \eta
\in {}^{k(*)}U$ and $\nu \in {}^n S$ we have $x_{k(*),\bar
\eta,\nu} \in G'$.  If $x_{k(*),\bar \eta,\nu} \in Y_2$ this is clear;
otherwise, by the definition of $Y_2$ for some $\ell < \omega$ (in
fact $\ell=n$) and $t
< \omega$ such that $\ell \ge s_t$ we have $\bar \eta = \bar \eta_t
\restriction k(*),\nu = \eta_{t,k(*)} \restriction \ell$.

Now
\mr
\item "{$(a)$}"   $y_{\bar \eta_{t,\ell +1}},
y_{{\bar \eta}_{t,\ell}}$  are in $Y_3 \subseteq G'$ 
\sn
\item "{$(b)$}"  $x_{m,\bar \eta_t \restriction \{i \le k(*):i \ne
m\},\nu}$ belong to $Y_1 \subseteq G'$ if $m < k(*)$.
\ermn
Hence by the equation $\boxtimes_{\bar \eta,n}$ in Definition
\scite{6.1}, clearly $x_{k(*),\bar \eta,\nu} \in G'$.  So as $Y_1
\subseteq G' \subseteq G'_U$, all the generators of the 
form $x_{m,\bar \eta,\nu}$
with each $\eta_\ell \in U$ are in $G'$.  

Now for each $t < \omega$ we prove that all the generators $y_{{\bar
\eta}_t,n}$ are in $G'$.
If $n \ge s_t$ then clearly $y_{\bar \eta_t,n} \in Y_3 \subseteq G'$.
So it suffices to prove this for
$n \le s_t$ by downward induction on $n$; for $n = s_t$ by an earlier
sentence, for $n < s_t$ by $\boxtimes_{\bar \eta,n}$.
The other generators are in this subgroup so we
are done.]
\mr
\item "{$(*)_2$}"  $Y_1 \cup Y_2 \cup Y_3 \cup \{z\}$ generates $G'_U$
freely. 
\nl
[Why?  Translate the equations, see more in \cite[\S5]{Sh:771}.] 
\nl
${{}}$  \hfill$\square_{\scite{af.28}}$
\endroster
\enddemo
\bigskip

\demo{Proof of \scite{af.42}}  0), 1) Obvious. 
\nl
2),3),4)   Follows.
\nl
5) Let $\langle \eta_\ell:\ell < m(*)\rangle$ list $u,U_\ell =
U \cup (u \backslash \{\eta_\ell\})$ so $G_{U,u} = G_{U^+_0} \ldots +
G_{U_{m(*)-1}}$.
First, $G_{U,u} \subseteq G_{U \cup u}$ follows by the definitions.
Second, we deal with proving $G_{U,u} \subseteq_{\text{pr}} G_{U \cup u}$.  
So assume $z^*
\in G,a^* \in \Bbb Z$ and $a^* z^*$ belongs to $G_{U_0} + \ldots +
G_{U_{m(*)}}$ so it has the form $\Sigma\{b_i x_{\bar \eta'
\upharpoonleft <m_i,n_i>}:i < i(*)\} + \Sigma\{c_j y_{\bar \eta_j,n_j}:j <
j(*)\} + az$ with $i(*) < \omega,j(*) < \omega$ and $a^*,b_i,c_j 
\in \Bbb Z$ and $\nu_i,\bar \eta^i,\bar \eta_j$ are suitable
sequences of members of $U_{\ell(i)},U_{\ell(i)},U_{k(j)}$ respectively where
$\ell(i),k(j) < m(*)$. We continue as in \cite{Sh:771}.
\nl
6) Easy.
\nl
7)  Clearly $U_1 \cup v = U_2 \cup u$ hence 
$G_{U_1 \cup u} \subseteq G_{U_1 \cup v} = G_{U_2 \cup u}$ hence
$G_{U,u} + G_{U_1 \cup u}$ is a subgroup of $G_{U,u} + G_{U_2 \cup u}$,
so the first quotient makes sense.

Hence $(G_{U,u} + G_{U_2 \cup u})/(G_{U,u} + G_{U_1 \cup u})$ is isomorphic
to $G_{U_2 \cup u}/(G_{U_2 \cup u} \cap (G_{U,u} + G_{U_1 \cup u}))$.
Now $G_{U_1,v} \subseteq G_{U_1 \cup v} = G_{U_2 \cup v} \subseteq
G_{U,u} + G_{U_2,u}$ and $G_{U_1,v} \subseteq G_{U,v} = G_{U,v
\backslash U} = G_{U,u} \subseteq G_{U,u} + G_{U_2,u}$.  Together
$G_{U_1,v}$ is included in their intersection, i.e. 
$G_{U_2 \cup u} \cap (G_{U,u} + G_{U_1 \cup u})$ include $G_{U_1,v}$ and
using part (1) both has the same divisible hull inside $G^+$.  But
as $G_{U_1,v}$ is a pure subgroup of $G$ by part (5) hence of $G_{U_1
\cup v}$.  So necessarily $G_{U_1 \cup u} \cap 
(G_{U,u} + G_{U_1,u}) = G_{U_1,v}$, so
as $G_{U_2 \cup u} = G_{U_1 \cup v}$ we are done.
\nl
8) See \cite{Sh:771}.      \hfill$\square_{\scite{af.42}}$
\enddemo
\bigskip

\demo{Proof of \scite{af.56}}  1) We prove 
this by induction on $|U|$; \wilog \, $|u|=k$
as also $k' = |u|$ satisfies the requirements.
\enddemo
\bn
\ub{Case 1}:  $U$ is countable.

So let $\{\nu^*_\ell:\ell < k\}$ list $u$ be with no repetitions, now if
$k=0$, i.e. $u = \emptyset$ then $G_{U \cup u} = G_U = G_{U,u}$ so the
conclusion is trivial.  Hence we assume $u \ne \emptyset$, and let $u_\ell :=
u \backslash \{\nu^*_\ell\}$ for $\ell < k$.

Let $\langle \bar \eta_t:t < t^* \le \omega \rangle$ list with no
repetitions the set $\Lambda_{U,u} := \{\bar \eta \in \Lambda^{\bold x} 
\cap {}^{k(*)+1}(U \cup u)$: for no
$\ell < k$ does $\bar \eta \in {}^{k(*)+1}(U \cup u_\ell)\}$.  Now
comes a crucial point: let $t < t^*$, 
for each $\ell < k$ for some $r_{t,\ell} \le
k(*)$ we have $\eta_{t,r_{t,\ell}} = \nu^*_\ell$ by the definition
of $\Lambda_{U,u}$, so 
$|\{r_{t,\ell}:\ell < k\}| = k < k(*)+1$ hence for some $m_t \le
k(*)$ we have $\ell < k \Rightarrow r_{t,\ell} \ne m_t$ so for each
$\ell < k$ the sequence $\bar \eta_t \restriction (k(*)+1 \backslash 
\{m_t\})$ is not
from $\{\langle \rho_s:s \le k(*)$ and $s \ne m_t \rangle:
\rho_s \in {}^\omega(U \cup u_\ell)$ for every $s \le k(*)$
such that $s \ne m_t\}$.  

For each $t < t^*$ we define $J(t) = \{m \le k(*):\{\eta_{t,s}:s
\le k(*) \and s \ne m\}$ is included in $U \cup u_\ell$ for no
$\ell \le k\}$.  So $m_t \in J(t) \subseteq \{0,\dotsc,k(*)\}$ and $m
\in J(t) \Rightarrow \bar \eta_t \restriction \{j \le k(*):j \ne m\}
\notin {}^{k(*)+1 \backslash \{m\}}(U \cup u_\ell)$ for 
every $\ell \le k$.  For 
$m \le k(*)$ let $\bar \eta'_{t,m} := \bar \eta_t
\restriction \{j \le k(*):j \ne m\}$ and $\bar \eta'_t := \bar \eta'_{t,m_t}$.
Now we can choose $s_t <
\omega$ by induction on $t$ such that
\mr
\item "{$(*)$}"  if $t_1 < t,m \le k(*)$ and $\bar \eta'_{t_1,m} = 
\bar \eta'_{t,m}$, then 
\nl
$\eta_{t,m} \restriction s_t \notin \{\eta_{t_1,m}
\restriction \ell:\ell < \omega\}$.
\ermn
Let $Y^* = \{x_{m,\bar \eta} \in G_{U \cup u}:x_{m,\bar \eta}
\notin G_{U \cup u_\ell}$ for $\ell < k\} \cup \{ y_{\bar \eta,n} \in
G_{U \cup u}:y_{\bar \eta,n} \notin G_{U \cup u_\ell}$ for $\ell <
k\}$. 
\nl
Let
\medskip

$\quad \qquad Y_1 = \{x_{m,\bar \eta,\nu} \in Y^*$: for no $t < t^*$ do we have
$m = m_t \and \bar \eta = \bar \eta'_t\}$.

$$
\align
Y_2 = \{x_{m,\bar \eta,\nu} \in Y^*:&\,x_{m,\bar \eta} \notin
Y_1
\text{ but for no} \\
 &\,t < t^* \text{ do we have } m = m_t \and \bar \eta = \bar \eta'_t
\and \\
  &\,\eta_{t,m_t} \restriction s_t \trianglelefteq \nu \triangleleft
\eta_{t,m_t}\}
\endalign
$$

$\qquad \quad Y_3 = \{y_{\bar \eta,n}:y_{\bar \eta,n} \in Y^*
\text{ and } n \in [s_t,\omega) \text{ for the } t < t^* \text{
such that } \bar \eta = \bar \eta_t\}$. 
\sn
Now the desired conclusion follows from
\mr
\item "{$(*)_1$}"  $\{y + G_{U,u}:y \in Y_1 \cup Y_2 \cup Y_3\}$
generates $G_{U \cup u} /G_{U,u}$
\sn
\item "{$(*)_2$}"  $\{y + G_{U,u}:y \in Y_1 \cup Y_2 \cup Y_3\}$
generates $G_{U \cup u} /G_{U,u}$ freely.
\endroster
\bn
\demo{Proof of $(*)_1$}  It suffices to check that all the generators
of $G_{U \cup u}$ belong to $G'_{U \cup u} =: \langle Y_1 \cup Y_2
\cup Y_3 \cup G_{U,u} \rangle_G$.

First consider $x = x_{m,\bar \eta,\nu}$ where $\eta \in {}^{k(*)+1}(U
\cup u),m < k(*)$ and $\nu \in {}^n S$ for some $n < \omega$.  If $x \notin
Y^*$ then $x \in G_{U,u_\ell}$ for some $\ell < k$ but $G_{U \cup
u_\ell} \subseteq G_{U,u} \subseteq G'_{U \cup u}$ so we are done, hence
assume $x \in Y^*$.  If $x \in Y_1 \cup Y_2 \cup Y_3$ we are done so
assume $x \notin Y_1 \cup Y_2 \cup Y_3$.  As $x \notin Y_1$ for some
$t < t^*$ we have $m = m_t \and \bar \eta = \eta'_t$.  As $x \notin
Y_2$, clearly for some $t$ as above we have 
$\eta_{t,m_t} \restriction s_t \trianglelefteq
\nu \triangleleft \eta_{t,m_t}$.  Hence by Definition \scite{6.1} the
equation $\boxtimes_{{\bar \eta}_t,n}$ from Definition \scite{6.1}
holds, now $y_{{\bar \eta}_t,n},y_{{\bar \eta}_t,n+1} \in 
G'_{U \cup u}$.  So in order to deduce from the equation that 
$x = x_{\bar \eta'_t \upharpoonleft <m_t,n>}$ belongs to $G_{U \cup u}$, it
suffices to show that $x_{\bar \eta'_{t,j} \upharpoonleft <j,n>} 
\in G'_{U \cup u}$ for each $j \le k(*),j \ne m_t$.  But each such
$x_{\bar \eta'_{t,j} \upharpoonleft <j,n>}$ belong to $G'_{U \cup u}$ 
as it belongs to $Y_1 \cup Y_2$. \nl
[Why?  Otherwise necessarily for some $r < t^*$ we have
$j = m_r,\bar \eta'_{t,j} = \bar \eta'_{r,m_r}$ 
and $\eta_{r,m_r} \restriction s_r \trianglelefteq \eta_t
\restriction n \triangleleft \eta_{r,m_r}$ so $n \ge s_r$ and as said
above $n \ge s_t$.  Clearly $r \ne t$ as $m_r
= j \ne m_t$, now as $\bar \eta'_{t,m_r} = \bar
\eta'_{r,m_r}$ and $\bar \eta_t \ne \bar \eta_r$ (as $t \ne r$) 
clearly $\eta_{t,m_r}
\ne \eta_{r,m_r}$.  Also $\neg(r < t)$ by $(*)$ above applied with
$r,t$ here standing for $t_1,t$ there  as $\eta_{r,m_r}
\restriction s_r \trianglelefteq \eta_{t,j} \restriction n \triangleleft
\eta_{r,m_r}$.  Lastly for if $t < r$, again $(*)$ applied with $t,r$
here standing for $t_1,t$ there as $n \ge m_t$ gives
contradiction.]
\nl
So indeed $x \in G'_{U \cup u}$.

Second consider $y = y_{\bar \eta,n} \in G_{U \cup u}$, if $y \notin
Y^*$ then $y \in G_{U,u} \subseteq G'_{U \cup u}$, so assume $y \in
Y^*$.  If $y \in Y_3$ we are done, so assume $y \notin Y_3$, so for
some $t,\bar \eta = \bar \eta_t$ and $n < s_t$.  We prove by downward
induction on $s \le s_t$ that $y_{\bar \eta,s} \in G'_{U \cup u}$, this
clearly suffices.  For $s=s_t$ we have $y_{\bar \eta,s} \in Y_3
\subseteq G'_{U \cup u}$; and if $y_{\bar \eta,s+1} \in G'_{U \cup u}$
use the equation $\boxtimes_{\bar \eta_t,s}$ from \scite{6.1}, in the 
equation
$y_{\bar \eta,s+1} \in G'_{U \cup u}$ and the $x$'s appearing in the
equation belong to $G'_{U \cup u}$ by the earlier part of the proof
(of $(*)_1$) so
necessarily $y_{\bar \eta,s} \in G'_{U \cup u}$, so we are done.
\enddemo
\bigskip

\demo{Proof of $(*)_2$}  We rewrite the equations in the new variables
recalling that $G_{U \cup u}$ is generated by the relevant variables
freely except the equations of $\boxtimes_{\bar \eta,n}$ from
Definition \scite{6.1}.  After rewriting, all the equations disappear.
\enddemo
\bn
\ub{Case 2}:  $U$ is uncountable.

As $\aleph_1 \le |U| \le \aleph_{k(*)-k}$, necessarily $k < k(*)$.

Let $U = \{\rho_\alpha:\alpha < \mu\}$ where $\mu = |U|$, 
list $U$ with no repetitions.
Now for each $\alpha \le |U|$ let $U_\alpha := 
\{\rho_\beta:\beta < \alpha\}$ and if $\alpha < |{\Cal U}|$ then 
$u_\alpha = u \cup \{\rho_\alpha\}$.  Now  
\mr
\item "{$\odot_1$}"  $\langle(G_{U,u} + G_{U_\alpha \cup u})/G_{U,u}:
\alpha < |U| \rangle$ is an 
increasing continuous sequence of subgroups of $G_{U \cup u}/G_{U,u}$.
\nl
[Why?  By \scite{af.42}(6).]
\sn
\item "{$\odot_2$}"  $G_{U,u} + G_{U_0 \cup u}/G_{U,u}$ is free.
\nl
[Why?  This is $(G_{U,u} + G_{\emptyset \cup u})/G_{U,u} = (G_{U,u} +
G_u)/G_{U,u}$ which by \scite{af.42}(8) is isomorphic to
$G_u/G_{\emptyset,u}$ which is free by Case 1.]
\ermn
Hence it suffices to prove that for each
$\alpha < |U|$ the group $(G_{U,u} + G_{U_{\alpha +1} \cup
u})/(G_{U,u} + G_{U_\alpha \cup u})$ is free.  But easily
\mr
\item "{$\odot_3$}"    this group is isomorphic to $G_{U_\alpha \cup
u_\alpha}/G_{U_\alpha,u_\alpha}$.
\nl
[Why?  By \scite{af.42}(7) with $U_\alpha,U_{\alpha +1},U,\rho_\alpha,u$
here standing for $U_1,U_2,U,\eta,u$ there.]
\sn
\item "{$\odot_4$}"  $G_{U_\alpha \cup u_\alpha}/G_{U_\alpha,u_\alpha}$
is free. 
\nl
[Why?  By the induction hypothesis, as $\aleph_0 + |U_\alpha| < |U|
\le \aleph_{k(*)-(k+1)}$ and $|u_\alpha| = k+1 \le k(*)$.] 
\ermn
2) If $k(*) = 0$ just use \scite{af.28}, so assume $k(*) \ge 1$.  Now
the  proof is similar to (but easier than) the proof of case (2)
inside the proof of part (1) above.
\nl
${{}}$  \hfill$\square_{\scite{af.56}}$
\newpage


\nocite{ignore-this-bibtex-warning} 
\newpage
    
REFERENCES.  
\bibliographystyle{lit-plain}
\bibliography{lista,listb,listx,listf,liste}

\enddocument